\documentclass[12pt,a4paper,openany]{report}
\usepackage[french]{babel}
\usepackage{amsxtra}
\usepackage{amssymb}
\usepackage{amscd}
\usepackage{makeidx}
\usepackage{amsthm}

\catcode `à = \active \def à{\`a}
\catcode `â = \active \def â{\^a}
\catcode `é = \active \def é{\'e}
\catcode `è = \active \def è{\`e}
\catcode `ê = \active \def ê{\^e}
\catcode `ë = \active \def ë{\"e}
\catcode `ä = \active \def ä{\"a}
\catcode `î = \active \def î{\^\i}
\catcode `ï = \active \def ï{\"\i}
\catcode `ô = \active \def ô{\^o}
\catcode `ö = \active \def ö{\"o}
\catcode `ù = \active \def ù{\`u}
\catcode `û = \active \def û{\^u}
\catcode `ü = \active \def ü{\"u}
\catcode `ç = \active \def ç{\c c}
\def\cin{ C^{\infty}}
\def\cov{\bigtriangledown}

\def\dbaret{\bar{\partial^*}}
\def\tr{{\bf{tr}}}
\def\pr{{\bf{P}}}
\def\dc{d^{c}}
\def\dbaretid#1{\bar{\partial}^{*}_{#1}}
\def\ddet{\bar{\partial}\bar{\partial^*}}
\def\ddetid#1{\bar{\partial}\bar{\partial}^{*}_{#1}}

\def\detdid#1{\bar{\partial}^{*}_{#1}\bar{\partial}}
\def\laplac#1#2{\Delta_{{#1},{#2}}}

\def\explapid#1#2{\exp(-t\Delta_{{#1},{#2}})}

\def\explapidu2#1#2{\exp(\frac{-u}{2}\Delta_{{#1},{#2}})}

\def\explap#1#2{\exp(-t\Delta_{{#1},{#2}})P_{n}^\perp}
\def\explapM#1#2{\exp(-M\Delta_{{#1},{#2}})P_{n}^\perp}
\def\du{\frac{d}{du}}
\def\dt{\frac{d}{dt}}
\def\ds{\frac{d}{ds}}

\def\espco{{\cal A}}

\def\im{$im$}

\def\jaugec{{\cal G}^{\mathbb {C}}}
\newtheorem{thm}{Th\'eor\`eme  }[section]
\newtheorem*{thm*}{Th\'eor\`eme  }

\newtheorem{dfn}[thm]{D\'efinition  }
\newtheorem{cor}[thm]{Corollaire  }

\newtheorem*{cor*}{Corollaire  }
\newtheorem{prop}[thm]{Proposition  }
\newtheorem*{prop*}{Proposition  }
\newtheorem{lem}[thm]{Lemme  }
\newenvironment {remq}{\noindent {\bf Remarque .}}{}

\newenvironment {demo}{\noindent {\bf D\'emonstration . }}{
$\square $\\}{}

\def\db{\bar{\partial}}
\def\lr{\longrightarrow}
\def\db{\bar{\partial}}
\newcommand{\norm}[1]{\Vert #1\Vert_{L^{2}_{1}}}
\newcommand{\nor}[1]{\Vert #1\Vert}
\newcommand{\mo}[1]{\vert #1\vert}

\newcommand{\norml}[1]{\Vert #1\Vert_{L^{2}}}

\newcommand{\normca}[1]{\Vert #1\Vert_{C^{0}(A(r',R))}}

\newcommand{\normlh}[1]{\Vert #1\Vert_{ L^{2}(h)}}

\newcommand{\normhsi}[1]{\vert #1\vert_{h}}
\newcommand{\normhox}[1]{\Vert #1\Vert_{h_{\ox(n)}}}
\newcommand{\hox}{{h_{\ox(n)}}}
\newcommand{\normk}[1]{\vert #1\vert_{k}}
\newcommand{\negli}[1]{\frac{1}{n^{#1}}}

\def\C{\mathbb C}
\def\P{\mathbb P}
\def\R{\mathbb R}
\def\N{\mathbb N}
\def\Z{\mathbb Z}
\def\ox{\mathcal{O}_{X}}

\def\lr{\longrightarrow}

\def\rg{\mbox{ rg }}
\def\im{\mbox{im}}

\makeindex
\begin{document}
\title{Approximation de m\'etriques de Yang-Mills pour un fibr\'e $E$ à partir de
m\'etriques induites de $H^{0}(X,E(n))$}
\author{Cécile Drouet\\
U.M.R. 5580 Laboratoire Emile Picard\\ Université
Paul Sabatier\\ U.F.R. M.I.G.\\
 118, route de Narbonne\\31062 Toulouse Cedex 4 \\
Mél : drouet@picard.ups-tlse.fr\\
 Tél  : (+33) 05 61 55 60 50 \\
 Fax : (+33) 05 61 55 82 00}
\date{\today}
%%%%%%%%%%%%%%%%%%%%%%%%%%%%%%%%%%%%%%%%%%%%%%%%%%%%%%%%%%%%%%%%%%%%%
\maketitle
%%%%%%%%%%%%%%%%%%%%%%%%%%%%%%%%%%%%%%%%%%%%%%%%%%%%%%%%%%%%%%%%%%%%%
\tableofcontents
%%%%%%%%%%%%%%%%%%%%%%%%%%%%%%%%%%%%%%%%%%%%%%%%%%%%%%%%%%%%%%%%%%%%%%%
%%%%%%%%%%%%%%%%%%%%%%%%%%%%%%%%%%%%%%%%%%%%%%%%%%%%%%%%%%%%%%%%%%%%%%%%
%%%%%%%%%%%%%%%%%%%%%%%%%%%%%%%%%%%%%%%%%%%%%%%%%%%%%%%%%%%%%%%%%%%%%%%
%%%%%%%%%%%%%%%%%%%%%%%%%%%%%%%%%%%%%%%%%%%%%%%%%%%%%%%%%%%%%%%%%%%%%%%%
%%%%%%%%%%%%%%%%%%%%%%%%%%%%%%%%%%%%%%%%%%%%%%%%%%%%%%%%%%%%%%%%%%%%%%%
%%%%%%%%%%%%%%%%%%%%%%%%%%%%%%%%%%%%%%%%%%%%%%%%%%%%%%%%%%%%%%%%%%%%%%%%
%%%%%%%%%%%%%%%%%%%%%%%%%%%%%%%%%%%%%%%%%%%%%%%%%%%%%%%%%introdu2.tex}
%%%%%%%%%%%%%%%%%%%%%%%%%%%%%%%%%%%%%%%%%%%%%%%%%%%%%%%%%%%%%%%%%%%%%%%
%%%%%%%%%%%%%%%%%%%%%%%%%%%%%%%%%%%%%%%%%%%%%%%%%%%%%%%%%%%%%%%%%%%%%%%
\chapter*{Introduction}\addcontentsline{toc}{section}{Introduction}
  Dans ce travail, nous considérons un fibré
vectoriel holomorphe, noté $E_{0}$ sur une courbe
algébrique $X$, de genre g, plongée dans $\C
P^{l}$.
 Le rang  de $E_{0}$ est $2$ et son degré est quelconque noté $d$.
  Soit  $E_{n}$ le fibré
$E_{0}\otimes
\ox(n)$.On munit l'espace $H^{0}(X,\ox(1))$ d'une métrique hermitienne.
 Le fibré $\ox(1)$ est muni de la métrique induite et notée $h_{\ox(1)}$.
  La forme volume de $X$, notée
 $\omega$ sera définie à partir de la forme de Chern de cette métrique et vérifiera
 $\int_{X}\omega=1$.\label{volume}

Nous pouvons munir $E_{n}$ de métriques
compatibles avec sa structure complexe et nous
notons $Met(E_{n})$ l'ensemble de ces métriques.
Remarquons que nous avons l'isomorphisme suivant :
  \begin{eqnarray*}\label{indep}
Met(E_{0})&\stackrel{\simeq}{\lr}& Met(E_{n})\\
h_{0}&\longmapsto& h_{0}\otimes
(h_{\ox(1)})^{\otimes n}=:h.
  \end{eqnarray*}
L'espace $Met(E_{n})$ ne dépend pas de $n$ via cet
isomorphisme.\\
 On  considère aussi l'espace  des
métriques que l'on pourra mettre sur
$W_n=H^{0}(X,E_{n})$ noté $Met(W_n)$.\\ On
remarque, par le théorème de Riemann-Roch, que la
dimension de l'espace $Met(W_n)$ augmente avec $n$
alors que l'espace $Met(E_{n})$ (indépendant de
$n$) est de dimension infinie. Cela nous incite à
étudier  les liens éventuels entre les deux
espaces  $Met(W_n)$ et $Met(E_n)$, quand $n$
devient grand .\\
 De plus,
les théorèmes de M. S. Narasimhan et C. S.
Seshadri \cite{nara} et de S.K. Donaldson
\cite{don1} donnent deux critères de stabilité du
fibré $E_{n}$ : l'un est algébrique et l'autre est
analytique faisant intervenir la courbure de Chern
du fibré $(E_{n},h)$, $h\in Met(E_{n})$.
 C'est pourquoi, on peut se poser la question
d'une caractérisation plus algébrique d'objets qui
pour l'instant sont définis de façon analytique
comme, par exemple, les métriques de Yang-Mills.
Ces métriques correspondent aux minima de la
fonctionnelle suivante :$${\mathcal
{YM}}(D):=(R(D), R(D))
=\int_{M} <R(D), R(D)>*\,(1),$$
où $D$ est la connexion de Chern du fibré
holomorphe hermitien $(E_{n},h)$ avec $h\in
Met(E_{n})$ et $R(D)$ la courbure de $D$ (cf.
paragraphe \ref{ymills}).\\
 On commence tout d'abord par relier les
espaces $Met(W_n)$ et $Met(E_n)$. Pour cela on
définit
 deux morphismes (cf. p. \pageref{l} et p. \pageref{i}).
 Le premier est
:
$$L_{n}\, :\, Met(E_{n})\lr Met(W_n),$$
qui à partir d'une métrique $h\in Met(E_n)$
fournit sa métrique $L^{2}$ dans $Met(W_n)$. Le
second morphisme est donné par :
\begin{eqnarray*}
I_{n} : Met (W_n)&\lr& Met(E_{n})\\ m&\longmapsto&
I_{n}(m)(e_{x},e_{x})= min\{m(v,v)/v\in W_n\mbox{
et }v(x)=e_{x}\},
\end{eqnarray*}
avec $e_{x}\in E_{x}$. L'application $I_{n}$
définit donc, à partir d'une métrique $m\in
Met(W_{n})$, une ``métrique induite" via les
surjections que nous avons entre $W_{n}$ et le
fibré $E_{n}$ lorsque $E_{n}$ est engendré par ses
sections (ce qui est le cas pour $n$ assez
grand).\\ Il est facile de voir que pour $n$ assez
grand, pour $h\in Met(E_{n})$, $I_{n}I_{n}(h)$ est
approximativement égal à $h$ à une constante
multiplicative près. On peut noter l'analogie avec
le résultat de G. Tian \cite{tian} pour les
métriques kälhériennes, et, au fond, nous obtenons
ainsi une généralisation directe du théorème de
Stone-Weierstrass. Notre approche de ce résultat,
via des sections concentrées, est similaire à
celle de G. Tian. Mais nous voulons approximer la
métrique de Yang-Mills par des métriques
$I_{n}(k)$ pour des points $k$ déterminés de façon
naturelle.
\\
Pour cela nous utiliserons un travail de S.K.
Donaldson \cite{don2} concernant la construction
de métriques de Quillen sur l'espace $Met(E_{n})$
donnant l'égalité liant $Met(E_{n})$ et $Met(W_n)$
:
$$Norm_{h, \omega}(\beta)=(\det_{W_n}(L_{n}(h)))^{{\frac{1}{2}}}\times T_{h}(X, E_{n}),$$
où $\beta$ est un élément fixé d'un fibré
particulier $L$ sur $Met(E_{n})$ et $T_{h}(X,
E_{n})$ est  la torsion analytique du fibré
$(E_{n},h)$ (cf. p. \pageref{norma})
. Les extrémums de
$Norm_{h,\omega}(\beta)$ correspondent aux
métriques de Yang-Mills  et  ses variations nous
fournissent la fonctionnelle de S.K. Donaldson,
sur $Met(E_{n})$ notée $\mathcal{M}$. En
parallèle, en s'inspirant de la caractérisation de
la stabilité du fibré $E_{n}$ décrite par G. Kempf
et L. Ness dans \cite{kempf} et dans \cite{git},
nous introduisons  la fonctionnelle
$\mathcal{KN}_n$ sur $Met(W_n)$ définie par (cf.
p. \pageref{kkkk}) :
\begin{eqnarray*}
\mathcal {KN}_{n} : Met(W_{n})&\longrightarrow&  {\mathbb R}\\
 m&\longmapsto&\frac{1}{2}ln(\det_{W_{n}}m)+\frac{\chi(n)}{2r} \int_{X}\ln(\det_{E_{n}}\circ I_{n}(m))\omega,
\end{eqnarray*}
où $\chi(n)$ est la caractéristique d'
Euler-Poincaré de $E_{n}$.
 Notre but
sera d'étudier le comportement  des fonctionnelles
$\mathcal {KN}_{n}$ et $\mathcal {M}$. Dans cette
optique, nous avons été amenés à étudier les
variations de la torsion analytique $\tau_h(X,
E_{n})=2\, ln
\, T_h(X, E_{n})$ (cf. p.
\pageref{toranalyt}) et cela de manière intrinsèque. L'étude de $T_{h}(X, E_{n})$ en fonction du degré
d'un fibré est d'ailleurs abordé dans \cite
{Bisvas} mais nous en généralisons l'étude à des
variations de métriques hermitiennes sur le fibré
holomorphe $E_{n}$. On peut remarquer que les
résultats obtenus sont  quelque peu différents de
ceux obtenus dans \cite{Bisvas} et présentent un
intérêt en soi. En reprenant les définitions du
paragraphe \ref{definit}, on a
 :
\begin{thm*}(cf. p. \pageref{fini})
 Etant
données deux m\'etriques hermitiennes $h_{0}$ et
$h_{1}$ sur $E_{n}$, on consid\`ere le chemin
$h(u)=uh_{1}+(1-u)h_{0}$ avec $u
\in
\lbrack0, 1\rbrack$, alors
\begin{eqnarray*}
{\du}\tau_{h(u)}(X, E_{n}) &=& O(\negli{})
\end{eqnarray*}
\end{thm*}
 Pour cette étude nous construisons ce que nous appelons des sections ``concentrées"
sur le fibré hermitien $(E_{n},h)$
. Pour cela,
étant donné $z_{0}\in X$, on choisit
convenablement (cf. chap 8, p. \pageref{baser})
une base locale holomorphe $\{e_{1},e_{2}\}$ de
$E_{n}$ dans un voisinage de  $z_{0}$ où $z_{0}\in
X$. Une section ``concentrée" est  alors une
section holomorphe, projection orthogonale pour
$h\in Met(E_{n})$ sur $H^{0}(X,E_{n})$ de la
section généralisée $u_{1}$ définie par :

$$u_{1}=\delta_{z_{0}}e_{1}\mbox{ (cf. p.\pageref{distribu}) }.$$

Cet outil nous fournira également les résultats
suivants :

\begin{enumerate}
  \item le calcul du noyau de l'opérateur  correspondant à la
   projection orthogonale pour la métrique $h$ des sections de $E_{n}$
    sur l'espace $W_n$ (cf. p. \pageref{cho} et p. \pageref{che}),

  \item l'image de $h\in Met(E_{n})$ par $I_{n}L_{n}$:

  \begin{thm*}(cf. p. \pageref{end}) Soit $(E_{n},h)$ un fibré holomorphe hermitien sur $X$, on a :
$$h^{-1}I_{n}L_{n}(h)=
\frac{\pi}{n^{2}}(R_{\C}(\ox(n)Id_{E_{0}}
-R_{\C,h}(E_{0})Id_{\ox(n)}+r Id_{E_{n}})+{\bf{O}(\negli{3})}.$$ Pour les notations on
regardera p. \pageref{courbb}.
\end{thm*}

\end{enumerate}
Ce  théorème nous permet alors d'établir le lien
entre les deux fonctionnelles $\mathcal {M}$ et
$\mathcal {KN}_n$ :
\begin{thm*}(cf. p. \pageref{princ})
$\exists C\,,\quad \forall h,k\in
Met(E_{n},C^{4})$
 $$\vert (\mathcal{M}-\mathcal{KN}_n\circ L_{n} )(h)-(\mathcal{M}-\mathcal{KN}_n\circ L_{n} )(k)\vert\leq \frac{C}{n}.$$
\end{thm*}

Soit maintenant $U$ un compact de $Met(E_{n},C
^{4})$ donné. On démontre également que  pour $n$
assez grand, si $m$ est un extrémum de
${\mathcal{KN}_{n}}{\vert}_{ L_{n}(U)}$ alors $ m$
est l'image d'un minimum local de $\mathcal M$. En
particulier, si le minimum de Yang-Mills est
contenu dans  $U$ alors le minimum de $\mathcal
{KN}_{n}{\vert}_{ L_{n}(U)}$ pour $n $ assez grand
correspond au minimum de Yang-Mills. Plus
précisément, on a le théorème suivant
:

\begin{thm*}(cf. p. \pageref{top})

 Soit  $E_{n}$ un fibré stable sur $X$. Soit
$U$ un compact convexe donné de $Met(E,C^{4})$.
Supposons que  $h_{YM}$, la métrique de
Yang-Mills(elle existe car $E_{n}$ est stable),
appartienne à $U$. On a alors
:
$$\forall \varepsilon >0 \quad \exists n_{U},\quad \forall n\geq n_{U},$$
il existe  m  extrémum de
${\mathcal{KN}_{n}}{\vert}_{ L_{n}(U)}$ et
$$\nor{h_{YM}- c.I_{n}(m)}_{L^{2}_{1}(Met(E))}\leq \varepsilon, $$
avec $c=\frac{n}{\pi}(1+O(1)).$

\end{thm*}
Notre projet de recherche consiste à généraliser
ce résultat en utilisant  non pas un minimum de
$\mathcal{KN}_{n}$ sur $L_{n}(U)$ mais le minimum
de $\mathcal{KN}_{n}$ sur tout $W_{n}$
. On peut songer à une écriture algorithmique de
ce processus. En effet si on se donne un nombre
suffisant de points sur $X$, on peut alors
rechercher le  minimum de $\mathcal{KN} $, via le
calcul de déterminants, et en déduire qu'il
provient approximativement d'un minimum de
Yang-Mills sous de bonnes conditions.
%%%%%%%%%%%%%%%%%%%%%%%%%%%%%%%%%%%%%%%%%%%%%%%%%%%%%%%%%%%%%%%%%%%%%%%
%%%%%%%%%%%%%%%%%%%%%%%%%%%%%%%%%%%%%%%%%%%%%%%%%%%%%%%%%%%%%%%%%%%%%%%%
%%%%%%%%%%%%%%%%%%%%%%%%%%%%%%%%%%%%%%%%%%%%%%%%%%%%%%%%%%%%%%%%%%%%%%%
%%%%%%%%%%%%%%%%%%%%%%%%%%%%%%%%%%%%%%%%%%%%%%%%%%%%%%%%%%%%%%%%%%%%%%%%
%%%%%%%%%%%%%%%%%%%%%%%%%%%%%%%%%%%%%%%%%%%%%%%%%%%%%%%%%%%%%%%%%%%%%%%
%%%%%%%%%%%%%%%%%%%%%%%%%%%%%%%%%%%%%%%%%%%%%%%%%%%%%%%%%%%%%%%%%%%%%%%%

%%%%%%%%%%%%%%%%%%%%%%%%%%%%%%%%%%%%%%%%%%%%%%%%%%%%%%%%%%%%%%%%%%%%%
%%%%%%%%%%%%%%%%%%%%%%%%%%%%%%%%%%%%%%%%%%%%%%%%%%%%%%%%%%%%%%%%%%%%%%
%%%%%%%%%%%%%%%%%%%%%%%%%%%%%%%%%%%%%%%%%%%%%%%%%%%%%%%%%%%%%part1h.tex}
%%%%%%%%%%%%%%%%%%%%%%%%%%%%%%%%%%%%%%%%%%%%%%%%%%%%%%%%%%%%%%%%%%%%%%%%%%

\part{Stabilité
 : deux approches}

\section*{Introduction}\addcontentsline{toc}{section}{Introduction}
Dans la suite du texte et sauf mention contraire
$X$\index{Notation!$X$} d\'esigne une courbe
alg\'ebrique lisse compacte projective sur
${\mathbb {C}}$ ou de façon \'equivalente un
surface de Riemann compacte lisse. Dans ce cas $X$
pourra \^etre consid\'er\'ee comme une vari\'et\'e
kälh\'erienne et on notera
$\omega$\index{Notation!$\omega$} la forme
kälh\'erienne associ\'ee
\`a la structure riemannienne $g$. La variété vue comme variété k\" alh\'erienne
sera notée $(X,\omega)$ et vue comme variété
riemannienne, $(X,g)$.\\ On s'intéresse à des
fibrés vectoriels sur $X$. Ils seront considérés
holomorphes si cela n'est pas précisé. On notera
$\deg\, E$ le degré  et  $\rg\, E$ son rang  . On
prendra en général $\deg\, E=d$ et on se limitera
au cas où $\rg\, E=2$. La dimension de
$H^{0}(X,E)$ sera notée $p$.
\\ Cette partie est consacrée à deux approches
différentes de la stabilité. \\ Le premier
chapitre traite de l'aspect algébrique de la
stabilité. Après avoir donné quelques unes des
propriétés des fibrés stables,  nous rappelons les
résultats de M. S. Narashiman et C. S. Seshadri.
 Ensuite, nous nous
intéressons à la caractérisation de la stabilité,
via un produit de grassmanniennes, préambule à la
construction dans la deuxième partie de la
fonctionnelle de G. Kempf et L. Ness.
 Dans le second
chapitre, le côté géométrie analytique est
développé. Après quelques définitions sur la
courbure, quelques généralités sur l'espace des
connexions et la définition de la fonctionnelle de
Yang-Mills, nous finissons par la relecture d'une
caractérisation des fibrés stables par S. K.
Donaldson dans \cite{don1}.  Cela nous conduit au
théorème suivant :
\begin{thm*}(cf. p. \pageref{donalds})
Un fibr\'e holomorphe, indécomposable, $E$ sur $X$
est stable si et seulement il y a sur $E$ une
connexion unitaire tel que sa courbure centrale
satisfasse
$$*F=-2\pi i \mu(E)\mbox{ avec } \mu(E)=\frac{\deg\, E}{rg\, E}. $$
Cette connexion est unique \`a isomorphisme
pr\`es.
\end{thm*}
 Nous donnons une esquisse de la
démonstration de ce théorème qui fait intervenir
une fonctionnelle sur l'espace des connexions du
fibré $E$, notée ici $J$  qui sert de base à la
fonctionnelle de S.K. Donaldson définie dans la
partie deux (cf. paragraphe \ref{cestfini}) et
notée $\mathcal M$.
\chapter{Fibr\'es stables et espace de module }

\section{Fibr\'es stables }\label{plongi}

Nous notons par la suite les fibr\'es vectoriels
par des lettres bâtons et les faisceaux  par des
lettres curvilignes.  Nous rappelons
\'egalement que nous attachons \`a la courbe $X$ le
faisceau ${\cal O}_{X}$\index{Notation!${\cal
O}_{X}$}.  Apr\`es avoir d\'efini un plongement,
noté $\theta$, dans ${\mathbb P}^{l}$, pour un $l$
donn\'e ($\mathbb {P}^{l}$  étant muni si
nécessaire de sa métrique de Fubini-Study (cf.
\cite{grif}), nous d\'efinissons les faisceaux
${\cal O}_{X}(d)$ et ce quelque soit $d$
appartenant à ${\mathbb Z}$.
 Soit  $X$ une vari\'et\'e complexe de dimension complexe un i.e une surface de
Riemann  compacte.  Soit E un fibr\'e holomorphe
de rang $ r $ sur X.  On pose :
$$\deg \,E =\deg {\cal E}=c_1(E), $$\index{Degré} \index{Classe de Chern} avec  $c_1(E)$  la premi\`ere classe de
 Chern de $E$. \\

\begin{dfn}

On d\'efinit la pente du fibr\'e $E$ par : $$ \mu (E)=c_1(E)/
 \rg(E) .$$
\end{dfn}\index{Notation! $\mu$}

\begin{dfn}\index{Fibré!stable}\index{Fibré!semi-stable}
 E est stable (resp.  semi-stable ) si pour tout sous fibr\'e E'  propre de E
avec $0<\rg E'<\rg E$  on a :
$$\mu (E')<\mu (E)  $$  (resp.  $\mu (E')\leq\mu (E)  $).
\end{dfn}
Les propri\'et\'es usuelles sont regroup\'ees dans la proposition
suivante :
\begin{prop}\label{stabil}
\begin{enumerate}
\item Tout fibr\'e en droite est stable.
\item Si F est stable et L est un fibr\'e en droite alors $F\otimes L$ est stable.
\item Si $F_{1}$ et $F_{2}$ sont stables et $\mu(F_{1})=\mu(F_{2})$ alors tout
homomorphisme non nul $h : F_{1}\longrightarrow F_{2}$ est un
isomorphisme.
\item Tout fibr\'e vectoriel stable est simple .
\item Soit F un fibr\'e semi-stable sur X de rang r et de degr\'e d.  Supposons que
$d\geq r(2g-1) $, alors
\begin{enumerate}
\item $H^{1}(F)=0$,
\item F est engendr\'e par ses sections.

 \end{enumerate}
  \end{enumerate}
 \end{prop}
 \begin{demo}
 \begin{enumerate}
 \item est \'evidente
 \item Comme $$\deg(F\otimes L)=\deg(F). \rg(L)+\deg(L). \rg(F)\,\, ,$$
 on obtient le r\'esultat.
 \item Consid\'erons un morphisme de faisceaux :
 $$h : {\cal F}_{1}\longrightarrow {\cal F}_{2}$$
  et soit les sous fibr\'es  $G_{1}$ et $G_{2}$ associ\'es respectivement \`a $\ker\,h$ et $\im\, h$
respectivement.  On a alors :
\begin{eqnarray*}
\rg(F_{1})&=&\rg(G_{1})+\rg(G_{2})\\
\deg(F_{1})&=&\deg(\ker\,h)+\deg(\im \,h)\\
 &\leq&\deg(G_{1})+\deg(G_{2})\\
&\leq& \rg(G_{1}). \mu(F_{1})+\rg(G_{2}).
\mu(F_{2})\\ &=&(\rg(G_{1})+\rg(G_{2})).
\mu(F_{1})=\deg(F_{1}),
\end{eqnarray*}
On en déduit que $\ker\,h=0$ et $\im\, h=F_{2}$
 \item
 Supposons que F soit stable et consid\'erons un endomorphisme
 $$h : F \longrightarrow F . $$
 Prenons un point $x\in X$.  On a alors une application lin\'eaire :
 $$h_{x} : F_{x}\longrightarrow F_{x}. $$
 Soit $\lambda$ une valeur propre  de $h_{x}$.
 Alors $h-\lambda Id_{F}$ n'est pas un isomorphime.  Donc d'apr\`es (3) $h-\lambda Id_{F}=0$,
 c'est \`a dire $h=\lambda Id_{F}.$
 \item
\begin{enumerate}
\item Supposons que $H^{1}(X,\mathcal{F})\not= \{0\}$ alors par dualit\'e il existe un morphisme de faisceau
non nul :
$$h : {\cal F}\longrightarrow \omega_{X},$$
où $\omega_{X}$ est le fibré en droite canonique
de X. Le sous fibr\'e G de F engendr\'e par $ker\,
h$ a pour rang $r-1$ et
$$\deg(G)\geq \deg(ker \, h)\geq \deg({\cal F})-\deg(\omega_{X})=d-(2g-2).$$
Puisque $F$ est semi-stable on a :
$$r(d-(2g-2))\leq r. \deg(G)\leq (r-1)d \,,$$
donc $$d\leq r(2g-2)$$ d'o\`u la contradiction.
\item Soit $x\in X$ et consid\'erons la suite exacte :
\begin{equation}\label {ideaux}
0\longrightarrow  {\mathfrak {m}}_{x}. {\cal F}\longrightarrow
{\mathcal F}
\longrightarrow{\mathcal F}_{x}\longrightarrow 0
\end{equation}
o\`u ${\mathfrak {m}}_{x}$ est le faisceau des id\'eaux d\'efinis
au point $x$, ${\mathcal F}_{x}$ le faisceau de torsion avec pour
support $x$.  Nous voulons d\'emontrer que l'application induite :
$$H^{0}(X,{\cal F})\longrightarrow H^{0}(X,{\cal F}_{x})$$
est toujours surjective.  D'apr\`es la suite
exacte longue associ\'ee à la suite \ref{ideaux},
il suffit de montrer que
$$H^{1}(X,{\mathfrak {m}}_{x}. {\mathcal F})=\{0\}.$$  Comme ${\mathfrak {m}}_{x}$
est sans torsion il est localement libre.  On peut consid\'erer le
fibr\'e en droite associ\'e $L_{x}$.  Alors, d'apr\`es (2),
$F\otimes L_{x}$ est semistable.  De plus \'etant donn\'e la suite
exacte suivante
:

$$0\longrightarrow{\mathfrak {m}}_{x} \longrightarrow {\mathcal {O}}_{X}
\lr \Bbbk
\longrightarrow 0$$
on a : $\deg(L_{x})=\deg({\mathfrak {m}}_{x})=-1$.  Donc :
$$\deg(F\otimes L_{x})=d-r>r(2g-2). $$
D'o\`u d'apr\`es ce qui pr\'ec\`ede :
$$H^{1}({\mathfrak {m}}_{x}. {\cal F})=H^{1}(F\otimes L_{x})=0. $$
 \end{enumerate}
 \end{enumerate}
 \end{demo}
 \section{R\'esultats de M. S. Narasimhan et C. S. Seshadri }
 On rapelle brièvement les résultats exposés dans \cite{nara} qui fournit
 un critère de stabilité.
On consid\`ere maintenant une surface de Riemann
compacte $X$ de genre $g\geq 2$\index{Genre}. Soit
$n\in {\mathbb N}$ et $q\in {\mathbb Z}$ tel que
$-n<q\leq 0$.  Soit $\pi$ un groupe discret
agissant effectivement, proprement et
holomorphiquement sur le disque unit\'e $Y$, de
telle sorte que $Y/\pi=X$.  On suppose
\'egalement que la projection naturelle $$p
: Y\longrightarrow X$$ est ramifi\'ee au dessus
d'un point $x_{0}\in X$ de degr\'e $n$  et non
ramifi\'ee ailleurs. Soient $y_{0}\in
p^{-1}(\{x_{0}\})$ et $\tau=\tau(q)$ le
caract\`ere du groupe d'isotropie $\pi_{y_{0}}$ de
$\pi$  tel que l'entier associ\'e \`a $\tau$ soit
$-q$. On a alors le th\'eor\`eme suivant :

 \begin{thm}\cite{nara}\index{Théorème! de M. S. Narasimhan et C. S. Seshadri}\label{se}
  Un fibr\'e vectoriel holomorphe $ F$ de rang
$n$ et de degr\'e $q$ sur $X$ est stable si et seulement si $F$ est
isomorphe \`a $p_{*}^{\pi}(E)$ o\`u $E$ est un $\pi$-fibr\'e
sp\'ecial sur Y de type $\tau(q)$ tel que la repr\'esentation
: $$\rho: \pi \longrightarrow GL(n,\mathbb C)\, ,$$ \`a laquelle E est associ\'e, soit irr\'eductible et unitaire.
\end{thm}\index{Fibré!stable}
Nous donnons les définitions nécessaires à la
compréhension du théorème \ref{se}.
\begin{dfn}Soit $\rho : \pi \longrightarrow GL(n,\mathbb C)$ une repr\'esentation.  Alors $\pi$ op\`ere sur le fibr\'e trivial
 $Y\times{\mathbb C}^{n}$ par $(y,v) \rightsquigarrow (\gamma y,\rho(\gamma)v)
 \quad y\in Y ,\  v\in {\mathbb C}^{n}
 \, ,\, \gamma \in \pi$.  On appelle $\pi$-fibr\'e le fibr\'e associ\'e \`a la repr\'esentation  $\rho$.
 \end{dfn}
 \begin{dfn}
 Soit $\rho : \pi \longrightarrow GL(n,\mathbb C)$ une repr\'esentation.  Elle est dite de type
 $\tau$ si $\forall \gamma \in \pi_{y_{0}}\,\  \rho(\gamma)=\tau(\gamma)Id$.
 Si de plus $\rho : \pi \longrightarrow U(n)$ elle est dite unitaire de type $\tau$.
 \end{dfn}
 \begin{dfn}
 un $\pi$-fibr\'e E sur Y est un $\pi$-fibr\'e sp\'ecial de type $\tau$ si E est associ\'e \`a une
 repr\'esentation $\rho : \pi \longrightarrow GL(n,\mathbb C)$ de type $\tau$.
 \end{dfn}
 \begin{dfn}
 Soit E un $\pi$-fibr\'e de rang n sur Y et $\cal E$ le faisceau associ\'e.  Alors $\pi$ agit
 sur $\cal E$ et donc sur l'image directe $p_{*}(\cal E)$. On consid\`ere donc le sous faisceau, not\'e
  $p_{*}^{{\pi}}(\cal E)$, de $p_{*}(\cal E)$, form\'e des \'el\'ements invariants sous l'action de $\pi$.
$p_{*}^{{\pi}}(\cal E)$ est un sous faisceau localement libre de
${\cal O}_{X}$-module de rang n.  On note $p_{*}^{{\pi}}( E)$ le
fibr\'e associ\'e.  On a donc un foncteur $p_{*}^{\pi}$ de la
cat\'egorie des $\pi$-fibr\'es sur Y vers la cat\'egorie des
fibr\'es vectoriels sur X.
 \end{dfn}
 Un corollaire de ce théorème est :
 \begin{prop}
 Un fibré vectoriel de degré zéro sur X est stable si et seulement
  si il provient d'une représentation irréductible unitaire
 du groupe fondamental  de X, $\pi_{1}(X)$.
 \end{prop}
 \begin{remq} On étend le théorème précédent pour les fibrés de
  tout degré en remarquant qu'il existe un fibré en droite $L$,
  dépendant de $n$ et de $q$ tel que $-n<\deg(F\otimes L)\leq 0$.
 \end{remq}

\section{Schéma de Hilbert et produit de grassmaniennes }\label{ke}
On pourra consulter \cite{seshadri} et
\cite{potier}.\\ Dans ce paragraphe nous nous
intéressons à la correspondance entre les points
d'un schéma de Hilbert qui correspondent aux
fibrés stables et des points, stables pour
l'action du groupe $SL(p)$, d'un produit de
grassmaniennes
.
\begin{prop}\cite{potier}
Soit ${\cal F}$ un faisceau coh\'erent sur X.  Les
quotients de ${\cal F}$ qui ont un polyn\^ome de
Hilbert fix\'e $P$ sont param\'etris\'es par les
points d'un sch\'ema projectif de type fini not\'e
$Q({\cal F}/P)$.  De plus il existe un faisceau
coh\'erent ${\cal U}$ sur $Q({\cal F}/P)\times X$
et un morphisme surjectif,   $p_{2}^{*}\, :\,
{\cal F}\longrightarrow {\cal U}$, tel que ${\cal
U}$ est plat sur $Q({\cal F}/P)$ et ``universel".
\end{prop}
Dans notre cas on prend $p\geq 1$:
$${\cal F}={\cal O}_{X}^{p} \mbox{ et } P(m)=p+rm\,. $$
\begin{dfn}\cite{potier}\label{rmod} On note ${\cal U}_{q},\  q\in Q({\cal F}/P)$, le faisceau coh\'erent sur X induit par $\cal U$ restreint \`a $
\{q\}\times X$.  On consid\`ere les points q de $Q({\cal F}/P)$ tel que :
\begin{enumerate}
\item ${\cal U}_{q}$ soit localement libre.
\item $H^{0}({\cal F})\longrightarrow H^{0}({\cal U}_{q})$ est un isomorphisme
 $\forall q \in  Q({\cal F}/P)$

\end {enumerate}
On notera ${\mathcal {R}}$ cet ensemble.
\end{dfn}
\begin{remq} Le choix du polyn\^ome de Hilbert P et la d\'efinition pr\'ec\'edente implique
$$\forall q \in Q({\cal F}/P)\qquad H^{1}(X_{q},{\cal U}_{q})=0\quad . $$
\end{remq}
\begin{remq}
Vues les hypoth\`eses,   le rang de ${\cal U}_{q}$
est r, son degr\'e est $p-(r(1-g))$  et
$h^{0}(X_{q},{\cal U}_{q})=p$ . En fait le choix
du polynôme de Hilbert $P$ est fait de telle sorte
que ${\cal U}_{q}$ est ces caractéristiques en
commun avec $E$.
\end {remq}\\
\begin {remq}
Si on identifie le groupe $GL(p)$ avec le groupe des automorphismes
de ${\cal O}_{X}^{p}$, alors $GL(p)$ agit sur $Q({\cal F}/P)$ et
sur $\cal U$.  De plus cette action induit une action de $PGL(p)$
sur $Q({\cal F}/P)$ mais pas sur $\cal U$ (car $\lambda \times Id$
agit sur $\cal U$ en multipliant par $\lambda$).
\end{remq}
\begin{thm}\cite{seshadri} On a :
\begin{enumerate}
\item ${\mathcal {R}}$ est un ouvert $PGL(p)$-invariant de $Q({\cal F}/P)$
et ${\mathcal U}_{\vert {\mathcal {R}}\times X}$
est localement libre.  Il correspond donc \`a un
fibr\'e $U$ sur ${\mathcal {R}}\times X$.
\item $U_{q_{1}}\cong U_{q_{2}}$ si et seulement si $q_{1}$ et $q_{2}$
 appartiennent \`a la m\^eme orbite
sous l'action de $PGL(p)$.
\item $\forall q \in {\mathcal {R}}$ le stabilisateur de q sous l'action de $PGL(p)$
est isomorphe au quotient $Aut(U_{q})/\{\lambda \times Id\}$.
\end{enumerate}
\end{thm}
\begin{prop} ${\mathcal {R}}$ est un sous-sch\'ema ouvert lisse de dimension
$p^{2}+r^{2}(g-1). $
\end{prop}

 On consid\`ere maintenant l'ensemble
${\mathcal {R}}^{^{ss}}$(resp. ${\mathcal
{R}}^{s}$) des points de ${\mathcal {R}}$ tels que
${\cal U}_{q}$ soit semi-stable (resp. stable).
\\ Comme pour tout $x$ appartenant \`a $X$ le
fibr\'e $U_{\vert {\mathcal {R}}\times
\{x\}}$ peut \^etre regard\'e comme une famille de quotients de
dimension $r$ de l'espace vectoriel $H^{0}({\cal
O}_{X}^{p})$ (de dimension p), on peut construire
le morphisme suivant, en notant
$H_{p,r}:=Grass(p,r)$  :
\begin{eqnarray*}
\psi_{x} : \,{\mathcal {R}}&\lr& H_{p,r}\mbox{, d\'efini par :}\\
q&\mapsto&\psi_{x}(q)=(U_{q})_{x}.
\end{eqnarray*}
$\psi_{x}$ est un $PGL(p)$-morphisme i. e.   $\psi_{x}(g. q)=g.
\psi_{x}(q)\quad\forall g\in PGL(p)$.  Prenons ensuite une suite de
points $x_{1},x_{2}\ldots x_{s}$ de $X$.  On peut
donc définir :
\begin{eqnarray*}
\psi : {\mathcal {R}}&\longrightarrow & \prod _{i=1}^{i=S}(H_{p,r})_{i}\mbox{ d\'efini par }\\
q&\mapsto&\psi(q)=((U_{q})_{x_{1}},(U_{q})_{x_{2}},\ldots,(U_{q})_{x_{s}}).
\end{eqnarray*}
Ce morphisme est aussi un $PGL(p)$-morphisme.  On a alors le th\'eor\`eme
suivant qui nous sera utile pour la construction de la fonctionnelle de
Kempf et Ness (cf. paragraphe \ref{ken}  ).
:
\begin{thm}\cite{seshadri} \label{psy}Pour r fix\'e, il existe $d_{0} \in \mathbb {N}$ tel que $\forall d>d_{0}$
, il existe une suite de $s:=s(d)$ points telle que le morphisme
$\psi
: {\mathcal {R}} \longrightarrow \prod\limits _{i=1}^{i=s}(H_{p,r})_{i}=:Z$ v\'erifie :
 \begin{enumerate}
\item $\psi$ est injective.
\item ${\mathcal {R}}^{ss}=\psi^{-1}(Z^{ss})$.
\item ${\mathcal {R}}^{s}=\psi^{-1}(Z^{s})$.
\item $\psi : {\mathcal {R}}^{ss}\longrightarrow Z^{ss}$ est propre.
\end{enumerate}
Ici l'action considérée sur $Z$ est celle de
$SL(p)$ et, comme pour ${\mathcal {R}}^{s}$ et
${\mathcal {R}}
^{ss}$,  les notations $Z^{s}$ et $Z^{ss}$
désignent les les points stables ou semi-stables
de $Z$ comme défini dans \cite{potier}.
\end{thm}
On peut caractériser la stabilité sur $Z$ de la
façon suivante :
\begin{prop}
Soit $v=(v_{1},v_{2},\ldots ,v_{s})\in Z$.  On
pose,
$$\delta(D)=\sum_{i=1}^{i=s}(\frac{dim \pi_{i}(D)}{dim D}-\frac{r}{p}),\,D\subset H^{0}({\cal
O}_{X}^{p})$$ où $\pi_{i}$ correspond à la
projection sur $(H_{p,r})_{i}$.
 Alors,

$$v\ est\ stable\ (resp. \
semi-stable)
$$si et seulement si  $$\forall D
\mbox{ sous-espace propre de } H^{0}({\cal
O}_{X}^{p}) , \qquad\delta(D)>0\,(\mbox {resp. }
\geq 0)\quad . $$

\end{prop}
Même si cela n'intervient pas par la suite, nous
remarquons que le th\'eor\`eme pr\'ec\'edent nous
permet de construire la  variété  des modules de
fibr\'es stables qui correspond au quotient
géométrique par $PGL(p)$ de ${\mathcal {R}}^{s}$.
La variété ${\mathcal {R}}^{ss}$ n'admettant quant
à elle, qu'un bon quotient par $PGL(p)$ qui est
une variété projective. Pour les définitions de
bon quotient et de quotient géométrique on
consultera \cite {newstead} p. 70. On rappelle
qu'un quotient géométrique est un bon quotient qui
est aussi un orbi-espace. On a également le
r\'esultat suivant
:
\begin{thm}
 Il existe un espace
de module grossier pour les fibr\'es vectoriels de rang r et de
degr\'e d.  Cet espace a une compactification naturelle en une
vari\'et\'e projective not\'ee ${\cal M}(r,d)$ avec $dim({\cal
M}(r,d))=r^{2}(g-1)+1$.
\end{thm}
Dans la seconde partie, cette construction
pr\'ecis\'ee ici et plus particulièrement  le
théorème \ref{psy}, nous permettra d'\'etablir une
d\'efinition de la fonctionnelle de G. Kempf et L.
Ness
.

\chapter{Espace des connexions sur $E$ et fibrés de Yang-Mills}
On pose :
$$\cin(X, \Lambda ^{q}T^{*}X\otimes F)=\Omega ^{q}(F)\,$$
avec $F$ un fibré sur $X$. Nous introduisons les
notations qui nous serons utiles par la suite
.
\section{Connexion et courbure}
 Soit E un fibr\'e
vectoriel $\cin$ sur $X$ de rang $r$ et de degr\'e
$d$ (ici on consid\`erera une vari\'et\'e X
complexe ou r\'eelle de dimension $n$, \'etant
entendu que le cas qui nous int\'eresse est celui
qui correspond \`a une vari\'et\'e complexe de
dimension complexe $1$).
\begin{dfn}\index{Connexion}
Une connexion D sur E est un op\'erateur
différentiel lin\'eaire D tel que :
\begin{eqnarray*}\index{Connexion}
D : \cin(X, \Lambda ^{q}T^{*}X\otimes E)\longrightarrow\cin(X,
\Lambda
^{q+1}T^{*}X\otimes E)\mbox{ avec , }\\ D(\alpha
\wedge u)=d\alpha\wedge u+(-1)^{\deg\,
\alpha}\alpha\wedge Du \mbox{    (R\`egle de Leibnitz)}
\end{eqnarray*}
 o\`u $\alpha \in \cin(X, \Lambda
^{q+1}T^{*}X)$ et $u\in
\cin(X, \Lambda ^{q}T^{*}X)$.
\end{dfn}
\begin{remq} \label{derico}

 On notera   $\cov^{E}$ la dérivée covariante de
$E$ définie par $D $ et également $\cov_{h}$ la
dérivée covariante de $(E,h)$.\label{deltah}
\end{remq}\\
Supposons de plus que $E$ soit muni d'une
m\'etrique $h$, hermitienne (ou euclidienne ), ce
que l'on notera $(E, h)$ et soit
  un champ de
repères $\cin$, $(e_{\lambda})_{\lambda=1. . r}$
de E.  Nous avons alors l'accouplement
sesquilin\'eaire (ou linéaire) canonique
:
  \begin{eqnarray*}
\cin(X, \Lambda
^{p}T^{*}X\otimes E)\times\cin(X, \Lambda
^{q}T^{*}X\otimes E)&\longrightarrow&
\cin(X, \Lambda ^{p+q}T^{*}X\otimes \mathbb C)\\
(u, v)&\longmapsto& \{u, v\}
  \end{eqnarray*}
o\`u
$$\{u,
v\}=\sum_{\lambda, \mu=1}^{r}u_{\lambda}\wedge
 v_{\mu}<e_{\lambda}, e_{\mu}>_{h}\mbox { si X est réelle }$$
et$$ \{u, v\}=\sum_{\lambda,
\mu=1}^{r}u_{\lambda}\wedge
\bar v_{\mu}<e_{\lambda}, e_{\mu}>_{h}\mbox { si X est complexe }.$$
On a posé :
\begin{eqnarray*}
u&=&\sum_{\lambda=1}^{r}u_{\lambda}\otimes
e_{\lambda}\\ v&=&\sum_{\mu=1}^{r}v_{\mu}\otimes
e_{\mu}
\end{eqnarray*}
\begin{dfn}\label{hermi}
Soit $D$ une connexion sur $E$. Si $$d\{u,
v\}=\{Du, v\}+(-1)^{\deg
\,  u}\{u, Dv\}$$
alors
\begin{itemize}
  \item D est dite  euclidienne , si  $h$ euclidienne . Nous parlerons dans  ce
   cas de fibré
euclidien $(E,h)$,
\item D est dite  hermitienne , si $h$ hermitienne. $(E,h)$ sera dit fibré
hermitien.
\end{itemize}
\index{Connexion!hermitienne}

\end{dfn}
Si maintenant le fibr\'e $E$ est holomorphe sur
$X$, vari\'et\'e complexe, on a
:
\begin{dfn}\label{holo}
Une connexion $D$ est compatible avec la structure
holomorphe si et seulement si
$$ \forall s \in
\cin(X, T^{*}X\otimes E)\qquad D^{''}s=\db s \quad,$$ o\`u $D^{''}$ est la (0, 1)-composante de $D$.
\end{dfn}
\begin{dfn}
Une connexion D est dite de Chern si elle
v\'erifie les d\'efinitions \ref{hermi} et
\ref{holo}. \index{Connexion!de Chern}
\end{dfn}
\begin{prop}\cite{koba}
La connexion de Chern  D sur un fibr\'e holomorphe
hermitien (E, h) existe et  est unique.
\end{prop}
On reprend le cas d'un fibré vectoriel $\cin$ sur
$X$.
\begin{dfn}\index{Courbure}
On appelle courbure de la connexion D,  la 2-forme
R(D) \`a valeurs dans $End(E)=E\otimes E^{*}$
d\'efinie par :
  \begin{eqnarray*}R(D) : \cin(X, T^{*}X\otimes
E)&\longrightarrow &\cin(\Lambda^{2}T^{*}X\otimes
E)\\
 R(D)&=&D\circ D.\\
  \end{eqnarray*}

On omettra de pr\'eciser la connexion en question
si le contexte est suffisamment clair.  Lorsqu'il
s'agira de la courbure de la connexion de Chern du
fibré $(E, h)$, nous utiliserons également la
notation $R_{h}(E)$.
\end{dfn}
\begin{dfn}\index{Connexion!plate}
Une connexion D sur E est dite plate si $R=0$
\end{dfn}
\begin{prop}\label{loc}
Soit D la connexion de Chern d'un fibr\'e
holomorphe hermitien (E, h) de rang $r$ sur $X$.
On a alors la d\'ecomposition suivante pour une
trivialisation donn\'ee $\phi : E_{\vert U}
\stackrel{\sim}{\longrightarrow} U\times
\mathbb{C}^{n}$:
$$\left\lbrace
\begin{array}{l}
D^{'}u=_{\vert \phi}\partial u+H^{-1}\partial H\wedge u\\
 D^{"}=_{\vert \phi}\db u
 \end{array}
 \right. $$
 o\`u $u\in \cin(X, E)$ et $H=(h_{\lambda\bar \mu})$ avec $h_{\lambda\bar\mu}=h(e_{\lambda}, e_{\mu})$ avec
 $(e_{\lambda})_{\lambda =1. . r}$ rep\`ere holomorphe local .

\end{prop}

\begin{remq}\label{remqu}
\begin {enumerate}
\item Plus g\'en\'eralement on a $D=_{\vert \phi}d+A$ avec
\begin{itemize}
\item  si $E$ est un fibr\'e vectoriel sur $X$ ,
vari\'et\'e r\'eelle, $A\in
\cin(\mathfrak{gl}(r, \mathbb{R})\otimes T^{*}X_{\vert U})$ o\`u $\mathfrak {gl}(r, \mathbb{R})$ correspond \`a l'alg\`ebre de Lie
de $GL(r, \mathbb{R})$,
  \item  si $E$ est un fibr\'e vectoriel complexe sur $X$ ,
vari\'et\'e complexe, $A\in
\cin(\mathfrak{gl}(r, \mathbb{C})\otimes T^{1,0}X_{\vert U})$ o\`u $\mathfrak {gl}(r, \mathbb{C})$ correspond \`a l'alg\`ebre de Lie
de $GL(r, \mathbb{C})$,

  \item  si $E$ est un fibr\'e vectoriel vari\'et\'e r\'eelle $X$, muni d'une m\'etrique
h, et D euclidienne pour le fibr\'e $(E,h)$ ,
$A\in
\cin(\mathfrak{o}(r,
\mathbb{R})\otimes T^{*}X_{\vert U})$ o\`u
$\mathfrak {o}(r, \mathbb{R})$ correspond \`a
l'alg\`ebre de Lie de $O(r, \mathbb{R})$,

  \item  si $D$ est la connexion de Chern du fibr\'e
  $(E,h)$ sur la vari\'et\'e complexe $X$,
  $A\in \cin(\mathfrak{u}(r)\otimes T^{1,0}X_{\vert U})$.

\end{itemize}

La diff\'erence de deux connexions est donc un
objet  de $\Omega^{1}(\mathfrak{gl}(E))$ et
 si il s'agit de fibré $(E,h)$ euclidien (resp. holomorphe hermitien) elle est
  dans $\Omega^{1}(\mathfrak{o}(E))$ (resp.
$\Omega^{1,0}(\mathfrak{u}(E)))$.

 De plus $$R(D)=_{\vert \phi}dA+A\wedge A.
$$
\item La courbure $R$ associ\'ee \`a une connexion $D$ compatible  avec la structure
holomorphe d'un fibr\'e $E$ sur $X$ est une (1,
1)-forme \`a valeurs dans $End(E)$ et dans ce cas
$$R(D)=_{\phi}D^{'}\circ D^{"}+D^{"}\circ D^{'}$$ o\`u
$D^{'}$ est la $(1, 0)$-composante de $D$.
D'apr\`es la  première remarque
$$R(D)=dA+A\wedge A=_{\vert
\phi}H^{-1}\db\partial H+H^{-1}\partial H\wedge H^{-1}\db H$$ en utilisant les notations
 de la proposition \ref{loc}.
\item Dans le cas d'une connexion  de Chern $D$ sur un fibr\'e holomorphe hermitien $(E, h)$,
la courbure R vérifie  que $$iR(D)\in
\cin(X, \Lambda^{1, 1}T^{*}X\otimes Herm(E, E)). $$
Ici $iHerm(E, E)$ correspond \'egalement  \`a
$Ad(E)$ et \`a $\mathfrak{u}(E)$,
 c'est \`a dire a\label{und}ux endomorphismes   hermitiens du fibr\'e.
\end{enumerate}
\end{remq}

\begin{prop}\cite{jost1}
La courbure R d'une connexion D d'un fibr\'e E sur X satisfait la
seconde identit\'e de Bianchi i. e\index{Identité de Bianchi}
$$DR=0. $$

\end{prop}
%%%%%%%%%%%%%%%%%%%%%%%%%%%%%%%%%%%%%%%%%%%%%%%%%%%%%%%%%%%%%%%%%%%%%%%%%
\section{Fibr\' es Einstein-Hermite }
Nous nous plaçons jusqu'\`a la fin du paragraphe
sur une vari\'et\'e  kälhérienne $X$ de dimension
r\'eelle $n$ o\`u $n=2m$. Soit $E$ un fibr\'e
holomorphe sur $X$ de rang $r$. Soit $h$ une
métrique sur $E$.
 Soit $(e_{\lambda})_{\lambda=1. . r}$ un champ de rep\`ere sur $E$ et
 $(z^1, z^2, . . , z^m)$ est un syst\`eme
de coordonn\'ees de $X$ . On note :
$$h_{\lambda\bar \mu}=h(e_\lambda, e_\mu)\mbox{ et }  g=
\sum {g_{i\bar j}}dz^{i} \otimes d{\bar z}^j . $$
La courbure $ R_{h}=(R_{\lambda}^{\mu})\in
\cin(End(E)\otimes \Lambda ^{2}T^{*}X)$  s'\'ecrit  par rapport
\`a $(e_{\lambda})$ :
$$R(e_\lambda)=\sum_{\mu=1..r}{R^\mu_\lambda} e_\mu \mbox{ et }
 {R^\mu_\lambda}=\sum_{i,j=1..m}
{R^\mu_{\lambda i\bar j}}dz^i \wedge d{\bar z}^j.
$$

La courbure peut se d\'ecomposer en deux
composantes dont l'une a pour définition suivante
:\\
\begin{dfn}\label{courbuprinc}
La courbure centrale K de (E,
h)\index{Courbure!principale} est d\'efinie par
:
$${K^\mu_\lambda} =\sum g^{i\bar j }{R^\mu_{\lambda i\bar j}}:=tr_{g}(R).
 $$ Alors $({K^\mu_\lambda})$
 d\'efinit un
endomorphisme $K$ sur $E$.
\end{dfn}
\begin{remq}
Si on prend $E=TM$ le fibr\'e tangent holomorphe  et $h=g$ alors
$$K=\mbox{ric}_{X}$$est la courbure de Ricci de X. \index{Courbure!de Ricci}
\end{remq}
\begin{dfn}
 Un fibr\'e $(E, h)$ sur $(X, \omega )$ satisfait la condition Einstein faible si
$$K=\varphi I_E$$ o\`u  $\varphi $ est une fonction r\'eelle . \\
 Si de plus $\varphi =cste$ alors\index{Fibré!Einstein Hermite}
le fibr\'e (E, h) est dit Einstein Hermite ou Yang-Mills Hermite.
\end{dfn}
\begin{remq}
Si $(E, h)$ est Yang-Mills Hermite alors la constante $\varphi$ est
\'egale \`a $$\varphi=\frac{2}{vol(X)}. \frac{\deg(E)}{\rg(E)}.$$
\end{remq}
\begin{remq} Si X est une courbe alg\'ebrique on a alors :
$g=g_{11}dz \wedge d{\bar z}$.  Ainsi,
$K^{\lambda}_{\mu}=g^{11}{R^\mu_{\lambda 11} }$.
Donc la courbure centrale est toute la courbure et
il n'y a qu'une composante.
\end{remq}\\
Cette notion de fibr\'e Yang-Mills Hermite sera utile dans la deuxi\`eme
partie avec l'introduction de la fonctionnelle de Donaldson.
%%%%%%%%%%%%%%%%%%%%%%%%%%%%%%%%%%%%%%%%%%%
\section{Généralités sur l'espace des connexions sur $E$}\label {symp}
Après les rappels du paragraphe précédent, nous
étudions l'espace des connexions d'un fibré $E$
holomorphe. On se place donc dans le cas où $X$
est une surface de Riemann. Soit $D_{\alpha},
D_{\beta}$ deux connexions sur un fibré $E$ qui
sont hermitiennes. Nous avons remarqué (cf.
paragraphe \ref{remqu})  que la différence de
$D_{\alpha}$ et de $D_{\beta}$ était une $1-$forme
à valeurs dans $\mathfrak {u}(E).$ Donc,  si on
considère l'ensemble des connexions sur $E$, que
l'on notera $\mathcal {A}$, il s'agit d'un espace
affine de direction vectorielle
$\Omega^{1}(\mathfrak {u}(E))$ où
$\Omega^{1}(\mathfrak {u}(E))$ désigne l'espaces
des sections $\cin$ de $\mathfrak {u}(E)\otimes
T^{*}X$.
 On  munit $\espco$ d'une métrique et la
forme associée définie par
\index{Notation!$\omega_{\espco}$}$$
\omega_{\espco}(D_{1}, D_{2})=\frac{1}{4\pi^{2}}\int_{X}tr(D_{1}
\wedge D_{2})\quad \forall D_{1}, D_{2}\in T\espco
$$ en fait un espace symplectique.\label{metria}
Regardons, dans le paragraphe suivant, l'action du
Groupe de Jauge.
\section{Groupe de Jauge}\index{Groupe de Jauge}\index{Notation!$\mathcal{G}$}
\label{jaug}
Soit $GL(E)$, le fibré des automorphismes
linéaires complexes du fibré $E$ holomorphe sur X
. Nous définissons le Groupe de Jauge par
 $$\mathcal{G}=\{s\in \cin(GL(E))/ ^{t}\bar s s=id_{E}\}.$$
 Soit $D\in \espco$ avec $D=\db+\partial$,  la décomposition habituelle.
 L'action de $s\in\mathcal{G}$ notée $s^{*}(D)=s^{-1}\circ D\circ s$ est donnée par :
\begin{equation*}
\begin {cases}
\db_{s(D)}&=s^{-1}\db_{D} s^{-1}\\
\partial _{s(D)} &=(^{t}\bar s)\partial _{D} (^{t}\bar s)^{-1}.
\end{cases}
\end{equation*}

 Les
éléments de ce groupe correspondent aux isométries
de la métrique mise sur $\espco$.\\ Par ailleurs
le groupe\index{Notation!$\mathcal{G}$}
$$\jaugec= \cin(GL(E))
$$ agit aussi sur $\espco$ de la façon suivante :
\begin{equation*}
\begin {cases}
\db_{g(D)}&=g^{-1}\db_{D} g
\\
\partial_{g(D)} &=(g^{*})\partial (g^{*})^{-1}
\end{cases}
\end{equation*}
avec $g\in \jaugec$ et $D\in \espco$.\\L'action du
groupe $\jaugec$ ne préserve pas la métrique mise
sur $\espco$.
 Cependant deux connexions définissent des structures holomorphes isomorphes pour $E$
si et seulement si  elles sont dans la même orbite
sous l'action de $\jaugec$. En effet par
définition deux opérateurs  $\bar\partial_{1}$ et
$\bar\partial_{2}$ définissent des structures
isomorphes si elles sont conjuguées par un élément
de $\cin(GL(E))$.

\section{Op\' erateur $*_{E}$  }

 Nous consid\'erons pour commencer
 une vari\'et\'e riemanienne orientée compacte lisse $(X,
g)$ de dimension réelle
 $m$.
 On  d\'efinit un op\'erateur lin\'eaire $*$ sur $\Lambda^{p}T_{x}^{*}X$
 de la façon suivante :

  \begin{eqnarray*}* :
   \Lambda^{p}T_{x}^{*}X
&\longrightarrow &\Lambda^{m-p}T_{x}^{*}X \mbox {
avec } \\ f_{i_{1}}\wedge \ldots \wedge
f_{i_{p}}&\longmapsto& f_{j_{1}}\wedge \ldots
\wedge f_{j_{m-p}}
  \end{eqnarray*} de telle sorte  que
$(f_{i_{1}}, \ldots, f_{i_{p}}, f_{j_{1}}, \ldots,
f_{j_{m-p}})$ soit un base orientée positive. On
considère maintenant une variété compacte
kälhérienne de dimension $n$ où $2n=m$. On
considère aussi un fibré holomorphe hermitien $E$
sur $X$.
 Soit $\tau :
E\longrightarrow E^{*}$  isomorphisme de fibr\'es
défini via $h.$ On d\'efinit alors,
  \begin{eqnarray*}
   *_{E} : \cin
(\Lambda^{p}T_{c}^{*}X\wedge E)&\longrightarrow&
\cin (\Lambda^{2n-p}T_{c}^{*}X\wedge E^{*})\mbox{ avec
}\\
\bar
*_{E}(\phi\otimes e)&=&* \bar \phi \otimes \tau
(e)\mbox{ pour }e\in E_{x} \mbox{  et  }\phi
\in \cin (\Lambda^{2n-p}T_{c,x}^{*}X).\\
  \end{eqnarray*}

%%%%%%%%%%%%%%%%%%
\begin{remq}\label{met}
On d\'efinit donc sur (X, g) une m\'etrique sur
les q-formes en posant : $<\mu_{1},
\mu_{2}>_{g}=*(\mu_{1}\wedge*\mu_{2})$.  On a
\'egalement
$$*(1)=(det(g_{i\bar j}))^{\frac{1}{2}}dx^{1}\wedge\ldots
\wedge dx^{m}$$ c'est \`a dire la forme volume.
\end{remq}\\
%%%%%%%%%%%%%%%%%%%%

\section{Op\'erateur de Laplace-Dolbeault}\label {lapl}
Soit $D$ une connexion d'un fibré  $E$. On
d\'efinit  formellement l'adjoint de D par :
$$D^*=-*D*$$ ainsi que dans le cas d'un fibr\' e holomorphe hermitien $(E, h)$
l'adjoint de $\db$.
$$\dbaretid{h}=-*_{E^*}\db*_E. $$
On d\'efinit enfin l'op\'erateur de Laplace
Dolbeaut par :
$$\Delta_{\db}=\dbaret_{h} \db +\db \dbaret_{h} . $$

%%%%%%%%%%%%%%%%%%%%%%%%%%%%%%%%%%%%%%%%%%%%%%%%%
\section{Fonctionnelle de Yang-Mills }\label{ymills}
Soit $X$, une variété complexe, $(E,h)$ un fibré
hermitien holomorphe sur $X$. Nous pr\'ecisons
d'abord le produit hermitien sur $\cin(X, \Lambda
^{p}T^{*}X\otimes Herm(E, E))$.
\begin{prop}
Soient $A, B \in \mathfrak{u}(n)$ on d\'efinit un produit hermitien
sur $\mathfrak{u}(n)$ grâce \`a :
$$A. B=-tr(A.\bar B). $$
Cela d\'efinit un produit hermitien dans $Herm(E, E). $
\end{prop}

\begin{dfn}\label{metri}(cf p. \pageref{met})
Le produit hermitien dans  les fibres de $
\Lambda
^{p}T^{*}X\otimes Herm(E, E)$ est le suivant :
$$<\mu_{1}\otimes \omega_{1}, \mu_{2}\otimes
\omega_{2}>=\mu_{1}. \mu_{2}<\omega_{1}, \omega_{2}>_{g}$$ o\`u
$\mu_{1}\otimes \omega_{1}, \mu_{2}\otimes
\omega_{2}\in  Herm(E_{x}, E_{x})\otimes \Lambda
^{p}T_{x}^{*}X\quad \forall x\in X$. \\
Cela induit donc un produit $L^2$ dans $\cin(X, \Lambda
^{p}T^{*}X\otimes Herm(E, E))$  :
$$( \mu_{1}\otimes \omega_{1}, \mu_{2}\otimes
\omega_{2})=\int_X <\mu_{1}\otimes \omega_{1}, \mu_{2}\otimes
\omega_{2}>*(1). $$
\end{dfn}
On peut maintenant d\'efinir la fonctionnelle de Yang-Mills.
\begin{dfn}\index{Fonctionnelle!de Yang-Mills}\label{yang}
Soient X une vari\' et\' e kälhérienne compacte
 $(E, h)$ un fibr\'e
hermitien holomorphe dont la courbure de Chern est
$R(D)$. La fonctionnelle de Yang-Mills est d\'
efinie par\index{Notation!${\mathcal YM}(D)$}
$${\mathcal {YM}}(D):=(R(D), R(D)) =\int_{X} <R(D), R(D)>*(1).$$
\end{dfn}
\begin{prop}
Une connexion $D$ est de
Yang-Mills\index{Connexion!Yang-Mills}
\index{Fibré!Yang-mills} si
elle est stationnaire pour ${\cal YM}(D)$ c'est \`a dire quelle
v\'erifie :
$$D*R(D)=0. $$
Dans ce cas le fibré $(E, D)$ est dit fibré de Yang-Mills.
\end{prop}
\begin{demo}
Comme l'espace de connexions compatibles avec la
m\'etrique du fibr\'e $(E, h)$ sur $X$ sur $E$ est
un espace affine model\'e sur
$\Omega^{1}(\mathfrak{u}(E))$, la diff\'erence de
deux connexions compatibles est un
\'el\'ement de $\Omega^{1}(\mathfrak{u}(E))$. Il s'ensuit
que pour
\'etudier les points critiques de ${\cal YM}$ on regarde :
$$\text{Pour } A\in \Omega^{1}(Ad(E)),\quad s\in \Gamma(E),$$
\begin{eqnarray*}
R(D+tA)(s)&=&(D+tA)(D+tA)(s)\\
 &=&D^{2}s+tD(As)+tA\wedge Ds+t^{2}(A\wedge A )s\\
 &=&(R(D)+t(DA)+t^{2}(A\wedge A))s
 \end{eqnarray*}
 car $D(As)=(DA)s-A\wedge Ds$
 On obtient donc :
 \begin{eqnarray*}
 \dt{\mathcal {YM}(D+tA)}_{\vert t=0}&=&\dt \int_{X}<R(D+tA), R(D+tA)>*(1)_{\vert t=0}\\
 &=&2\int_{X}<\dt(R(D)+t(DA)+t^{2}(A\wedge A)), R(D+tA)_{\vert t=0}>*(1)\\
 &=&2\int_{X}<DA, R(D)>*(1)\\
 &=&2(DA, R(D)).
 \end{eqnarray*}
 D'o\`u,   $\displaystyle {\dt{\mathcal {YM}(D+tA)}_{\vert t=0}=0\Leftrightarrow  D^{*}R(D)=0. }$
\end{demo}

 \begin{thm}La fonctionnelle de Yang-Mills est invariante sous
l'action du groupe de Jauge ${\cal G}$.   De plus,
$\forall s \in {\cal G}$, si $D$ est une connexion
de Yang-Mills alors $s^*(D)=s^{-1}\circ D\circ s $
est aussi une connexion de Yang-Mills.
\end{thm}
\begin{demo}
Soit $s\in {\cal G}$,  $s$ agit sur l'espace des
connexions sur E de la façon suivante :
$$s^{*}(D)=s^{-1}\circ D\circ s $$
on a donc
\begin{eqnarray*}
s^{*}R(D)&=&s^{-1}\circ D\circ s \circ s^{-1}\circ D\circ s \\
&=&s^{-1}\circ R(D)\circ s
\end{eqnarray*}
Comme de plus $s$ est  un automorphisme du fibr\'e
$(E, h)$ qui  est une isom\'etrie pour $<. , . >$
il s'ensuit
:
$$<s^{*}R(D), s^{*}R(D)>=<R(D), R(D)>$$
\end{demo}
%%%%%%%%%%%%%%%%%%%%%%%%%%%%%%%%%%%%%%%%%%%%%%%%%%%

%%%%%%%%%%%%%%%%%%%%%%%%%%%%%%%%%%%%%%%%%%%%%%%%%%%%%%%%%%%%%%
\section{Relecture du th\'eor\`eme de M. S. Narasimhan et C. S. Seshadri par S.K.Donaldson}
Dans \cite{don1},  S. K.  Donaldson a \'etabli un
th\'eor\`eme liant les fibr\'es stables et les
fibr\'es dont la courbure centrale est constante.
Pour cela il utilise une fonctionnelle notée ici
$J$ qui correspond pour les fibrés de degré zéro
et de rang deux à la fonctionnelle de Yang-Mills.
Elle préfigure la fonctionnelle que nous appelons
fonctionnelle de S.K. Donaldson et noterons
$\mathcal{M}$ (cf. paragraphe \ref{cestfini}).
Voici donc le théorème en question.
\begin{thm}\cite{don1}\index{Théorème! de S.K Donaldson}\label{donalds}
Un fibr\'e holomorphe, indécomposable, $E$ sur $X$
est stable si et seulement il y a sur $E$ une
connexion unitaire $D$ telle que sa courbure
centrale satisfasse :
$$*F=-2\pi i \mu(E)\mbox{ avec } \mu(E)=\frac{\deg\, E}{rg\, E}.
 $$Cette connexion est unique \`a isomorphisme
pr\`es.
\end{thm}
\index{Courbure!centale}
 La d\'emonstration , donnée dans \cite{don1},  de ce th\'eor\`eme
s'organise de la façon suivante :
\begin{itemize}
\item Pour les fibr\'es en droite il s'agit d'une cons\'equence de la th\'eorie de Hodge.
\item On suppose que le r\'esultat a \'et\'e prouv\'e pour des fibr\'es de rang $r-1$.
On le d\'emontre pour des fibr\'es de rang $r$.
\end{itemize}
Pour cela on consid\`ere une suite de connexions
compatibles avec la m\'etrique de $(E, h)$,
fibr\'e $\cin$,  qui minimisent une fonctionnelle
$J$ définie par :
$$J(D)=N(\frac{*F(D)}{2i\pi}+\frac{\mu(E)}{vol\, X}Id_{E}) \,\mbox{ avec }D \mbox{ connexion sur }E,$$
 et $N(s)=\int_{X}(s^{*}s)^{\frac{1}{2}}\, ,s\in End(E)$.  On
extrait de cette suite une sous-suite convergeant
faiblement et on obtient deux cas :
\begin{itemize}
\item La connexion limite est dans l'orbite, pour l'action de $\jaugec$, des connexions unitaires
de E. Dans ce cas,  si $Inf J_{\vert Orb(E)}$ est
atteint dans $Orb(E)$, notons le $D_{0}$. Alors on
montre
 que la courbure $R(D_{0})$ de la connexion pour laquelle le minimum est atteint
 vérifie $*R(D_{0})=-2\pi i\mu(E)$ et
 que, quitte \`a changer de jauge, elle est continue
 et unique \`a l'action du groupe de jauge près.
\item La connexion limite est dans l'orbite d'un fibr\'e $F$ distinct de $E$
et dans ce cas le fibr\'e $E$ n'est pas stable.
\end{itemize}
La deuxi\`eme situation ne se produit pas si le fibr\'e est
irréductible.  Pour cela on utilise d'abord un résultat de
Uhlenbeck
:
\begin{prop}(\cite{uhlenbeck}) Supposons que $(D_{i})_{i\in \mathbb N}$ soit une suite de connexions unitaires
 $L^{2}_{1}$ sur E avec $\norm {R(D_{i})}$ born\'ee.  Alors il existe une sous-suite
 $(D_{\sigma(i)})_{i\in \mathbb N}$ avec $\sigma :\mathbb N\longrightarrow \mathbb N $ croissante
 et une transformation de Jauge $L^{2}_{2}$,  $u_{\sigma(i)}$ telle que  $(u_{\sigma(i)}
 (D_{\sigma(i)}))_{i\in \mathbb N}$ converge faiblement dans $L^{2}_{1}$.
 \end{prop}
 Cela nous permet d'obtenir le résultat suivant :
 \begin{prop}
 Soit $E$ un fibré holomorphe sur $X$.  Alors on a deux cas :
 \begin{itemize}
 \item $Inf\,  J_{\vert Orb(E)} $ est atteint dans $Orb(E). $
\item $Inf\,  J_{\vert Orb(E)} $ n'est pas atteint dans $Orb(E)$.  Il existe alors un fibré
holomorphe $F$ qui est de même degré et de même
rang que $E$ mais tel que $F\not\cong E$,   $Inf\,
J_{\vert Orb(F)}\leq Inf\, J_{\vert Orb(E)}$ et
$Hom(E, F)\not =0$.
 \end{itemize}
  \end{prop}
 Si on choisit un morphisme $f$ non trivial de $E$ dans $F$ on a alors le diagramme suivant :

\begin{equation*}
  \begin{CD}
 0@>>>  P @>>>\\
     \\
 0@<<< N @<<<
  \end{CD}\begin{CD}
 E   @>>>Q @>>> 0 \\
     @VV{f}V  @VV{g}V \\
 F   @<<<M @<<< 0
  \end{CD}
\end{equation*}

 avec $\deg Q\leq \deg M$ et $\rg Q =\rg M$.
 Or nous avons les deux propositions suivantes :
 \begin{prop}Si $F$ est un fibré holomorphe sur $X$ qui est une extension :
 $$0\longrightarrow  M\longrightarrow F \longrightarrow N \longrightarrow 0$$ et si
 $\mu(M)\geq \mu(F)$ alors pour toute connexion unitaire $D $ sur $F$ on a :
 $$J(D)\geq J_{0}, $$
 avec $$ J_{0}=\rg\, M(\mu(M)-\mu(F))+\rg \, N(\mu(F)-\mu(M)). $$
 \end{prop}
 \begin{prop}Si $E$ est un fibré stable et si le théorème a été prouvé pour des fibrés de
 rang inférieur,
 si il existe une extension :
$$0\longrightarrow P\longrightarrow E \longrightarrow Q \longrightarrow 0\, ,$$
alors il existe une connexion $D$ sur E telle que
$$J(D)<J_{1}, $$ avec $$J_{1}=\rg \, P. \mu(E)+
\rg \, Q(\mu(Q)-\mu(E)). $$
 \end{prop}
 Donc si on suppose que $E$ est stable et si le théorème est prouvé pour les fibrés de rang inférieur
 alors dans le cas de la seconde situation , il existe un fibré $F$ de même degré et de même rang que $E$
 avec $Hom(E, F)\not =0$ et $J_{0}\leq Inf\,  J_{\vert Orb(F)}\leq Inf\,  J_{\vert Orb(E)}$.
 Or d'après la dernière proposition $Inf\,  J_{\vert Orb(E)}<J_{1}$ mais comme
 $$\left\lbrace
 \begin{array}{l}
 \rg(Q)=\rg(M)\\
 \rg(P)=\rg(N)\\
 \deg(Q)>\deg(M)\\
 \deg(P)<\deg(N)
 \end{array}\right . \mbox{ on a } J_{1}\leq J_{0}. $$
 Donc le deuxième cas ne se produit pas si $E$ est stable.\\

\section{Conclusion}
Nous avons vu dans cette partie  deux visions,
l'une plutôt algébrique et l'autre plus analytique
de la stabilité. On va dans la prochaine partie
étudier plus profondément les liens entre ces deux
visions. Les outils mis en lumière dans cette
première partie seront d'ailleurs utilisés dans
les deux parties suivantes.
%%%%%%%%%%%%%%%%%%%%%%%%%%%%%%%%%%%%%%%%%%%%%%%%%%%%%%%%%%%%%%%%%%%%%%%
%%%%%%%%%%%%%%%%%%%%%%%%%%%%%%%%%%%%%%%%%%%%%%%%%%%%%%%%%%%%%%%%%%%%%%%%
%%%%%%%%%%%%%%%%%%%%%%%%%%%%%%%%%%%%%%%%%%%%%%%%%%%%%%%%%%%%%%%%%%%%%%%
%%%%%%%%%%%%%%%%%%%%%%%%%%%%%%%%%%%%%%%%%%%%%%%%%%%%%%%%%%%%%%%%%%%%%%%%
%%%%%%%%%%%%%%%%%%%%%%%%%%%%%%%%%%%%%%%%%%%%%%%%%%%%%%%%%%%%%%%%%%%%%%%
%%%%%%%%%%%%%%%%%%%%%%%%%%%%%%%%%%%%%%%%%%%%%%%%%%%%%%%%%%%%%%%%%%%%%%%%
%%%%%%%%%%%%%%%%%%%%%%%%%%%%%%%%%%%%%%%%%%%%%%%%%%%%%%%%%%%%%%%%%%%%%
%%%%%%%%%%%%%%%%%%%%%%%%%%%%%%%%%%%%%%%%%%%%%%%%%%%%%%%%%%%%%%%%%%%%%%

%%%%%%%%%%%%%%%%%%%%%%%%%%%%%%%%%%%%%%%%%%%%%%%%%%%%%%%%%p2bis3.tex}

\part{Liens entre les  espaces $Met(E_{n})$ et $Met(W_{n})$ : comparaison de métriques}

\section*{Introduction}\addcontentsline{toc}{section}{Introduction}
Nous désignons par $X$  une courbe alg\'ebrique
lisse compacte projective sur ${\mathbb {C}}$
plongée dans $\C P^{l}$. On
 munit  $\C P^{l}$ de la métrique
classique de Fubini-Study. On pourra considérer
$X$ comme une vari\'et\'e kälh\'erienne et on
notera $\omega$ la forme kälh\'erienne associ\'ee.
On considère un fibré vectoriel holomorphe  sur
$X$, noté $E_{0}$, de degré quelconque et de rang
deux. Soit $E_{n}$ le fibré $E_{0}\otimes
\ox(n)$.\\Pour comprendre les relations entre les
deux caractérisations de la stabilité décrites
dans la première partie, nous abordons, dans cette
seconde partie, les liens qu'il y a entre l'espace
des métriques sur $E_{n}$, noté $Met(E_{n})$ et
celui des métriques sur $H^{0}(X,E_{n})$ noté
$Met(W_{n})$.On note $p_{n}$ la dimension  de
$H^{0}(X,E_{n})$. Nous rappelons que l'espace
$Met(E_{n})$ reste inchang\'e pour tout $n$, alors
que $Met(W_{n})$ voit sa dimension augmenter avec
$n$.
\\Nous
construisons tout d'abord de façon naturelle, les
morphismes reliant les deux espaces $Met(W_{n})$
et $Met(E_{n})$ que nous appelons $I_{n}$ et
$L_{n}$.
 Nous d\'ecrivons les objets mis en
en jeu et mettons en lumière  les deux
fonctionnelles qui vont nous intéresser
: une, analytique, de S.K. Donaldson notée $\mathcal M$, ne d\'ependant pas de $n$
et d\'efinie sur $Met(E_{n})$ et l'autre,
algébrique, $\mathcal {KN}_{n}$ d\'efinie sur
$Met(W_{n})$ et construite  explicitement  p.
\pageref{kkkk}. L'idée initiale de cette
construction provient du travail décrit au
paragraphe \ref{cotal}. La première  fonctionnelle
permet de déterminer les points de $Met(E_{n})$
qui sont des métriques de  Yang-Mills Hermite et
la seconde les fibr\'es stables.\\
  Ensuite nous examinons  un côté géométrique de l'espace $\espco$ des
connexions   en utilisant essentiellement la
construction de la métrique de Quillen. Il en
découlera en fait assez naturellement le résultat
de S.K. Donaldson donné dans le chapitre
\ref{normap}. Ce résultat relie en partie les deux
espaces $Met(E_{n})$ et $Met(W_{n})$ via la
torsion analytique abordée dans  le chapitre p.
\pageref{torsi}.
C'est enfin dans le dernier chapitre de cette
partie que nous énonçons et démontrons les
résultats de ce travail dont un des plus
intéressant est :
\begin{thm*}(cf. p. \pageref{top})
Soit  $E_{n}$ un fibré stable sur $X$. Soit $U$ un
compact convexe donné de $Met(E,C^{4})$. Supposons
que $h_{YM}$, la métrique de Yang-Mills (elle
existe car $E_{n}$ est stable), appartienne à $U$.
On a alors
:
$$\forall \varepsilon >0 \quad \exists n_{U},\quad \forall n\geq n_{U},$$
il existe  m  extrémum de
${\mathcal{KN}_{n}}{\vert}_{ L_{n}(U)}$ et si on
pose
 $c=\frac{n}{\pi}(1+O(1))$ alors
$$d(h_{YM},cI_{n}(m))_{L^{2}_{1}(Met(E_{n}))}\leq \varepsilon. $$

\end{thm*}

Autrement dit la métrique de Yang-mills est
approximée par des métriques de la forme
$cI_{n}(m)$ où $m\in Met(W_{n})$ est un minimum de
la fonctionnelle algébrique $\mathcal KN_{n}$ sur
$L_{n}(U)$.

\chapter{L'espace des m\'etriques sur $E_{n}$ et l'espace des m\'etriques sur $W_{n}$ : les fonctionnelles ``alg\'ebriques"
et ``analytiques "}

\section{Morphismes entre $Met(W_{n})$ et $Met(E_{n})$.}

On pose :
$$W_{n}=H^{0}(X,E_{n})$$\index{Notation!$W_{n}$}
On d\'esigne par $Met(W_{n})$ (respectivement
$Met(E_{n})$) l'ensemble des m\'etriques sur
l'espace $W_{n}$ (resp. sur le fibré $E_{n}$). On
peut remarquer que :
 $$Met(E_{n_{1}})=Met(E_{n_{2}})\,
\quad
\forall n_{1},n_{2}\in \N \mbox{cf. p. \pageref{indep}}.$$
Nous ne distinguerons pas  ces espaces par la
suite.\\

 On s'intéresse au passage d'une
m\'etrique sur le fibr\'e $E_{n}$ à une m\'etrique
sur $H^{0}(X,E_{n})$. Pour cela nous introduisons
deux op\'erateurs :
\begin{dfn}\label{l}On note $$L_{n} :Met(E_{n})\lr Met(W_{n}),$$\index{Notation!$L_{n}$}
l'opérateur défini par :

 $$L_{n}(h)=\int_{z\in X}h(.,.)(z)\omega \ ,$$ avec
 \begin{eqnarray*}
\forall z\in X\,,\quad  h(.,.)(z) : W_{n}\times W_{n}&\,& \longrightarrow {\mathbb C}\\
(u,v)&\,& \longmapsto h_{z}(u(z),v(z))\quad
.
\end{eqnarray*}
Cette application fait correspondre à un choix de
métrique  $h$ sur $E_{n}$  la métrique $L^{2}$ sur
$W_{n}$ définie par $h$.
\end{dfn}
\begin{remq} Soit $\bar E_{n}$ le fibré conjugé du fibré $E_{n}$
et $\bar W_{n}$ l'espace conjugué de $W_{n}$. On
considère l'application naturelle :
 $$A_{n} : W_{n}\otimes \bar W_{n}\longrightarrow \cin(E_{n}\otimes \bar E_{n}).$$
Comme $Met(W_{n})\subset W_{n}^{*}\otimes \bar
W_{n}^{*}$ et $Met(E_{n})\subset
 \cin(E_{n}^{*}\otimes \bar E _{n}^{*})$, l'application $L_{n}$
 est la transpos\'ee du morphisme $A_{n}$.
 C'est un morphisme lin\'eaire.
\\
 \end{remq}\\
  Maintenant, si le degré du fibré $E_{n}$ est assez grand, alors le fibré $E_{n}$ est engendré par ses sections c'est à dire
 nous avons le morphisme surjectif (cf. paragraphe \ref{phix}):
 $$W_{n}\otimes \ox \overset{\phi}{\lr} E_{n} .$$
 Ainsi, chaque fibre $E_{n,z}$ avec $z\in X$ du fibré peut être vue, via le morphisme $\phi_{z}$,
  comme  quotient de l'espace vectoriel
 $W_{n}$. Si on se donne une métrique $m\in Met(W_{n})$,
   $\phi_{z}$ donne une métrique induite sur $E_{n,z}$.
 Plus précisément on a la définition suivante :
\begin{dfn}\index{Notation!$I_{n}$}\label{i}
Le morphisme $I_{n}$ est défini  de la façon
suivante
:
\begin{eqnarray*}
I_{n} : Met (W_{n})&\lr& Met(E_{n})\\
 I_{n}(m)(e_{z},e_{z})&=& min\{m(v,v)/v\in W_{n}\mbox{ et }v(z)=e_{z}\},
\end{eqnarray*}
avec $e_{z}\in E_{z}$, $m\in Met(W_{n})$.
\end{dfn}

\section{Fonctionnelles sur $Met(E_{n})$}
Nous commençons par fixer une application
déterminant dans $Met(E_{n})$. Soit $h\in
Met(E_{n})$\label{nopp}. Pour tout $z\in X$,
$h_{z}$ est une forme hermitienne. Par conséquent
on peut choisir une application déterminant,
$$det_{E_{n}}(h) :X\lr
\R_{*}^{+},$$ à partir d'un choix  d'un
isomorphisme $\Lambda^2E\cong \mathcal
{O}_X(2d_{n})$ où $d_{n}=deg(E_{n})$.\\ Dans toute
la suite du texte nous prenons donc un point de
référence $k_{0}\in Met(E_{n})$ tel que $$\det
_{E_{n}}(k_{0})=1.$$
Remarquons que si on se donne $h \in Met(E_{n})$
on a alors
: $k_{0}^{-1}h=H$ avec H endomorphisme autoadjoint
\label{metriend} défini positif du fibré $E_{n}$. Il vérifie :
$$\forall u,v \in E_{n}\quad<Hu,v>_{k_{0}}=<u,v>_{h}$$
Nous avons défini dans la première partie
 la fonctionnelle de Yang-Mills notée $\mathcal{YM}$\index{Fonctionnelle!de Yang-Mills}
   invariante sous l'action du groupe de Jauge ${\cal G}$ mais pas sous l'action du groupe
   complexifié ${\cal G}^{\mathbb C}$ (cf. p. \pageref{yang}).
 Nous
définissons maintenant  d'autres  fonctionnelles
sur $Met(E_{n})$.

\subsection{La fonctionnelle de Donaldson $M$}
On pourra consulter  \cite{don2} ou \cite{marg}.
\subsubsection {Morphisme de Chern-Weil.}
Soit $E_{n}$ un fibré vectoriel holomorphe de rang
$r$ sur $X$. Le cas qui nous intéresse est $r=2$,
mais les définitions qui suivent restent valables
pour $r$ quelconque. On note $Met(E_{n})$ l'espace
de toutes les métriques hermitiennes définies sur
$E_{n}$.
\begin{dfn}
Soit $P$
 une application $l$-linéaire sur $\mathfrak
{gl}(r,
\mathbb{C})$ symétrique invariante par $Gl(r,
\mathbb{C})$. On considère le morphisme
$\tilde P$ défini par
:
\begin{eqnarray*}
\tilde P
\, :&\,&Met(E_{n})\lr \Omega ^{(l,l)}(X)\\
&\,& h\,\mapsto P(iR(D),iR(D),...,iR(D))
\end{eqnarray*}
où $D$ est l'unique connexion associée à $h$. Ce
morphisme s'appelle morphisme de Chern-Weil.
\index{Morphisme de Chern-Weil.}

\end{dfn}
\begin{remq}
\begin{enumerate}
  \item Soient $h\in Met(E_{n})$ et  $H\in T_{h}Met(E_{n})$. Alors on a l'application tangente :

$$(\delta\tilde P)_{h}(H)=i.l.\db\partial P(\theta(H),iR(D),\ldots,iR(D))$$
 avec $\theta$ l'identification canonique  de $T_{h}Met(E_{n})$ avec  $Herm(E_{n},E_{n})$ :

$H(u,v)=h(\theta(H)u,v)$ où $u,v\in E_{n}$.
\item  Si on pose
$$\Phi(H)= l.P(\theta(H),iR(D),\ldots,iR(D))\, ,$$
 $\Phi$ est $\delta$-fermée  modulo
$Im\,\partial\oplus Im\,
\db $ (cf. \cite{koba} chap 6). Donc $\Phi$ est une
$1$-forme sur $Met(E_{n})$ à valeurs dans
$\Omega^{(l-1,l-1)}(X).$
\end{enumerate}
\end{remq}
Alors il existe ${\bf {R}}\,:\, Met(E_{n})\lr
\Omega^{(l-1,l-1)}(X)$ tel que $$
(\delta{\bf{R}})_{h_{0}}(H)=\Phi(H) \quad \forall
h_{0}\in Met(E_{n}) \mbox{ et } H\in
T_{h_{0}}Met(E_{n}).$$ Vu la définition de $\Phi$,
 on a également $i\db\partial
\delta{\bf {R}}=\delta\tilde P.$
On pose
\begin{eqnarray*}
{\bf{R}} \,:\,Met(E_{n})\times Met(E_{n})&\lr&
\Omega^{(l-1,l-1)}(X)\\
(h,k)&\longmapsto&{\bf {R}}(h,k)={\bf {R}}(h)-{\bf
{R}}(k) .
\end{eqnarray*}
Soient  $h,h'\in Met(E_{n})$. On a alors
$i\db\partial {\bf {R}}(h,h')=\tilde P(h)-\tilde
P(h').$ Etudions maintenant les cas où $l=1$ et
$l=2$.
\begin{enumerate}
  \item $l=1$  :
  on note $P_{1}$ l'application $1$-linéaire sur $\mathfrak
{gl}(r,
\mathbb{C})$.
$$P_{1}(A)=tr(A)\quad\forall A  \in \mathfrak
{gl}(r,
\mathbb{C})$$
alors $ \tilde P_{1}(h)=2\pi c_{1}(E_{n},h)$. On
note $\bf{R_{1}}$ l'application correspondant à
$\bf{R}$ dans le cas où $l=1$. On a
$$i\db\partial
{\bf{R_{1}}}(h,h')=2\pi(c_{1}(E_{n},h)-c_{1}(E_{n},h'))$$
et
  $$\delta {\bf{R_{1}}}(H)=\Phi(H)=tr (\theta(H))\, .$$
  C'est à dire $${\bf{R_{1}}}(h',h)=ln \,\det (h ^{-1}h')\,.$$

  \item $l=2$ : (Ceci est donné pour $(X,\omega)$  variété kälhérienne compacte de dimension complexe $m$) \\
on note $P_{2}$ l'application $2$-linéaire sur
$\mathfrak {gl}(r,
\mathbb{C})$,
$$P_{2}(A,B)=tr(A.B)\quad\forall A,B
\in
\mathfrak {gl}(r,
\mathbb{C})\mbox{(forme de Killing)}.$$
  $\tilde P_{2}(h)=4\pi^{2}(c_{1}(E_{n},h)^{2}-2c_{2}(E_{n},h))$. On note $\bf R_{2}$
l'application correspondant à $\bf{R}$ dans le cas
où $l=2$.
   $$\delta {\bf{R_{2}}}(h',h)=2i.tr(h_{t}^{-1}\delta h_{t}.R(h_{t}))$$
où $h_{t}$ est un chemin différentiable de
métriques entre $h$ et $h'$. Nous avons en fait :
  $$i\db\partial {\bf{R_{2}}}(h,h')=4\pi^{2}(\pi_{1}(E_{n},h)-\pi_{1}(E_{n},h'))$$ où $\pi_{1}(E_{n
  },h)$
   est la première classe de Pontryagin de
  $(E_{n},h).$
  On pourra consulter   \cite{koba} pour une autre présentation.

\end{enumerate}

La construction précédente permet de poser la
définition suivante :

\begin{dfn}\index{Fonctionnelle!de S. K  Donaldson}\label{foncd}
Soit $E_{n}$ un fibré vectoriel holomorphe de rang
$r$ sur $(X,\omega)$, variété  kälhérienne
compacte de dimension un. La première
fonctionnelle de S. K. Donaldson est définie par
$$M=\int_{X}{\bf{R_{2}}}+\lambda {\bf R_{1}}\omega ,$$
avec  $\lambda=-4\pi {\mu(E_{n})}$ et $\omega$ est
la forme kälhérienne de $X$.\label{lambda}
\end{dfn}
\begin{prop}\label{r11}(cf. \cite{don3}) Soient $h,h',h''\in Met(E_{n})$, on a :
\begin{enumerate}
  \item $M(h,h')+M(h',h'')+M(h'',h)=0\,.$

  \item $M(h,ah)=0\quad\forall a>0$

  \item  Soient $h(u)$ est une courbe différentiable dans $Met(E_{n})$ et
   $v_{u}=h(u)^{-1}\partial_{u}h(u),$
   $$\du M(h(u),k)=2i\int_{X}tr(v_{u}R(h(u)))+2\lambda.tr(v_{u})\omega.$$
En particulier $\du
\int_{X}{\bf{R_{2}}}(h(u),k)=2i\int_{X}tr(v_{u}R(h(u)))$,
où $k\in Met(E_{n})$ fixé.

\end{enumerate}
\end{prop}
Nous avons alors le théorème suivant :
\begin{thm}(\cite{koba} p .203)
La métrique $h\in Met(E_{n})$ est Einstein-Hermite
si et seulement si elle est un point critique de
la fonctionnelle $M(h,k_{0})$.
\end{thm}

Nous définissons maintenant une fonctionnelle
intermédiaire qui va nous permettre de définir la
deuxième fonctionnelle de S.K Donaldson construite
à partir de $M$.
 \subsection{La fonctionnelle $Norm_{-,\omega}(\beta)$}\index{Notation!$Norm_{h,\omega}(\beta)$}
On peut définir  la fonctionnelle
$Norm_{-,\omega}(\beta)$  :
  \begin{eqnarray*}
  Norm_{-,\omega}(\beta)\, :\, Met(E_{n})&\lr&\R^{*}_{+}\\
  h&\longmapsto&Norm_{h,\omega}(\beta).
  \end{eqnarray*} Si $h,k
\in  Met (E_{n})$  sont telles que $det_{E_{n}}(h)=det_{E_{n}}(k)$, $Norm_{h,\omega}(\beta)$ vérifie (cf. \cite{don2} p. 239 (19)):
\begin{equation}\label{normalpha}
  Norm_{h,\omega}(\beta)=Norm_{k,\omega}(\beta)\times\exp\frac{1}{8\pi}(\int_{X}{\bf R_{2}}(h,k)),
\end{equation}

avec $\beta$ fixé dans le fibré $L$ (cf. p.
 paragraphe \ref{morphi}).
 En fait si  $h,k
\in  Met (E_{n})$, on a
\begin{equation}\label{normalpha1}
  \ln \frac{Norm_{h,\omega}(\beta)} {Norm_{k,\omega}(\beta)}=
  (\frac{\chi(n)}{2r}\int_{X}\ln det_{E_{n}}(k^{-1}h)\omega)+
\frac{1}{8\pi}M(h,k),
\end{equation} où $\chi(n)$  est la caractéristique d'Euler-Poincaré de $E_{n}$. Tout comme $M$
elle permet de déterminer les métriques qui sont
Yang-Mills Hermite. En effet on a la proposition
suivante :
\begin{prop}(cf. \cite{don2})
Soit $E_{n}$ un fibré holomorphe sur $X$. Soit
$k_{0}\in  Met (E_{n})$ fixé. Une métrique de
$E_{n}$ est Yang-Mills Hermite si et seulement si
elle est un extrémum de la fonctionnelle
$Norm_{-,\omega}(\beta)$.
\end{prop}
\index{Fibré!Yang-Mills Hermite}
\subsection{La fonctionnelle $\mathcal  M $}
La formule \ref{normalpha} va nous permettre de
définir la fonctionnelle de S.K Donaldson. En
effet nous avons alors :

\begin{equation}\label{normlog}
  \ln (Norm_{h,\omega}(\beta))=\ln(Norm_{k,\omega}(\beta))+
  \frac{\chi(n)}{2r}\int_{X}\ln det_{E_{n}}(k^{-1}h)\omega+\frac{1}{8\pi}M(h,k).
\end{equation}

\begin{dfn}\label{cestfini}Nous  définissons la
fonctionnelle de S.K. Donaldson par
\begin{eqnarray*}
 \mathcal  M : Met(E_{n})&\longrightarrow&{\mathbb R}\hbox{ avec
}\\
 h &\longmapsto&\ln (Norm_{h,\omega}(\beta))-
  \frac{\chi(n)}{2r}\int_{X}\ln det_{E_{n}}(k_{0}^{-1}h)\omega.
\end{eqnarray*}
\end{dfn}
Nous rappelons que $k_{0}$ est la métrique que
nous avons fixée dans $Met(E_{n})$. \\Nous avons
également :
\begin{thm}(\cite{koba} p .203)
La métrique $h\in Met(E_{n})$ est Einstein-Hermite
si et seulement si elle est un point critique de
la fonctionnelle $\mathcal  M$.
\end{thm}
\begin{remq}\label{mnod}
\begin{enumerate}
%\item On a donc $\forall h\in Met(E_{n})\quad \mathcal  M (h)=
%\ln (Norm_{h,\omega}(\beta))+ \frac{\lambda}{8\pi}\int_{X} R_{1}(h,k_{0})\omega$.
  \item  Dans le paragraphe \ref{prob} nous montrons que :
$\forall h\in Met(E_{n})\, ,$
$$\mathcal{M}(h)={\frac{1}{2}}ldet_{W_{n}}\circ
L_{n}(h)-
\frac{\chi(n)}{2r} kn_{n}(h)+\frac{1}{2}\tau_{h}(X,E_{n}). $$

\end{enumerate}

\end{remq}

\section{Fonctionnelle sur $Met(W_{n})$ : Construction de $\mathcal{KN}_{n}$}
 On
 choisit un isomorphisme
$\Lambda^{dim(W_{n})}(W_{n})\cong
\C$. Ceci revient à définir un déterminant $\det_{W_{n}}$ sur $W_{n}$.\\
 Nous notons $p_{n}$, la
dimension de $W_{n}$. Soit $E_{n}$ un fibré
holomorphe sur $X$.
 \subsection{Rappels}\label{cotal}
On se référera à \cite{kempf} ou \cite{Kirwan} ou
\cite{don3}.

Soit V un espace vectoriel complexe de dimension
finie sur lequel agit un groupe connexe réductif
$G$. Soit $K$ un sous groupe compact de $G$. On
fixe une norme $\Vert\quad\Vert$ sur $V$ de telle
sorte que l'action de $K$ sur $V$ préserve cette
norme.\\ Soit $v\in V$. On
\'etudie la fonction :
$$p_v(g)=\Vert g.v\Vert \quad\forall g\in G\, .$$
Pour mémoire, notons que cette application est
reliée avec la notion d'application moment cf.
\cite{ness}. On a le résultat suivant
:
\begin{thm}\cite{kempf}
Tout point critique de $p_v$ est un point o\`u $p_v$ atteint une
valeur minimale.
Si $p_v$ atteint une valeur minimale alors,
\begin{enumerate}
\item L'ensemble M o\`u $p_v$ atteint cette valeur est une seule classe
$KgG_v\in K\backslash G/G_{v}$ .
\item La variation au second ordre de $p_v$ en un point de M dans
toutes les directions non tangentes \`a M est
définie positive.
\end {enumerate}
\end{thm}
On a également la définition suivante :
\begin{dfn}
$$v\mbox { est semi-stable }\iff 0 \mbox{ n'est pas adhérent à l'orbite }Orb_{G}(v).$$
Si $v$ est semi-stable, on a :
\begin{eqnarray*}
  v \mbox { est stable }\iff&& dim\,{G}(v)=dim\,
Orb_{G}(v)\\ &&Orb_{G}(v) \mbox{ est fermée }
 .
\end{eqnarray*}

\end{dfn}
Prenons $G=SL(p_{n})$ le th\'eor\`eme  suivant
relie la fonction précédente avec la stabilit\'e:
\begin{thm}(cf .\cite{kempf})
\index{Stabilité}
le vecteur $v$ est stable $\iff p_v $ atteint une
valeur minimale.
\end{thm}

\subsection{Construction sur $E_{n}$ de $\mathcal{KN}^{alg}$}\label{ken}
Nous voulons construire  une fonction qui permette
de repérer les fibrés  stables. Pour cela nous
utilisons les résultats donnés le paragraphe
\ref{ke} pour établir le lien entre le fibré
$(E_{n},h)$ et un espace de dimension finie
$V=\Lambda^{2}\C^{p_{n}}$ sur lequel nous pouvons
repérer les points stables  sous l'action d'un
groupe  $G=SL(p_{n})$ grâce à la construction de
la fonction $p_{v}\, ,\, v\in V$.
\\ Nous
supposons que $\deg\,E_{n}\geq d_{0}$ où $d_0$ est
défini comme dans le théorème \ref{psy} p.
\pageref{psy}. Nous posons $W_{n}=H^{0}(X,E_{n})$.
Notons
 $p_{n}=dim\, W_{n}$. Nous avons
également $Met(W_{n})=GL(W_{n})/U(W_{n})$.

 En utilisant le
paragraphe \ref{plongi}, on a alors le morphisme
induit
\label{phix}
:
$$\varphi\, :\,\ox\otimes W_{n} \lr E_{n}\lr 0. $$
avec :
$$\forall x_{i}\in X,\qquad\varphi\,_{x_{i}}
:\,\mathcal{O}_{x_{i}}\otimes W_{n} \lr E_{n,x_{i}}\lr 0.$$
 Nous rappelons que nous nous limitons au cas
où $\rg(E_{n})=2$. \\ Nous pouvons exprimer
$E_{n,x_{i}}$ comme un quotient de dimension deux
de $W_{n}$, c'est à dire comme un élément de
$Grass(2,W_{n})$.
 De plus $SL(p_{n})$ agit sur
$Grass(2,W_{n})$ de façon évidente.
\\Considérons le plongement  de Plücker :
  \begin{eqnarray*}
Grass(2,W_{n})&\lr&\mathbb{P}(\Lambda^{2}W_{n}),\\
F=Vect\{e_{1},e_{2}\}&\longmapsto&[e_{1}\wedge
e_{2}]\mbox{(nous avons identifié $W_{n}^{*}$ à
$W_{n}$)}.
  \end{eqnarray*}

 Ici $E_{n,x_{i}}\in Grass(2,W_{n})$ et
$E_{n},x_{i}$ possède une structure hermitienne
$h_{x_{i}}$. Prenons $\{e_{1},e_{2}\}$ une famille
libre de $\Lambda^{2}W_{n}$ telle que
$\{e_{1},e_{2}\}$ soit une base de    de
$E_{n,x_{i}}$ $h_{x_{i}}$-orthonormée.
 Considérons également une   famille libre $\{v_{1},v_{2}\}$ de
 $\Lambda^{2}W_{n}$ telle que $\{v_{1},v_{2}\}$ soit une base de $E_{n},x_{i}$
$k_{0_{x_{i}}}$-orthonormée. On a alors :
  \begin{eqnarray*}
e_{1}\wedge e_{2} &=&\begin{pmatrix}
<e_{1},v_{1}>_{k_{0_{x_{i}}}}&<e_{2},v_{1}>_{k_{0_{x_{i}}}}\\
<e_{1},v_{2}>_{k_{0_{x_{i}}}}&<e_{2},v_{2}>_{k_{0_{x_{i}}}}\\
\end{pmatrix}v_{1}\wedge v_{2}\\
&=& \begin{pmatrix}
<H_{x_{i}}e_{1},v_{1}>_{h_{x_{i}}}&<H_{x_{i}}e_{2},v_{1}>_{h_{x_{i}}}\\
<H_{x_{i}}e_{1},v_{2}>_{h_{x_{i}}}&<H_{x_{i}}e_{2},v_{2}>_{h_{x_{i}}}\\
\end{pmatrix}v_{1}\wedge v_{2} \quad \mbox{(cf. p.\pageref{metriend})}\\
&=&\det_{E_{n_{x_{i}}}} H_{x_{i}}.det
(P_{e,v})v_{1}\wedge v_{2}\\
&=&\det_{E_{n_{x_{i}}}}k^{-1}_{0_{x_{i}}}h_{x_{i}}.det
(P_{e,v}). v_{1}\wedge v_{2}.
\end{eqnarray*}
Alors
$\det_{E_{n_{x_{i}}}}k^{-1}_{0_{x_{i}}}h_{x_{i}}.det
(P_{e,v}). v_{1}\wedge v_{2}$ correspond au relevé
dans $\Lambda^{2}W_{n}$ de $[e_{1}\wedge e_{2}]\in
\mathbb{P}(\Lambda^{2}W_{n}).$\\
\begin{remq}
Ceci est indépendant du choix des deux bases
$\{e_{1},e_{2}\}$ et $\{v_{1},v_{2}\}$ pourvu
quelles soient $h_{x_{i}}$-orthonormée et
$k_{0_{x_{i}}}$-orthonormée. Remarquons également
que l'action de $g\in Sl(p_{n})$ se traduit par la
modification de $H_{x_{i}}$ en $g.H_{x_{i}}$. Il
existe alors $h^{'}_{x_{i}}$ telle que
$g.H_{x_{i}}=k^{-1}_{0_{x_{i}}}h^{'}_{x_{i}}$
\end{remq}

Maintenant faisons la même construction mais pour
$N$ points.
 Soient  $x_{1},\ldots,x_{N}$, $N$ points,
 sur $X$ avec $N\geq s(\deg\,E_{n})$ où $s$ comme défini dans le chapitre \ref{ke}
 . Considérons les morphismes $\varphi_{x_{i}}$ associés à chacun de ces $N$
points,  nous avons alors :
$$E_{n,x_{1}}\times\ldots\times E_{n,x_{N}}\in
\prod_{\substack{x_{i}\in X\\i=1. . N}}Grass(2,W_{n}). $$
Considérons ensuite les $N$ plongements de Plücker
on a alors
:
$$\prod_{\substack{x_{i}\in X\\i=1. . N}}Grass(2,W_{n})\lr
\prod_{i=1}^{N}(\mathbb{P}\Lambda^{2}{W_{n}}),$$
et $SL(p_{n})$ agit sur chacune de ces variétés.
Nous notons $Z$, la variété
$$\prod_{\substack{x_{i}\in X\\i=1.
. N}}Grass(2,W_{n}),$$comme  nous l'avons fait dans le chapitre \ref {ke}.
 Nous avons remarqué,
au chapitre \ref{ke}, que les points stables (ou
semi-stables) de l'espace  $\mathcal R$ (cf.
définition \ref{rmod} p. \pageref{rmod})
correspondaient aux points stables (ou
semi-stables) de $Z. $ Cette fois-ci le travail de
G. Kempf et L. Ness permet de caractériser ces
points particuliers non pas dans $\mathcal R$ mais
dans $\prod_{i=1}^{N}(\Lambda^{2}W_{n})$. A chaque
point $E_{n,x_{i}}$, $i=1..N$ on peut associer un
élément
$$\det_{E_{n_{x_{i}}}}k^{-1}_{0_{x_{i}}}h_{x_{i}}
. v_{1}\wedge v_{2}\in\Lambda^{2}W_{n}. $$Nous identifions maintenant $W^{n}$
avec l'espace hermitien $\C^{p_{n}}$ muni de sa
base canonique
$\{\varepsilon_{1},\varepsilon_{2},...,\varepsilon_{p_{n}}\}$
en posant $\varepsilon_{1}=v_{1}$ et
$\varepsilon_{2}=v_{2}$. Nous appliquons
maintenant les résultats de G. Kempf et L. Ness.\\
Dans notre situation nous avons $G=SL(p_{n})$, $K$
est tel que $K^{\C}=G$ c'est à dire $K=SU(p_{n})$
et $V=\Lambda^{2}\C^{p_{n}}$.
%%%%%%%%%%%%%%%%%%%%%%%%%%%%%%%%%%%%%%%%%%%%%%
On munit  $V$ de la métrique usuelle $\nor{\quad
}_{V}$ qui
 est $U(p_{n})$-invariante.  De plus
$$ln\prod_{i=1}^{N}\nor{\det_{E_{n_{x_{i}}}}k^{-1}_{0_{x_{i}}}h_{x_{i}}
. v_{1}\wedge v_{2}}_{V}=
\sum_{i=1}^{N}\ln \det_{E_{n_{x_{i}}}}k^{-1}_{0_{x_{i}}}h_{x_{i}}.$$
Nous obtenons alors :
\begin{prop}\label{kalg}
 Soit $\mathcal{KN}^{alg}$ la fonctionnelle de G. Kempf
et L. Ness algébrique  définie  par
:\label{kn3}
$$\forall g\in SL(p_{n}), \quad \mathcal{KN}^{alg}(h,g; x_{1}, x_{2}, . . . x_{N})
=\sum_{i=1}^{N}\ln\det_{E_{n_{x_{i}}}}g.k^{-1}_{0_{x_{i}}}h_{x_{i}}. $$
Elle coïncide avec l'application $p_{v}$ définie
p.\pageref{cotal} avec
$V=\prod_{i=1}^{N}\Lambda^{2}\C^{p_{n}}_{i}$,
$G=SL(p_{n})$.
\end{prop}
\begin{remq}
\begin{enumerate}
  \item
Dans le paragraphe suivant nous définissons une fonctionnelle plus générale
 que nous noterons $\mathcal{KN}$ définie à partire de
 $$\mathcal {KN}_{n}^{\mu}(m) =\frac{1}{2}
 \ln(\det_{W_{n}})(m)-\frac{\chi(n)}{2r} \int_{X}\ln(\det_{E_{n}}\,I_{n}(m))d\mu ,$$
 avec $m\in Met(W_{n})$
et où $\mu$ sera une mesure sur $X$. De plus on
aura:
 $$\mathcal{KN}^{alg}(h,g; x_{1}, x_{2}, . . . x_{N})=
 \mathcal KN^{\delta_{x_{1}}+\ldots+\delta_{x_{N}}}(I_{n}(m)),$$
 avec $\delta_{x_{1}}+\ldots+\delta_{x_{N}}=\mu$ et $I_{n}(m)=h$ .
  \item Le fibré
  $E_{n}$ sera donc stable si $\mathcal{KN}^{alg}(h,g; x_{1}, x_{2}, . . . x_{N})$ atteint une valeur minimale
  lorsque $g$ décrit $Sl(p_{n})$.

\end{enumerate}

\end{remq}

\subsection{Fonctionnelle de G. Kempf et L. Ness généralisée $\mathcal{KN}_{n}$ }

L'énoncé de la   proposition \ref{kalg} nécessite
un nombre fini de points $x_{1},x_{2},...,x_{N}$
avec $N>s(\deg E_{n})$. Nous allons généraliser la
fonctionnelle $\mathcal{KN}^{alg}$,
 en sommant non plus sur un nombre
fini de points mais sur $X$ tout entier. Pour
toute mesure $\mu$ sur $X$ avec $\mu(X)=1$, on
définit :
\begin{eqnarray*}\index{Notation!$kn$}
 kn_{n}^{\mu} :  Met(E_{n}) &\longrightarrow& {\mathbb R}\hbox{
donné par :}\\
 h&\longmapsto& \int_{X}\ln(\det_{E_{n}}\,h)\mu.
 \end{eqnarray*}

 Nous rappelons que nous
avons choisi un isomorphisme
$\Lambda^{dim(W_{n})}(W)\cong
\C$. Nous posons :
\begin{eqnarray*} \index{Notation!$ldet_{W}$}
ldet_{W_{n}} : Met(W_{n})&\lr&  \mathbb {R}\\
 m&\longmapsto&
ln(\det_{W_{n}}m)\, .
\end{eqnarray*}
 Nous avons alors une  fonctionnelle de Kempf et Ness
``généralisée" définie par
:\index{Notation!$\mathcal {KN}_{n}$}
\begin{dfn}\label{kkkk}
\begin{eqnarray*}
\mathcal {KN}_{n}^{\mu} : Met(W_{n})&\longrightarrow&  {\mathbb R}\hbox{ définie par }\\
\mathcal {KN}_{n}^{\mu} &=&{\frac{1}{2}}ldet_{W_{n}}-\frac{\chi(n)}{2r} kn_{n}^{\mu}\circ I_{n}.
\end{eqnarray*}
Si, $d\mu$ est $\omega$ ce qui sera notre cas,
nous la noterons plutôt $\mathcal {KN}_{n}$.
\end{dfn}
La proposition \ref{kalg} nous permet d'obtenir le
lemme suivant :
\begin{lem}
Soit $m\in Met(W_{n})$. Posons $h=I_{n}(m)$. On a
alors :
 $\mathcal{KN}_{n}^{\mu}(m)$
est la limite (en adaptant les constantes) de
$\mathcal{KN}^{alg}(h,g;x_{1},...,x_{N} )$ lorsque
$N\lr
\infty$, pour une distribution convenable des
points $x_{1}, x_{2}, . . . x_{N}$.\\
\end{lem}

Cette fonctionnelle utilise la construction
décrite dans le paragraphe \ref {kn3}. En effet,
puisque nous voulons faire augmenter le degré du
fibré $E_{n}$, nous voulons, en fait, travailler
avec des espaces $Met(W_{n})$ de dimension de plus
en plus grande. De ce fait le nombre de points
pour définir la fonctionnelle $\mathcal{KN}^{alg}$
est lui-aussi de plus en plus grand. C'est
pourquoi nous avons choisi de définir la
fonctionnelle $\mathcal {KN}_{n}$ de cette façon.
Nous avons ainsi les  fonctionnelles qui nous
intéressent.

\chapter{Relations entre les  espaces $Met(E_{n})$ et $Met(W_{n})$ : quelques rappels}

Nous travaillons avec des fibrés $E_{n}$ de rang
$2$ sur $X$ de dimension complexe $1$, mais les
rappels donnés dans ce chapitre sont parfois
donnés dans un cadre plus général. La  relation
entre $Met(E_{n})$ et $Met(W_{n})$ nécessitant la
notion de torsion analytique, nous commençons par
en donner une description.
 On pourra consulter
\cite{atybot}et \cite{don2}.
 \section{Fonction Zêta et Torsion analytique}\label{analyt}
 \label{torsi}
 La torsion analytique est une version analytique de la
torsion de Reidmeister, invariant topologique de
deuxième ordre (elle généralise le volume d'une
transformation linéaire). Ray et Singer l'ont
introduite en 1971 dans \cite{ray1}.Beaucoup de
travaux comme \cite{berline}, \cite{bismut},\cite
{Bisvas},  ont permis de mieux connaître cet
invariant.\\ Reprenons le cadre dans lequel nous
nous situons.
  Soit $X$ une surface de Riemann compacte ou de façon équivalente une
courbe projective  irréductible lisse sur $\mathbb
{C} $. Soit $(E_{n}, h)$ un fibré vectoriel
holomorphe  hermitien. Nous nous intéressons  à
l'étude de la torsion analytique pour ce fibré
$(E_{n}, h).$\\

Soit la suite suivante :
\begin{equation}\label{eqan}
C^0\stackrel {\db}{\longrightarrow} C^1
\longrightarrow 0
\end{equation}
avec $C^i=\Omega_{X}^{0, i}(E_{n} )$. Ici $\db$
correspond
\`a la  partie $(0, 1)$ de la connexion de Chern de $(E_{n},h)$.\\
  La suite pr\'ec\'edente
(\ref{eqan}) d\'efinit un complexe. Notons
$\laplac{h}{q}$ l'opérateur de Laplace ou
Laplacien agissant sur $C^{q}$ défini par :
$$\laplac{h}{q}=\ddetid{h}+\detdid{h}.$$\label{zzzzz}
Nous avons alors :
\begin{dfn}\label{toranalyt}Soit ($E_{n}, h)$ un fibré holomorphe hermitien sur $X$.
La torsion analytique $T_h(X,
E_{n})$\index{Torsion
analytique}\index{Notation!$T_h(X, E_{n})$}
\index{Notation!$\tau_h(X, E_{n}
)$}est d\'efinie par
\label{definit}:
$$\tau_h(X, E_{n}
)=2\, ln \, T_h(X, E_{n}
)=\sum_{q=0}^{dim\,X}(-1)^q q\zeta_{h, q}'(0)
\mbox { avec }$$
$$\zeta_{h, q}(s)={\frac{1}{\Gamma (s)}}\int_0^\infty t^{s-1}
\tr(\explap{h}{q})dt .$$
(on notera également $\zeta_{h, q}$ par
$\zeta_{\Delta_{h,q}}$ pour préciser l'opérateur
de Laplace).
 Ici $P_{n}^\perp=1-P_{n}$ o\`u $P_{n}$
  est la projection orthogonale pour $h$  sur  $\ker\Delta _{h, q}$
.\label{cien}

De plus on   a  posé
$$\tr(\explap{h}{q})=\int_{x\in X}tr(p_{t}(x, x))d(vol(x)),$$
où $p_{t}(x, y)$ correspond au noyau de la chaleur
de $\explap{h}{q}$ (cf. ci-dessous).
\end{dfn}

 Dans la troisième partie nous étudions les
variations de cette torsion analytique lorsque
nous faisons varier la métrique h du fibré
holomorphe hermitien $(E_{n},h)$.
\section{Noyau de la chaleur}\label{nych}
On se rapportera par exemple à \cite{berline}.
L'opérateur
$\laplac{h}{q}=\ddet_{h}+\bar{\partial^*}_{h}\bar{\partial}$
est un opérateur elliptique et l'opérateur
$\explap{h}{q}$ existe et possède un noyau que
l'on notera $p_{t}(x, y)$, avec $x, y\in X$. Ce
noyau s'appelle noyau de la chaleur. On sait de
plus que $\explap{h}{q}$ admet un spectre discret
dont les valeurs propres non nulles sont
positives. On notera $\lambda_{1}$ la première
valeur propre  et $\lambda_{j}\quad j\geq 2$ les
autres. A ces valeurs propres, on peut associer
une base hilbertienne, formée de sections propres
correspondantes que l'on note $\psi_{j}\in
\cin(X, E_{n})$. Le noyau de la chaleur est donné par la formule
suivante :
$$p_{t}(x, y)=\sum_{j\geq 1}exp(-t\lambda_{j})\psi_{j}(x)\otimes \psi^{*}_{j}(y).$$
\begin{remq}(cf. \cite{berline})
Le noyau de la chaleur a les propriétés suivantes
:
\begin{enumerate}
\item $p_{t}(x, y)\in\cin(]0, \infty[\times X\times X, Hom(E, E))$,
\item$(\dt + \laplac{h}{q})p_{t}(x, y)=0$,
\item$p_{t}(x, y)\lr \delta_{y}$ si $t\lr 0$,
\item $p_{t}(x, y)=p^{*}_{t}(y, x)$.
\end{enumerate}
\end{remq}
 Nous utiliserons par la
suite pour $t$ proche de zéro un développement
asymptotique du noyau de la chaleur. Si nécessaire
nous indiquerons l'opérateur en question.

\section{Fibr\'es en droite particuliers sur $\espco$ : métrique de Quillen}\label{morphi}
Soit $E_{n}$ un fibré sur $X$  de degré zéro. On
consultera \cite{don3} et \cite{don2}.\\ Sur
$\espco $, on considère la forme définie (cf.
\pageref{metria}) par
:
$$\omega_{\espco}(a,b)=\frac{1}{4\pi^{2}}\int_{X}tr(a\wedge b).$$ Soit  $D\in\espco$, on
considère le fibré en droite sur $\espco$ défini
par :
$$L_{D}=\Lambda^{max}(Ker\db_{D})^{*}\otimes \Lambda^{max}(Ker\db^{*}_{D}).$$
L est un fibré holomorphe appelé fibré déterminant
(cf. \cite{quillen}). Nous allons maintenant
expliciter une métrique $\mathfrak h$ hermitienne
bien particulière définie par D. G. Quillen (cf.
\cite{quillen}).
 Pour cela considérons un autre
fibré défini sur $\espco$ que l'on notera
$\mathcal L$ définit comme suit :\\ Soit
$\alpha>0$, soit $U_{\alpha}$ l'ouvert de $\espco$
ne contenant que des connections $D$ telles que
$\alpha $ ne soit pas valeur propre de l'opérateur
de Laplace $\db_{D}\db_{D}^{*}$. Au dessus de
$U_{\alpha}$ nous avons le sous espace propre des
valeurs propres de $\db_{D}^{*}\db_{D}$
correspondant au valeurs propres $\lambda<\alpha$
noté $\mathcal {H}^{+}_{\alpha}$ et le sous espace
propre de $\db_{D}\db_{D}^{*}$ défini de même pour
$\alpha>\lambda$ et noté $\mathcal
{H}^{-}_{\alpha}$. Sur $U_{\alpha }$ on définit
donc
:
$$\mathcal {L}=\Lambda^{max}\mathcal{H}^{+}_{\alpha}\otimes \Lambda^{max}(\mathcal {H}
^{-}_{\alpha})^{*}.$$
  Sur
$U_{\alpha}$, $\mathcal {L}_{\vert_{U_{\alpha}}}$
est muni d'une métrique $L^{2}$ naturelle, notée
$\vert
\,.\,\vert_{\alpha }$. Pour ``recoller" ces métriques on
les modifie  en les multipliant par le déterminant
régularisé de l'opérateur de Laplace
\mbox{c.-à-d.} $\Pi_{\sigma}
\lambda_{\sigma}=\exp(-\zeta^{'}_{\Delta_{D}}(0))$ ou $\lambda_{\sigma}$
désigne une valeur propre de l'opérateur de
Laplace
$\Delta_{D}=\db_{D}\db_{D}^{*}+\db_{D}^{*}\db_{D}$(cf.
paragraphe \ref{torsi}). Il existe un isomorphisme
canonique entre $L$ et $\mathcal {L}^{*} $(cf.
\cite{quillen}). On définit alors
:

\begin{equation}\label{quil}
  \mathfrak {h}(., .)=\vert .,.
 \vert_{\alpha }. exp(-\zeta^{'}_{\Delta_{.}}(0))\,  .
\end{equation}

C'est la métrique de Quillen du fibré $L$. Elle a
pour courbure $-2i\pi\omega_{\espco}$(cf. \cite
{quillen} ). Cette métrique reste valable si le
degré du fibré $E_{n}$ est quelconque.

L'isomorphisme entre $L$ et $\mathcal L^{*}$ et en
particulier la formule \ref{quil} nous permettent
de munir le fibr\'e holomorphe $L$ d'une norme
faisant intervenir la torsion analytique . En
effet, on a  le r\'esultat suivant :
\begin{thm}\label{norma}(\cite{don2}) Soit $E_{n}$ un fibré  sur $X$. %Fixons le point $D_{0}\in \espco
%$ qui correspond à la métrique $k_{0}$,
Il existe $n_{0}$ tel que pour $n \geq n_{0}$ on
ait :\\
 $\forall h\in Met(E_{n})$
\begin{equation}\label{normap}
Norm_{h, \omega}(\beta )=
(\det_{W_{n}}(L_{n}(h)))^{{\frac{1}{2}}}\times
T_{h}(X, E_{n})
\end{equation}
avec $T_{h}(X, E_{n})$ la torsion analytique qui
est d\'efinie au paragraphe \ref{analyt}.
\end{thm}

 Cette
norme sur $L$ nous permet de d\'etecter les
m\'etriques hermitiennes Yang-Mills en regardant
les extr\'emas de la fonctionnelle suivante :
$$h\longrightarrow Norm_{h, \omega}(\beta)$$
avec $\beta $ un \'el\'ement du fibr\'e en droite
L. Ce résultat est le premier pont que l'on peut
faire entre $Met(E_{n})$ et $Met(W_{n})$.  Dans
notre situation par l'étude du terme $ T_{p}(X,
E_{n})$ et grâce aux deux morphismes $L_{n}$ et
$I_{n}$, nous obtiendrons divers résultats résumés
dans le chapitre suivant.

\chapter {Relation entre les fonctionnelles}
Nous rappelons que $Met(E_{n})$ ne dépend pas de
$n$ (cf. \pageref{indep}).

Définissons la  variété banachique
$Met(E_{n},C^{4})$ : nous fixons une métrique
lisse de référence notée
 $k_{0}$, l'espace $Met(E_{n},C^{4})$ est l'espace des métriques $h$,
 telles que $h=k_{0}H$ avec $H\in C^{4}(End(E_{n})).$
Sur $Met(E_{n},C^{4})$ on définit alors  une
distance notée $d(.,.)_{L^{2}_{1}(Met(E_{n}))}$
définie par :$$\forall h,k\in Met(E_{n},C^{4})
\quad d(h,k)_{L^{2}_{1}(Met(E_{n}))}=d(k_{0}^{-1}h,k_{0}^{-1}k)_{L^{2}_{1}(End(E))}.$$
\section{Enonc\'es des résultats}
On a  le résultat suivant (pour les notations cf.
p. \pageref{kkkk} et p. \pageref{cestfini})
:
\begin{thm}\label{th1}
 Soient $h_{0},h_{1}\in Met(E_{n}).$  On consid\`ere
le chemin $$h(u)=uh_{1}+(1-u)h_{0}$$ avec $u
\in
\lbrack0, 1\rbrack$. Alors :

$$\forall \varepsilon >0 \quad \exists n_{0},\quad \forall n\geq n_{0},$$

$$\vert(\du(\mathcal{M}+\frac{\chi(n)}{2r}  kn_{n}-{\frac{1}{2}}{ldet_{W_{n}}}\circ L_{n} ))(h(u))\vert\leq \varepsilon.$$
\end{thm}
\begin{remq}
Si $det(h_{u})=1$, alors $
\frac{\chi(n)}{2r} kn_{n}=0$, d'où le théorème donne avec les mêmes hypothèses :
$$\vert(\du(\mathcal{M}-{\frac{1}{2}}{ldet_{W_{n}}}\circ L_{n} ))(h(u))\vert\leq \varepsilon.$$

\end{remq}
\label{prob}
La preuve de ce théorème est donnée p.
\pageref{th1}. Le second théorème va relier en fait
les fonctionnelles $\mathcal  M$ et
${\mathcal{KN}_{n}}.$ Pour un $n$ donné, on peut
écrire la  fonctionnelle $\mathcal{M}$  sous la
forme
:\index{Notation!$\mathcal  M$}
\begin{lem}\label{lem12}$$\mathcal  M(h)={\frac{1}{2}}ldet_{W_{n}}\circ
L_{n}(h)-
\frac{\chi(n)}{2r} kn_{n}(h)+\frac{1}{2}\tau_{h}(X,E_{n}).$$
\end{lem}
\begin{demo}
En effet d'apr\`es le théorème \ref{norma} (p.
\pageref{norma}) on a  :

\begin{eqnarray}\label{fff}
Norm_{h, \omega}(\beta )&=&(\det_{W_{n}}
L_{n}(h))^{\frac{1}{2}}\times T_{h}(X,
E_{n})\nonumber\\
\ln(Norm_{h,\omega}(\beta ))&=&\ln((\det_{W_{n}} L_{n}(h))^{\frac{1}{2}}\times T_{h}(X, E_{n}))\nonumber\\
\ln (Norm_{h, \omega}(\beta))&=&{\frac{1}{2}}ldet_{W_{n}}\circ L_{n}+\ln(T_{h}(X, E_{n}))\nonumber\\
 \ln (Norm_{h, \omega}(\beta))&=&{\frac{1}{2}}ldet_{W_{n}}\circ L_{n}+\frac{1}{2}\tau_{h}(X,E_{n}),
\end{eqnarray}
par définition de $ldet_{W_{n}}$ et par définition
de la torsion analytique p.\pageref{toranalyt}.\\
Vu la définition p. \pageref{mnod}  de $\mathcal
{M}$ et en utilisant l'équation \ref{fff} , on a
:$$\mathcal {M}(h)={\frac{1}{2}}ldet_{W_{n}}\circ
L_{n}(h)-
\frac{\chi(n)}{2r} \int_{X}\ln
\det_{E_{n}} (k_{0}^{-1}h)\omega+\frac{1}{2}\tau_{h}(X,E_{n}).$$ En notant que
$$\int_{X}\ln
\det_{E_{n}} (k_{0}^{-1}h)\omega= kn_{n}(h)$$(cf. p.
\pageref{kkkk}), on a donc le résultat souhaité.
  \end{demo}
  Les deux fonctionnelles, $\mathcal  M$ et ${\mathcal{KN}_{n}}$ vérifient les propri\'et\'es suivantes :
\begin{enumerate}
  \item les points où $\mathcal  M$ atteint son minimum correspond
  aux métriques du fibré holomorphe E dont la connexion de Chern est Yang-
  Mills.

  \item
  Si $E_{n}$ est tel que ${\mathcal{KN}_{n}}$ atteint un extrémum alors $E_{n}$ est stable.

\end{enumerate}
Nous voudrions ``comparer" les points où les
extréma de  ${\mathcal{KN}_{n}}$ et $\mathcal  M$
sont  atteints. Soit $n\in \N$, soit $h\in
Met(E_{n})$ alors :
\label{forul}
$$\mathcal{KN}_{n}\circ L_{n}(h)=-\frac{\chi(n)}{2r} kn_{n}\circ I_{n}\circ L_{n}(h)
+{\frac{1}{2}}{ldet_{W_{n}}}\circ L_{n}(h)\mbox{
(d'après la déf. \ref{kkkk} p.
\pageref{kkkk})}$$
et
$$\mathcal{M}(h)={\frac{1}{2}}ldet_{W_{n}}\circ
L_{n}(h)-
\frac{\chi(n)}{2r} kn_{n}(h)+\frac{1}{2}\tau_{h}(X,E_{n})$$

(d'après le lemme \ref{lem12} ci dessus). On a
donc
:
$$\mathcal{KN}_{n}\circ L_{n}(h)-\mathcal{M}(h)=-
 \frac{\chi(n)}{2r} kn_{n}\circ I_{n}\circ L_{n}(h)+\frac{\chi(n)}{2r}
kn_{n}(h) -\frac{1}{2}\tau_{h}(X,E_{n}).$$ On
utilise la définition de $kn_{n}$ (cf.
p.\pageref{kkkk}) et on obtient :
\begin{equation}\label{forul1}
(\mathcal{KN}_{n}\circ L_{n}-\mathcal{M})(h) =
-\frac{\chi(n)}{2r} \int_{X}\ln \det_{E_{n}}(h^{-1}I_{n}L_{n}(h))\, \omega-\frac{1}{2}\tau_{h}(X,E_{n}).
\end{equation}

  Le théorème suivant relie $\mathcal{KN}_{n}\circ
L_{n}$ et $\mathcal{M}$.
\begin{thm}\label{princ}
$\quad \forall h,k\in
Met(E_{n},C^{4})\,\quad\exists C\in \R^{+}$
 $$\vert \mathcal{M}(h)-\mathcal{KN}_n\circ L_{n} (h)-
 \mathcal{M}(k)+\mathcal{KN}_n\circ L_{n}(k)\vert\leq \frac{C}{n}.$$
\end{thm}
Considérons maintenant $U$ un compact convexe dans
$Met(E_{n},C^{4})$. On en déduit alors, pour $n$
assez grand que  si $m$ est un extrémum de
${\mathcal{KN}_{n}}\vert_{ L_{n}(U)}$ alors $ m$
est approximativement l'image d'un minimum local
de $\mathcal M$. En particulier, si le minimum de
Yang-Mills est contenu dans  $U$ alors  pour $n $
assez grand, le minimum de $\mathcal
{KN}_{n}\vert_{ L_{n}(U)}$ correspond
approximativement au minimum de Yang-Mills. Plus
précisément on a le corollaire suivant
:
\begin{thm}\label{top}
Soit  $E_{n}$ un fibré stable sur $X$. Soit $U$ un
compact convexe donné de $Met(E,C^{4})$. Supposons
que $h_{YM}$, la métrique de Yang-Mills(elle
existe car $E_{n}$ est stable), appartienne à $U$.
On a alors
:
$$\forall \varepsilon >0 \quad \exists n_{U},\quad \forall n\geq n_{U},$$
il existe  m  extrémum de
${\mathcal{KN}_{n}}{\vert}_{ L_{n}(U)}$ et
$$d(h_{YM},cI_{n}(m))_{L^{2}_{1}(Met(E_{n}))}\leq \varepsilon, $$
avec $c=\frac{n}{\pi}(1+O(1)).$
\end{thm}
\section{Résultats  de la troisième partie utilisés }
Nous donnons ici de façon synthétique les deux
résultats qui seront démontrés dans la troisième
partie.

\subsection{Résultat 1 (cf. th. \ref{torsiuon} )}
\begin{thm}\label{ao} Soient $h_{0},h_{1}\in
Met(E_{n}).$  On consid\`ere le chemin
$$h(u)=uh_{1}+(1-u)h_{0}$$ avec $u
\in
\lbrack0, 1\rbrack$. On considère la famille
$(E_{n},h(u))_{n\in N}$ de fibrés holomorphes
hermitiens sur X. Alors
:

$${\du}\tau_{h(u)}(X, E_{n})
=O(\negli{}).$$

\end{thm}
\begin{remq} On peut remarquer notamment que si on se place sur un
compact convexe $U$ de $Met(E_{n},C^{4})$, alors
ce résultat est uniforme en $h_{0}\in U$ et
$h_{1}\in U$.
\end{remq}
\subsection{Résultat 2 (cf. th. \ref{iln})}
Nous notons  par la suite pour $a\in \N$,
 $\bf{O(\negli{a})}$ les
matrices dont les coefficients appartiennent à
$\cin(X,E_{n})$ et sont comparables, en norme
$C^{0}$, à $\negli{a}$.
\begin{thm}\label{iln2} (cf. pour les notations  p. \pageref{courbb})\\
Soit $(E_{n},h)$ un fibré hermitien sur $X$ on a :
$$h^{-1}I_{n}L_{n}(h)=
\frac{\pi}{n^{2}}(R_{\C}(\ox(n)Id_{E_{0}}
-R_{\C,h}(E_{0})Id_{\ox(n)}+r Id_{E_{n}})+{\bf {O}(\negli{3})}.$$
Ici $r$ est une fonction ne dépendant que de
$X$(cf. formule \ref{vol} p.
\pageref{vol}).\end{thm}
\section{Exemple}
L' exemple choisi illustre d'une part le résultat
du théorème \ref{ao}, et d'autre part permet un
calcul de $\mathcal M$. Dans ce cas particulier,
la fonctionnelle $\mathcal M$ est linéaire et donc
ne présente pas d'extrémum car le fibré choisi est
instable.\\
 Soit $\P^{1}$ muni de la métrique
de Fubini Study. Si on se donne le fibré  suivant
:
$$E_{0}=\mathcal{O}_{\P^{1}}(1)\oplus
\mathcal{O}_{\P^{1}}(-1),$$
 On a alors
$$E_{n}=\mathcal{O}_{\P^{1}}(1+n)\oplus
\mathcal{O}_{\P^{1}}(n-1).$$ On considère $h(u)$ un chemin de
métriques de $E_{0}$ de la forme $\exp(u)h_{\mathcal{O}_{\P^{1}}(1)}\oplus
\exp(-u)h_{\mathcal{O}_{\P^{1}}(-1)}$, alors
$\tau_{h(u)}(X, E_{n})$ est indépendant de $u$,
car

\begin{enumerate}
  \item $\Delta_{E_{0}}=\Delta_{\mathcal{O}_{\P^{1}}(1)
  }\oplus\Delta_{\mathcal{O}_{\P^{1}}(-1)}$, donc
les valeurs propres  de $\Delta_{E_{0}}$ sont
indépendantes de $u$ de même pour
$\Delta_{E_{n}}$, et de même $R_{h(u)}(E_{n})$.
\item  $\tr\alpha=0$. Ici $\alpha =h_{u}^{-1}\partial_{u}h_{u}=\begin{pmatrix} 1&0\\ 0&-1\\
\end{pmatrix}$.
\end{enumerate}
On obtient de plus que :
$$\du \mathcal {M}(h(u))=\frac{2i}{8\pi}\int_{X}tr(\alpha
R_{h(u)}E_{n})=
c_{1}(\mathcal{O}_{\P^{1}}(1))=1.$$
 Donc  $\mathcal {M}$ est linéaire (on
constate d'ailleurs l'instabilité de $E_{0}$). De
même on remarque que si on a écrit
$h(u)=\exp(u)h_{\mathcal{O}_{\P^{1}}(1)}\oplus
\exp(-u)h_{\mathcal{O}_{\P^{1}}(-1)}$
$$ ldet_{W_{n}}(L_{n}(h(u)))=u.dim
H^{0}(X,E_{n+1})-u.dim H^{0}(X,E_{n-1}).$$ Cette
fonctionnelle ne présente pas non plus d'extrémum
car elle est linéaire. De plus
$$\du ldet_{W_{n}}(L_{n}(h(u)))=2.$$
Donc cet exemple illustre le théorème \ref{ao}, de
plus il permet de vérifier les signes utilisés.
\section{ Preuves }\label{squel}
 \subsection{Preuve du théorème \ref{th1}}
  La
d\'emonstration du théorème \ref{th1} découle du
théorème \ref{ao} et du lemme \ref{lem12}.\\
 Soient $h_{0},h_{1}\in
Met(E_{n}).$  On consid\`ere le chemin
$h(u)=uh_{1}+(1-u)h_{0}$ avec $u
\in
\lbrack0, 1\rbrack$.

$$\mathcal  M (h(u))={\frac{1}{2}}ldet_{W_{n}}\circ L_{n}(h(u)) -
\frac{\chi(n)}{2r}  kn_{n}(h(u))+\frac{1}{2}\tau_{h(u)}(X,E_{n}).$$

Donc
$$\du[\mathcal  M (h(u))-{\frac{1}{2}}\du ldet_{W_{n}}\circ
L_{n}(h(u))+\frac{\chi(n)}{2r}
kn_{n}(h(u))]=\frac{1}{2}\du\tau_{h(u)}(X,
E_{n}).$$

Nous utilisons alors  le théorème \ref{ao} qui
nous  donne le premier résultat (théorème
\ref{th1}).
\subsection{Preuve du théorème \ref{princ}}\label{dprinc}
Soit $n\in \N$, nous avons vu la formule
\ref{forul1}:
 $$ (\mathcal{KN}_{n}\circ
L_{n}-\mathcal{M})(h)=-\frac{\chi(n)}{2r}
_{1}\int_{X}\ln
\det_{E_{n}}(h^{-1}I_{n}L_{n}(h))\,
\omega-\frac{1}{2}\tau_{h}(X,E_{n}).$$ Le théorème \ref{iln2} permet d'écrire
:
\begin{eqnarray}\label{eqqq}&&\frac{\chi(n)}{2r}\int_{X}\ln
\det_{E_n}(h^{-1}I_{n}L_{n}h)\, \omega
=\nonumber \\ &&\frac{\chi(n)}{2r}\int_{X}\ln \det_{E_n}
(\frac{\pi}{n^{2}}(R_{\C}(\ox(n))Id_{E_{0}}-R_{\C,h}(E_{0})Id_{\ox(n)}+r
Id_{E_{n}})\nonumber\\&+&{\bf {O}(\negli{3})} )
\omega.
\end{eqnarray} On rappelle que le
fibré en droite $\ox(n)$ est stable (cf. prop.
\ref{stabil} p. \pageref{stabil}). Donc
$R_{\C}(\ox(n))= n.Id_{\ox(n)}$. On a donc
:

$$\ref{eqqq}=\frac{\chi(n)}{2r}\int_{X}\ln
\det_{E_n} ((\frac{\pi}{n^{2}}R_{\C}(\ox(n)Id_{E_{0}})(Id_{E_{n}}-
\frac{R_{\C,h}(E_{0})Id_{\ox(n)}}{n}+\frac{r Id_{E_{n}}}{n}+{\bf
{O}(\negli{2})})\omega.$$ Cela donne :

  \begin{eqnarray*}
  \ref{eqqq}&=&\frac{\chi(n)}{2r}\int_{X}\ln
\det_{E_n} (\frac{\pi}{n^{2}}R_{\C}(\ox(n))Id_{E_{0}})\omega\\
&+&\frac{\chi(n)}{2r}\int_{X}\ln
\det_{E_n} (Id_{E_{n}}-\frac{R_{\C,h}(E_{0})Id_{\ox(n)}}{n}+\frac{r Id_{E_{n}}}{n}+{\bf
{O}(\negli{2})})\omega.
  \end{eqnarray*}
On a alors :
\begin{equation}\label{inty}
  \ref{eqqq}= \frac{\chi(n)}{2r}\int_{X}\ln
\det_{E_n}((1+\frac{r}{n})Id_{E_{n}}-\frac{R_{\C,h}(E_{0})Id_{\ox(n)}}{n}+{\bf
{O}(\negli{2})})\omega .
\end{equation}
 L'équation \ref{inty} donne :
  \begin{eqnarray*}
 \frac{\chi(n)}{2r}\int_{X}\ln(1&+&tr
\frac{R_{\C,h}(E_{0})}{(1+\frac{r}{n})n}\\
&&+\frac{1}{4}\big{(}(tr
\frac{R_{\C,h}(E_{0})}{(1+\frac{r}{n})n})^{2}
-tr(\frac{R_{\C,h}(E_{0})}{(1+\frac{r}{n})n}\wedge\frac{\C,R_{h}(E_{0})}{(1+\frac{r}{n})n})+{{o}(\negli{})}\big{)}\omega\,
  \end{eqnarray*}

 et nous avons alors :
$$\ref{inty}=\frac{\chi(n)}{2r}\int_{X}tr
\frac{R_{\C,h}(E_{0})}{(1+\frac{r}{n})n}+\frac{1}{4}\big{(}
(tr \frac{R_{\C,h}(E_{0})}{(1+\frac{r}{n})n})^{2}
-tr(\frac{R_{\C,h}(E_{0})}{(1+\frac{r}{n})n}\wedge
\frac{R_{\C,h}(E_{0})}{(1+\frac{r}{n})n})\big{)}+{{O}(\negli{2})}\big{)}\omega\,
.$$
Nous avons :
  $$\frac{\chi(n)}{2r}\int_{X}tr\frac{1}{((1+\frac{r}{n})n)^{2}}(R_{\C,h}(E_{0})\wedge R_{\C,h}(E_{0}))\omega=O(\negli{})$$
  et $$\frac{\chi(n)}{2r}\int_{X}\
\big{(}tr
\frac{R_{\C,h}(E_{0})}{(1+\frac{r}{n})n}\big{)}^{2}\omega=O(\negli{}).$$
Il s'ensuit que
:\begin{eqnarray}&&\frac{\chi(n)}{2r}\int_{X}\ln
\det_{E_n}(h^{-1}I_{n}L_{n}h)\,\omega\nonumber \\ &=&\frac{\chi(n)}{2r}\int_{X}\ln
\det_{E_n} (\frac{\pi}{n^{2}}R_{\C}(\ox(n))Id_{E_{0}})\omega+\frac{\chi(n)}{2r}\int_{X}tr
\frac{R_{\C,h}(E_{0})}{n}\omega+O(\negli{}).\nonumber
\end{eqnarray}
Nous obtenons de même pour $k\in Met(E_{n})$
\begin{eqnarray}&&\frac{\chi(n)}{2r}\int_{X}\ln
\det_{E_n}(k^{-1}I_{n}L_{n}(k))\, \omega=\nonumber \\ &=&\frac{\chi(n)}{2r}\int_{X}\ln
\det_{E_n} (\frac{\pi}{n^{2}}R_{\C}(\ox(n))Id_{E_{0}})\omega+\frac{\chi(n)}{2r}\int_{X}tr
\frac{R_{\C,k}(E_{0})}{n}\omega+O(\negli{}).\nonumber
\end{eqnarray}

Notons que
 $R_{\C}(\ox(n))$ est indépendant du choix de la métrique sur $E_{n}$.
D'où finalement :
  \begin{eqnarray*}
   (\mathcal{KN}_{n}\circ
L_{n}-\mathcal{M})(h)&-&( (\mathcal{KN}_{n}\circ
L_{n}-\mathcal{M})(k)\\&=&-\frac{\chi(n)}{2r}
_{1}\int_{X}\ln
\det_{E_{n}}(h^{-1}I_{n}L_{n}(h))\,
\omega-\frac{1}{2}\tau_{h}(X,E_{n}))\\
&-&(-\frac{\chi(n)}{2r}
_{1}\int_{X}\ln
\det_{E_{n}}(h^{-1}I_{n}L_{n}(h))\,
\omega-\frac{1}{2}\tau_{h}(X,E_{n})))\\ &=&-\frac{\chi(n)}{2r}\int_{X}tr
\frac{R_{\C,h}(E_{0})}{n}\omega+\frac{\chi(n)}{2r}\int_{X}tr
\frac{R_{\C,k}(E_{0})}{n}\omega\nonumber\\&+&\frac{1}{2}\tau_{k}(X,E_{n})
-\frac{1}{2}\tau_{h}(X,E_{n})+O(\negli{})\nonumber\\
&=&-\frac{\frac{\chi(n)}{2r}}{n}c_{1}(E_{0})+
\frac{\frac{\chi(n)}{2r}}{n}c_{1}(E_{0})+\frac{1}{2}\tau_{k}(X,E_{n})
-\frac{1}{2}\tau_{h}(X,E_{n})+O(\negli{})\\
&=&\frac{1}{2}\tau_{k}(X,E_{n})
-\frac{1}{2}\tau_{h}(X,E_{n})+O(\negli{}).
  \end{eqnarray*}
Or le théorème \ref{ao} donne :
$$\du\tau_{h(u)}(X,E_{n})=O(\negli{}).$$
Le théorème en découle alors.
\subsection{Preuve du théorème \ref{top}}
Nous utilisons le lemme suivant :
\begin{lem}\label{vfvf}Soit $E_{n}$ un fibré stable sur $X$. Soit $U$ un compact de $Met(E_{n},C^{4})$.
Alors il existe des constantes  $C_{1},C_{2}$
telles que : $\forall h\in U$
$$C_{1}.\mathcal{M}(h)\leq d(h,h_{YM})^{2}_{L^{2}_{1}(Met(E_{n}))}\leq C_{2}.\mathcal{M}(h)$$
\end{lem}
\begin{demo}Le point de référence est ici $h_{YM}$. Nous pouvons supposer que
$\det_{E_{n}}(h)=\det_{E_{n}}(h_{YM})$. Posons
$h=h_{YM}e^{s}$, où $s$ est une section de trace
nulle de $End(E)$. Nous pouvons remarquer que les
valeurs propres de $s$ sont bornées en normes
$C^{0}$ car $U$ est un compact de
$Met(E_{n},C^{4})$. L'inégalité de Sobolev donne
pour $s$ section à trace nulle ($\nabla_{h}$ est
défini p.\pageref{deltah}):
\begin{equation}\label{sobo}
\nor{s}_{L^{2}(End(E_{n}))}\leq
\nor{\nabla_{h}s}_{L^{1}}.
\end{equation}
D'après \cite{don2} p. 242 on a
: $\exists\, C_{3 }>0$ tel que
$$\nor{\nabla_{h}s}^{2}_{L^{2}}\leq
C_{3}.\mathcal{M}(h).$$ En combinant ce résultat
et l'inégalité \ref{sobo} on obtient : $\exists\,
C_{2 }>0$ tel que
$$d(h,h_{YM})^{2}_{L^{2}_{1}(Met(E_{n}))}\leq
C_{2}.\mathcal{M}(h).$$
 De plus, comme
$U$ est compact de $Met(E_{n},C^{4})$, les valeurs
propres de $h_{YM}^{-1}h$ sont bornés en norme
$C^{0}$. On a également : $\exists\, C_{4 }>0$ tel
que
$$\mathcal M(h)\leq C_{4} \nor{\nabla s}_{L^{2}}^{2}\mbox{ cf. \cite{don2} p. 242 }.$$
d'où : $\exists\, C_{1 }>0$ tel que
$$\mathcal M(h)\leq C_{1}. d(h,h_{YM})^{2}_{L^{2}_{1}(Met(E_{n}))}.$$

\end{demo}
Reprenons la démonstration du théorème \ref{top}.
 Comme $E_{n}$ est stable, la
fonctionnelle de Donaldson $\mathcal M$ admet une
borne inférieure qui est l'image du point unique
$h_{YM}$ indépendant de $n$. Comme $U$ est compact
$\mathcal{KN}_n\circ L_{n} $ admet un minimum sur
$U$ noté $h_{KNn}$. On considère alors la suite
$(h_{KNn})_{n\in
\N}$. Soit $k\in U$. Posons $A=(\mathcal{M}-\mathcal{KN}_n\circ L_{n} )(k)$.
On utilise maintenant le théorème \ref{princ}
. On a alors :
$$\exists C>0
\, ,\quad \forall h\in U$$
\begin{equation}\label{rrrt}
  \vert(\mathcal{M}-\mathcal{KN}_n\circ L_{n} )(h)-A\vert \le \frac{C}{n}.
\end{equation}

On obtient donc : $\forall \, n\in N$
  \begin{eqnarray*}
  (\mathcal{KN}_n\circ L_{n} )(h_{YM})&-&(\mathcal{KN}_n\circ L_{n} )(h_{KNn})
  =(\mathcal{KN}_n\circ L_{n} )(h_{YM})+A-\mathcal{M}(h_{YM})\\
  &+&\mathcal{M}(h_{KNn})-
  (\mathcal{KN}_n\circ L_{n} )(h_{KNn})-A+\mathcal{M}(h_{YM})-\mathcal{M}(h_{KNn}).
  \end{eqnarray*}
Or $h_{YM}$ est un extrémum de $\mathcal{M}$ donc
:
$$\mathcal{M}(h_{YM})-\mathcal{M}(h_{KNn})\leq 0.$$
D'où en utilisant l'inégalité \ref{rrrt} :
$\forall \, n\in N$
$$\vert(\mathcal{KN}_n\circ L_{n} )(h_{YM})
-(\mathcal{KN}_n\circ L_{n} )(h_{KNn})\vert\leq \frac{2c}{n}$$
Si on note $h_{KN}$ la limite de la suite(ou d'une
sous-suite convergente) $(h_{KNn})_{n\in N}$, On
obtient  :
$$\mathcal{M}(h_{KN})=\mathcal{M}(h_{YM})=0.$$
D'après le lemme \ref{vfvf}, on a donc :
$$d(h_{YM},h_{KN_{n}})_{L^{2}_{1}(Met(E_{n}))}\lr 0 \mbox{ quand } n\lr \infty.$$
De plus le théorème \ref{iln2} donne :
$$h_{KN}^{-1}I_{n}L_{n}(h_{KN})=
\frac{\pi}{n^{2}}(R_{\C}(\ox(n)Id_{E_{0}}
-R_{\C,h_{KN}}(E_{0})Id_{\ox(n)})+{\bf {O}(\negli{3})}$$
D'où si on pose $L_{n}(h_{KN})=m_{KN_{n}}$, on
obtient
:
$$h_{YM}^{-1}I_{n}(m_{KN_{n}})=\frac{\pi}{n^{2}}(R_{\C}(\ox(n)Id_{E_{0}}
-R_{\C,h_{KN}}(E_{0})Id_{\ox(n)})+{\bf {O}(\negli{3})}.$$
On définit alors $c=\frac{n}{\pi}(1+O(1))$ tel que
:
\\$\forall
\varepsilon\quad \exists n_{U}\quad \forall n>n_{U}\,
 ,$$$d(h_{YM},c.I_{n}(m_{KN_{n}}))_{L^{2}_{1}(Met(E_{n}))}\leq \varepsilon.$$
 Ceci achève la démonstration.

%%%%%%%%%%%%%%%%%%%%%%%%%%%%%%%%%%%%%%%%%%%%%%%%%%%%%%%%%%%%%%%%%%%%%%%%
%%%%%%%%%%%%%%%%%%%%%%%%%%%%%%%%%%%%%%%%%%%%%%%%%%%%%%%%%%%%%%%%%%%%%%%
%%%%%%%%%%%%%%%%%%%%%%%%%%%%%%%%%%%%%%%%%%%%%%%%%%%%%%%%%%%%%%%%%%%%%%%%
%%%%%%%%%%%%%%%%%%%%%%%%%%%%%%%%%%%%%%%%%%%%%%%%%%%%%%%%%%%%%%%%%%%%%
%%%%%%%%%%%%%%%%%%%%%%%%%%%%%%%%%%%%%%%%%%%%%%%%%%%%%%%%%%%%%%%%%%%%%%
%%%%%%%%%%%%%%%%%%%%%%%%%%%%%%%%%%%%%%%%%%%%%%%%%%%%%%%%%%%%%%%%%%%%%
%%%%%%%%%%%%%%%%%%%%%%%%%%%%%%%%%%%%%%%%%%%%%%%%%%%%%%%%%%%%%%%%%%%%%%
%%%%%%%%%%%%%%%%%%%%%%%%%%%%%%%%%%%%%%%%%%%%%%%%%%%%%%%%p3bis6.tex}

\part{Torsion analytique et sections ``concentrées"}

\section*{Introduction}\addcontentsline{toc}{section}{Introduction}
Dans la suite de cette partie, sauf mention
contraire, $E_n$ désignera un fibré holomorphe
hermitien  sur $X$ faisant partie de la famille
$\{E_{n}\, / n \in \Z \}$ avec $E_{0}$ qui sera un
fibré de degré quelconque, de rang deux, et
$E_{n}$ sera un fibré  défini en tensorisant le
faisceau $\mathcal E$ associé à $E_{0}$, par le
faisceau $\ox(n).$ On effectue les calculs en
coordonnées explicites pour obtenir des calculs
précis. On choisit donc un système de coordonnées
sur X noté $z$ tel que $z_{0}$ soit zéro. On pose
$$R_{h}(E_{n}):=R_{h,\C}(E_{n})dz\wedge d\bar
z:=R_{h,\R}(E_{n})dx\wedge dy,$$ avec
$$R_{h,\R}(E_{n})=-2iR_{h,\C}(E_{n}).$$
On choisit aussi un repère local holomorphe de
$E_{0}$ autour de $z_{0}$ tel que dans cette
situation  on a pour $z_{0}\in X$
:
\begin{equation}\label{courbb}
 R_{h,\C}(E_{0})(z_{0})=\begin{pmatrix} d_{1}&0\\ 0&d_{2}\\
\end{pmatrix},
\end{equation}
 alors
:$$R_{h,\C}(E_{n})(z_{0})=\begin{pmatrix} d_{1}+n&0\\
0&d_{2}+n\end{pmatrix}\\
=\begin{pmatrix} n_{1}&0\\
0&n_{2}\\
\end{pmatrix}.$$

\label{introd}
 Cette troisième partie est la partie
technique de notre travail . Les notations,
définitions des parties précédentes seront
conservées. Cette partie se divise en deux
sous-parties :
\begin{itemize}
  \item La première consiste en l'obtention précise de la variation  de la fonction zêta lorsque
  l'on fait varier les métriques hermitiennes sur le fibré $E_{n}$
  ce qui nous donne le résultat suivant :
  \begin{prop*}(cf. p. \pageref {afin})
  Soient deux m\'etriques hermitiennes $h_1$ et $h_{0}$ sur $E_{n}$  . On consid\`ere le chemin $h(u)=uh_{1}+(1-u)h_{0}$ avec $u
\in
\lbrack0, 1\rbrack$. Alors \\
  $${\du}\tau_{h(u)}(X, E_{n})=\frac{1}{2i\pi}
\int_{X}tr(\alpha R_{h(u)}(E_{n}))+
\int_{X}\tr(\alpha P_{n}).$$
  \end{prop*}
  ( On pose
$\alpha :=-{*^{-1}_{E_n}}{\du}*_{E_n}$).
  On peut remarquer que, puisque la dimension de $H^{1}(X,E_{n})$ est nulle pour $n$ assez grand,
   nous avons le théorème de Riemann-Roch qui s'écrit pour un fibré $(E_{n},h)$ avec $n$ assez grand :

  \begin{eqnarray*}
\dim H^{0}(X,E_{n})&=&\frac{-1}{2i\pi}
\int_{X}tr(R_{h}(E_{n}))+2(1-g).\\
\end{eqnarray*}

C'est pourquoi, de façon heuristique, nous pouvons
penser, puisque $tr (P_{n})$ correspond à la
dimension de $H^{0}(X,E_{n})$, que
$$\lim_{n\rightarrow \infty}{\du}\tau_{h(u)}(X,
E_{n})$$ s'annule.
      Nous montrons, en effet, cela (cf. ci
dessous).

  \item La deuxième vise à construire des sections holomorphes ``concentrées" du fibré $E_n$
  à partir de sections  particulières. Cette terminologie,
  employée, dans un autre cadre, dans \cite{don3} ,
  nous paraît adaptée à notre situation.Notons également l'utilisation de ce type d'outil dans \cite{tian}.
 Cette construction géométrique nous permet d'obtenir deux résultats intéressants dont l'un est, comme prévu :
  \begin{thm*}(cf. p. \pageref{fini})

$${\du}\tau_{h(u)}(X, E_{n})=O(\negli{})
.$$
\end{thm*}
\begin{remq}\label{rempans}
Après lecture d'une version préliminaire, P. Pansu
a remarqué que cette formule coïncidait avec la
variation de la formule de S.K. Donaldson (cf.
\cite{don2}ou p. \pageref{norma}), et que l'on
pouvait assez facilement identifier le terme
$\int_{X}tr(\alpha R_{h}(E_{n}))$ (resp.
$\int_{X}tr(\alpha P_{n})\omega$ avec $P_{n}$
projection orthogonale pour $h$ sur
$W_{n}$)\label{notapro} comme correspondant à la
variation de $\ln Norm_{h,
\omega}(\beta)$ (resp. de $\ln\det_{W_n}(L_{n}(h))$)(cf p. \pageref{remarquep}). Cette
remarque permet, en s'appuyant sur la formule de
S.K Donaldson \cite{don2}, de se limiter à la
lecture du chapitre huit (qui cependant utilise
des résultats du chapitre sept) dans cette
troisième partie pour la démonstration du théorème
\ref{princ}. Cependant, grâce à cette remarque,
nous obtenons, en fait, une nouvelle démonstration
de la formule de S. K. Donaldson \cite{don2} que
nous exposons dans les chapitres six et sept.
L'utilité de ceux-ci réside dans le fait que
l'énoncé de S .K . Donaldson \cite{don2} ne
s'applique qu'aux fibrés de caractéristique d'
Euler nulle; notre démonstration est valable en
degré quelconque.
\end{remq}\\
L'approche géométrique de la deuxième sous partie,
nous permet également de mieux analyser le
``passage" de $Met(W_n)$ à $Met(E_n)$. Ainsi, le
dernier paragraphe bénéficie de la construction
précédente pour obtenir un résultat sur
$I_{n}\circ L_{n}(h)$ en fonction de la courbure
de $(E_{n},h)$ qui nous permet d'obtenir le second
théorème
:
\begin{thm*}(cf. p. \pageref{end}) Soit $(E_{n},h)$ un fibré hermitien sur $X$, on a :
$$h^{-1}I_{n}L_{n}(h)=
\frac{\pi}{n^{2}}(R_{\C}(\ox(n)Id_{E_{0}}
-R_{\C,h}(E_{0})Id_{\ox(n)}+rId_{E_{n}})+{\bf{O}(\negli{3})}$$
(ici $r$ ne dépend que de $X$).De plus cela est
vrai uniformément en $h\in Met(E,C^{4})$.

\end{thm*}

Ces deux résultats achèvent les démonstrations des
théorèmes \ref{th1} et \ref{princ}.

\end{itemize}
\chapter{Etude de la variation de la torsion analytique de $(E_{n},h)$}
Les notations sont celles du paragraphe
\ref{cien}. \\ Nous donnons, dans ce premier
paragraphe, un résultat préalable utile au calcul
de ${\du}\tau_{h(u)}(X, E_n)$
.

\section{Variation du Laplacien $\Delta_{h, q}$ le long d'un chemin $h(u)$ }
Soient deux m\'etriques hermitiennes  $h_{0}$ et
$h_1$ sur $E_{n}$. On consid\`ere le chemin $h
(u)=uh_{1}+(1-u)h_{0}$ avec $u
\in
\lbrack0, 1\rbrack$.\\
Nous reprenons les notations du paragraphe
\ref{cien}.Nous posons $\alpha
:=-{*^{-1}_{E_n}}{\du}*_{E_n}$. Nous avons alors la
proposition suivante :
\begin{prop}\label{chap8}
\begin{eqnarray*}
{\du}\lbrack-\tr(\explapid{h(u)}{1})\rbrack& =&t.
\tr((\db\alpha^*\dbaretid{u} +\ddetid{u}\alpha
) (-\explapid{h(u)}{1}))
\end{eqnarray*}
\end{prop}
\begin{demo}
Etudions tout d'abord les variations du Laplacien
$$\Delta_{h(u), q}=\ddetid{u}+\detdid{u}$$ sur $C^{q}$
(cf. paragraphe \ref{cien})  en fonction des
variations de  celle de $(E_{n},h(u))$.
 Alors (cf. paragraphe
\ref{lapl}): $$\dbaret_{u}
=-*_{{E_n}^*}\db*_{E_n}\,  .$$ Mais $*_{{E_n}^*}$ et $*_{E_n}$
d\'ependent des m\'etriques choisies. Cependant
comme $$*_{{E_n}^*}*_{E_n}
=cste$$ on obtient :
$${\du}*_{{E_n}^*}*_{E_n}+*_{{E_n}^*}{\du}*_{E_n}=0 .$$ On pose
$$\alpha :=-{*^{-1}_{E_n}}{\du}*_{E_n}\hbox{ et
}\alpha^*:=-{\du}*_{{E_n}^*}{*^{-1}_{{E_n}^*}}.\label{araj}$$
On obtient alors, puisque  $\db$ est invariant par
changement de m\'etrique  :
\begin{eqnarray}
\du \dbaretid{u} &=&\du (-*_{{E_n}^*}\db * _{E_n})\nonumber\\
&=&-\du(*_{{E_n}^*})\db *_{E_n}
-*_{{E_n}^*}\db\du ( *_{E_n})\nonumber\\
&=&\alpha^{*} \dbaretid{u} +\dbaretid{u}
\alpha\nonumber
\end{eqnarray}
\begin{eqnarray}\label{formul}
{\du}\Delta_{h(u), q}&=&
\db{\du}\dbaretid{u}+{\du}\detdid{u}\nonumber\\ &=&
\db\alpha^{*}\dbaretid{u}+\ddetid{u}\alpha +\alpha ^* \detdid{u}+\dbaretid{u}\alpha
\db\hbox{ sur }C^{q}.
\end{eqnarray}
On utilise alors la proposition suivante
:
\begin{prop}(cf\cite{ray1})\label{formid}
$${\du}\tr(exp-t \Delta_{h(u), q})=-t.\tr(\dot{\Delta}_{h(u), q}(exp-t
\Delta_{h(u), q})$$
avec $\dot{\Delta}_{h(u), q}={\du}\Delta_{h(u),
q}$\\
 (ici $q=0$ ou $q=1$).
\end{prop}
Si $q=1$ on obtient d'après l'équation
\ref{formul}$${\du}\Delta_{h(u),
1}=\db\alpha^{*}\dbaretid{u}+\ddetid{u}\alpha.$$
Ceci et la proposition \ref{formid} achèvent la
démonstration de la proposition \ref{chap8}.
\end{demo}
\begin{remq}\label{ttt}
\end{remq}
Si $q=0$ on obtient $${\du}\Delta_{h(u), q}=\alpha
^* \detdid{u}+\dbaretid{u}\alpha
\db.$$
\begin{remq}
\label{tttt}
 On pose $$h(u)=h_{0}H_{u},$$  avec
$(\phi,\psi)_{h(u)}=(H_{u}\phi,\psi)_{h_{0}}=(\phi,H_{u}\psi)_{h_{0}},$
où $\phi,\psi\in\cin(X,\Lambda^{0,0}T^{*}X\otimes
{E_n}) $ ou
$\phi,\psi\in\cin(X,\Lambda^{1}T^{*}X\otimes
{E_n}) $.
 Soient
$\eta\in
\cin(X,\Lambda^{0,0}T^{*}X\otimes {E_n})$ et $\zeta\in
\cin(X,\Lambda^{0,1}T^{*}X\otimes {E_n}).$
On calcule alors :
  \begin{eqnarray*}
\int_{X}(\dbaretid{u}\zeta,\eta)_{h(u)}&=&1.\int_{X}(\zeta,\db\eta)_{h(u)}=
\int_{x}(H_{u}\zeta,\db\eta)_{h_{0}}\\
&=&\int_{X}(\dbaretid{0}(H_{u}\zeta),\eta)_{h_{0}}\\
&=&\int_{X}(H_{u}^{-1}\dbaretid{0}(H_{u}\zeta),\eta)_{h(u)}.
  \end{eqnarray*}
D'où :
$$\dbaretid{u}=H_{u}^{-1}\dbaretid{0}H_{u}.$$
Comme $$H_{u}H_{u}^{-1}=I,$$ on a
:$$\du H_{u}.H_{u}^{-1}+H_{u}.\du H_{u}^{-1}=0,$$
donc $$\du H_{u}^{-1}=-H_{u}^{-1}\du
H_{u}.H_{u}^{-1}.$$ On en déduit
:
$$\du\dbaretid{u}=-H_{u}^{-1}\du H_{u}.H_{u}^{-1}\dbaretid{0}H_{u}
+H_{u}^{-1}\dbaretid{0}\du H_{u},$$ d'où
$$\du\dbaretid{u}=-H_{u}^{-1}\du H_{u}.H_{u}^{-1}\dbaretid{0}H_{u}
+H_{u}^{-1}\dbaretid{0}H_{u}H_{u}^{-1}\du H_{u}.$$
Posons $\beta=-H_{u}^{-1}\du H_{u}$ alors :
$$\du\dbaretid{u}=-\beta.H_{u}^{-1}\dbaretid{0}H_{u}+H_{u}^{-1}\dbaretid{0}H_{u}\beta=
\dbaretid{u}\beta-\beta.\dbaretid{u}.$$

On remarque alors que : $$\beta=\alpha\mbox { et
surtout }\alpha^{*}=-\beta.$$ On a également
\begin{equation}\label{alphavu}
h(u)^{-1}\du h(u)=H_{u}^{-1}\du H_{u}=-\alpha.
\end{equation}

\end{remq}

\section{Opérateurs de Hilbert-Schmidt et opérateurs bornés :
 Reformulation de ${\du}\tau_{h(u)}(X,
E_n)$  }
 \label{etutor}
 On se propose, dans ce paragraphe, d'étudier la
variation de la torsion analytique par rapport aux
métriques mises sur $E_{n}$, avec $n$ assez grand.
Soient  deux m\'etriques
  hermitiennes $h_1$ et $h_{0}$ sur $E_{n}$. On
 consid\`ere le chemin $h (u)=uh_{1}+(1-u)h_{0}$
                     avec $u\in\lbrack0, 1\rbrack$.
                      On reprend les définitions du paragraphe \ref{definit}.
                  On a alors :
 \begin{eqnarray*}
  {\du}\tau_{h(u)}(X,E_n)&=&\du\lbrack\sum_{q=0}^{q=1}(-1)^q
q
 {\ds}\lbrack{\frac{1}{\Gamma (s)}}\int_0^\infty
        t^{s-1} \tr(\explap{h(u)}{q})\,
         dt\rbrack\rbrack_{\vert s=0}\\
                &=&-{\du}\lbrack
{\ds}{\frac{1}{\Gamma (s)}}\int_0^\infty t^{s-1}
 \tr(\explap{h(u)}{1})dt \rbrack _{\vert s=0} \\
        &=&-{\ds}\lbrack{\frac{1}{\Gamma
         (s)}}\du\int_0^\infty t^{s-1} (
\lbrack \tr\explap{h(u)}{1}\rbrack )dt\rbrack_{\vert s=0}\\
&=&-{\ds}\lbrack{\frac{1}{\Gamma (s)}}\int_0^\infty t^{s-1} (\du
\lbrack \tr\explap{h(u)}{1}\rbrack )dt\rbrack_{\vert s=0}.\\
                 \end{eqnarray*} Nous  supposons $n$ assez grand. Le noyau
de $\tr\explapid{h(u)}{1}$ est alors réduit à
zéro, donc $P_{n}^{\perp}$(cf. définition
\ref{toranalyt}) correspond à l'identité, nous
omettrons alors l'écriture de $P_{n}^{\perp}$ dans
le terme $\tr\explap{h(u)}{1}$.
\begin{prop}\label{rrrrr}On a
\begin{eqnarray*}
{\du}\tau_{h(u)}(X, E_n)&=&
 {\ds}\lbrack{\frac{1}{\Gamma
(s)}}\int_0^\infty t^{s}\lbrack\tr\alpha{\dt}
(\explapid{h(u)}{1})\\&+&\tr\alpha^{*}{\dt}\explap{h(u)}{0})\rbrack
dt\rbrack_{\vert s=0}.
\end{eqnarray*}
\end{prop}
\begin{demo}
Nous avons par définition
$${\du}\tau_{h(u)}(X,
E_n)=-{\ds}\lbrack{\frac{1}{\Gamma
(s)}}\int_0^\infty t^{s-1} (\du
\lbrack \tr\explapid{h(u)}{1}\rbrack )dt\rbrack_{\vert s=0}.$$
Nous utilisons la proposition \ref{chap8}. On doit
donc étudier la quantité
:
 \begin{equation}\label{eqa1}
{\du}\lbrack-\tr(\explapid{h(u)}{1})\rbrack=t.
\tr((\db\alpha^*\dbaretid{u} +\ddetid{u}\alpha^*) (-\explapid{h(u)}{1})) .
\end{equation}
 Nous allons pour cela utiliser les propriétés de la trace  pour commuter  des opérateurs particuliers.
 Remarquons   maintenant le fait que l 'op\'erateur de Laplace permute
avec l'op\'erateur $\db$ et son adjoint
$\dbaretid{u}$:

\begin{eqnarray*}
\dbaretid{u}\Delta_{h(u), q}&=&\Delta_{h(u), q-1}\dbaretid{u},\\
\db \Delta_{h(u), q}&=&\Delta_{h(u), q+1}\db,
\end{eqnarray*}
(ici $q=0$ ou $q=1$).  \\ De plus les op\'erateurs
consid\'er\'es sont Hilbert-Schmidt c'est à dire
:
\begin{dfn}Soit J un opérateur sur un espace de Hilbert T muni d'une base hilbertienne orthonormale
$(e_{i})_{i}$. J est un opérateur de
Hilbert-Schmidt si et seulement si
$$\Vert J \Vert^{2}_{H.S}=\sum_{i}\Vert Je_{i}\Vert^{2}_{T}=
\sum_{i, j}\vert(Je_{i}, e_{j})_{T}\vert^{2}<\infty .$$
\index{Op\'erateur!de Hilbert-Schmidt}
\end{dfn}On a également la
propri\'et\'e suivante :
\begin{prop}(\cite{berline})
Etant donn\'es un op\'erateur de Hilbert-Schmidt J
sur T et un opérateur borné K sur T on a :
$$\tr(JK)=\tr(KJ).$$
\end{prop}
On applique cette propriété à l'opérateur
$exp-t\frac{1}{2}\Delta_{h(u), 1}$ qui est de
Hilbert Schmidt et à $\Delta_{h(u),
1}\db\alpha^*\dbaret
\exp-t\frac{1}{2}
\Delta_{h(u), 1}$ qui est un opérateur borné. On
\'etudie donc les deux quantit\'es suivantes qui proviennent
de \ref{eqa1}.
\begin{eqnarray}
\label{eqa2} \tr(\db\alpha^*\dbaretid{u}\explapid{h(u)}{1}),& &\\
\label{eqa3} \tr(\ddetid{u}\alpha\explapid{h(u)}{1}).
\end{eqnarray}
 On effectue le calcul pour \ref{eqa2},  l'
autre s'obtenant de la même façon.

\begin{eqnarray*}
(\ref{eqa2})&=&\tr(\exp(-t\frac{1}{2}\Delta_{h(u),
1})(\db\alpha^*\dbaretid{u} \exp(-t\frac{1}{2}
\Delta_{h(u), 1})),\\ &=&\tr(\alpha
^* \dbaretid{u} \exp-t\Delta_{h(u), 1}\db),\\
 &=&\tr(\alpha^*\detdid{u} \explap{h(u)}{0}).
\end{eqnarray*}

De m\^eme on obtient :

\begin{eqnarray*}
(\ref{eqa3})&=&\tr(\alpha\ddetid{u}\explapid{h(u)}{1}).\\
\end{eqnarray*}

D'o\`u finalement :

\begin{eqnarray*}
(\ref{eqa1})&=&t.\tr(\alpha^*\detdid{u}\explap{h(u)}{0})+t.
\tr(\alpha\ddetid{u}\explapid{h(u)}{1}).\\
 \end{eqnarray*}

On obtient donc (cf. remarque p. \pageref{ttt}):

\begin{eqnarray*}
(\ref{eqa1})&=&{\du}\lbrack-\tr(\explapid{h(u)}{1})\rbrack\\
&=&t.\tr(\alpha^{*}(\detdid{u}+\ddetid{u})\explap{h(u)}{0})+t.\tr(\alpha(\detdid{u}+
\ddetid{u})\explapid{h(u)}{1})\\
&=&t.(\tr\alpha^{*}{\dt}
(\explap{h(u)}{0})+\tr\alpha{\dt}
(\explapid{h(u)}{1})).
\end{eqnarray*}
Ceci achève la démonstration de la proposition
\ref{rrrrr}.
\end{demo}

\section{Etude de de ${\du}\tau_{h(u)}(X,
E_n)$ } Nous reprenons avec les notations données
p. \pageref{araj} :
\begin{eqnarray*}
  {\du}\tau_{h(u)}(X,
E_n)&=&{\ds}\lbrack{\frac{1}{\Gamma
(s)}}\int_0^\infty t^{s-1}.t.\lbrack
\tr\alpha^{*}{\dt} (\explap{h(u)}{0})\\&+&\tr\alpha{\dt}
(\explapid{h(u)}{1})\rbrack dt\rbrack_{\vert s=0}.
\end{eqnarray*} Nous posons, pour alléger les notations
 :
$$A={\du}\tau_{h(u)}(X,
E_n).$$ Soient $\epsilon>0$ et $M>0$, nous allons
étudier dans ce chapitre la quantité
$A_{\epsilon,M}$ décrite ci-dessous. Comme $A$
correspond à $A_{0,\infty}$ nous en déduisons une
écriture simplifiée de $A$ qui sera exploitée au
chapitre suivant. En effet, nous obtenons la
proposition suivante :
\begin{prop}\label{resul}
$$A={\du}\tau_h(u)(X,
E_n)=\frac{1}{4\pi}(\tr\alpha^{*}\Phi_1(x,
x;\Delta
_{h(u),0})+\tr\alpha\Phi_1(x,x;\Delta
_{h(u), 1}))-\tr(\alpha^{*}P_{n}),$$
où $\Phi_1$ est défini p. \pageref{chaleur}.
\end{prop}
 On commence d'abord par intégrer
par parties la quantité
:
$$A_{\epsilon, M}:={\ds}\lbrack {\frac{1}{\Gamma
(s)}}\int_\epsilon^M
t^{s}\tr\alpha^{*}({\dt}(\explap{h(u)}{0})
+{\dt}(\tr\alpha(\explapid{h(u)}{1})dt\rbrack_{\vert
s=0}.$$
\begin{eqnarray*}
A_{\epsilon, M}&=&{\ds}\lbrack {\frac{1}{\Gamma
(s)}}\int_\epsilon^M
t^{s}\tr\alpha^{*}({\dt}(\explap{h(u)}{0})
+{\dt}(\tr\alpha(\explapid{h(u)}{1})dt\rbrack_{\vert
s=0} \\
 &=&{\ds}\lbrack {\frac{1}{\Gamma (s)}}
 \lbrack t^s \tr\alpha^{*}(\explap{h(u)}{0})+\tr\alpha(\explapid{h(u)}{1}) \rbrack_\epsilon ^M
 \rbrack_{\vert s=0} \\
  & -&{\ds}\lbrack{\frac{1}{\Gamma (s)}}s\int_\epsilon^M t^{s-1}\lbrack
  \tr\alpha^{*} (
\explap{h(u)}{0})+ \tr\alpha(\explapid{h(u)}{1})\rbrack dt\rbrack_{\vert
s=0} \\
   &=& A_{1, \epsilon, M}+A_{2, \epsilon, M}.
\end{eqnarray*}
Dans cette étude,  en ce qui concerne les
notations,  le deuxième et le troisième indice
indiquent les bornes d'intégration ordonnées. Nous
montrons dans le paragraphe suivant que :
$$\lim_{s\lr 0
}A_{1, 0, \infty}= 0,$$ le terme $A_{2, \epsilon,
M}$ sera étudié en utilisant  un développement
asymptotique du noyau de la chaleur.
 \section{Etude de $A_{1, \epsilon, M}$}
%%%%%%%%%%%%%%%%%%%%%%%%%%%%%%%%%%%%%%%%%%%%%%%%%%%%%
\begin{prop}$$\lim_{s\lr 0
}A_{1, 0, \infty}= 0$$
\end{prop}

\begin{demo}On rappelle que $$\lim_{s\rightarrow 0}\frac{1}{\Gamma (s)}=0$$
et que $$A_{1, \epsilon, M}={\ds}\lbrack
{\frac{1}{\Gamma (s)}}
 \lbrack t^s \tr\alpha^{*}(\explap{h(u)}{0})+\tr\alpha(\explapid{h(u)}{1}) \rbrack_\epsilon ^M
 \rbrack_{\vert s=0}.$$
Nous avons alors :
\begin{itemize}
 \item $\lim _{t \rightarrow
0}\frac{1}{\Gamma (s)}t^s (\tr\alpha^{*}
(\explap{h(u)}{0})+t^s
\tr\alpha(\explapid{h(u)}{1}))=0$
\item On a également :
  \begin{equation*}
{\ds}(\lim_{t\rightarrow \infty}\lbrack
{\frac{1}{\Gamma (s)}}
 \lbrack  t^s \tr\alpha^{*}(-\explap{h(u)}{0})+ t^s  \tr\alpha(\explapid{h(u)}{1})
 \rbrack)_{\vert s=0}=0.
  \end{equation*}

En effet
 %%%%%%%%%%%%%%%%%%%%%%%%%%%%%%%%%%%%%%%%%%%%%%%%%%%%%%%
 la proposition suivante nous permet l'étude du terme \\
 ${\ds}\underset{t\to
\infty}{Lim}\lbrack {\frac{1}{\Gamma (s)}}
 \lbrack t^s \tr\alpha(\explapid{h(u)}{1})
 \rbrack_{\vert s=0}$.
  \begin{prop}(cf.\cite{berline})\label{majo}
Il existe $L>0$ tel que : $\forall t
\geq L$
$$\parallel\explapid{h(u)}{1}\parallel_{\ell}<c(\ell)e^{-t\frac{{\bf{\lambda }}
_1}{2}}$$ o\`u
${\bf{\lambda _1}}$ correspond \`a la premi\`ere
valeur propre non nulle de l'op\'erateur $\Delta
_{h(u), 1}$
.\\ Ici $\Vert
\, .\, \Vert_{l}$ désigne une norme $C^{l}$ sur les sections $C^{l}$ du fibré $E_{n}$ et $c(l)$
une constante dépendant de l.

\end{prop}
  Donc le
terme $\underset{t\to
\infty}{Lim}\lbrack {\frac{1}{\Gamma (s)}}
 \lbrack t^{s}.
 \tr\alpha(\explapid{h(u)}{1})
 \rbrack$ est nul.
Il reste alors
%%%%%%%%%%%%%%%%%%%%%%%%%%%%%%%%%%%%%%%%%%%%%%%%%%%%%
le terme : $${\ds}\underset{M\to
\infty}{Lim}\lbrack {\frac{1}{\Gamma (s)}}M^{s}
 \lbrack  \tr\alpha^{*}(\explapM{h(u)}{0})
 \rbrack_{\vert s=0}.$$On  utilise la  proposition suivante :
 \begin{prop}(cf.\cite{berline})\label{majoi}
 il existe $L>0$ tel que : $\forall t
\geq L$
$$\parallel\explapid{h(u)}{0}P_{n,\frac{1}{2}}^{\perp}\parallel_{\ell}
<c(\ell)e^{-t\frac{{\bf{\lambda }}
_1}{2}}
$$ o\`u
${\bf{\lambda}} _1$ correspond \`a la premi\`ere
valeur propre non nulle de l'op\'erateur $\Delta
_{h(u), 0}$.\\ Ici $\Vert
\, .\, \Vert_{l}$ désigne une norme $C^{l}$ sur les sections $C^{l}$ du fibré et $c(l)$
une constante dépendant de l.

\end{prop}
On obtient alors de la même façon  :
$${\ds}\underset{M\to
\infty}{Lim}\lbrack {\frac{1}{\Gamma (s)}}M^{s}
 \lbrack  \tr\alpha^{*}(\explapM{h(u)}{0})
 \rbrack\rbrack_{\vert s=0}=0.$$
  \end{itemize}
  \end{demo}
\section{Etude de $A_{2, \epsilon, M}$}
Par définition, on a :
$$ A_{2, \epsilon, M}={\ds}\lbrack{\frac{1}{\Gamma (s)}}s\int_\epsilon^M
t^{s-1}\lbrack
\tr\alpha^{*}(\explap{h(u)}{0})+\tr\alpha(\explapid{h(u)}{1})\rbrack
dt\rbrack_{\vert s=0}$$ Posons
\begin{eqnarray*}
  I_{2, \epsilon, M}&=&\frac{1}{\Gamma (s)}
s\int_\epsilon^M t^{s-1}\lbrack
\tr\alpha^{*}(\explap{h(u)}{0})+\tr\alpha(\explapid{h(u)}{1})\rbrack
 dt\\
 &=&\frac{s}{\Gamma (s)}
\int_\epsilon^M t^{s-1}
\tr\alpha^{*}(\explap{h(u)}{0})
 dt+\frac{s}{\Gamma (s)}
\int_\epsilon^M
t^{s-1}\tr\alpha(\explapid{h(u)}{1})
 dt\\
 &=&I_{2, \epsilon, M, 0}+I_{2, \epsilon, M, 1}.
\end{eqnarray*}
Etudions chacun de ces termes.
\subsection{Etude de $I_{2, \epsilon, M, 1}$. }\label {i1}
\begin{itemize}
\item Au voisinage de zéro,
nous utilisons un développement asymptotique du
noyau de la chaleur\index{Noyau de la chaleur}. On
a, en effet, le r\'esultat suivant
:
\begin{prop}(cf. \cite{berline})\label{chaleur}\index{Th\'eor\`eme!Asymptote du noyau de la chaleur}Asymptote
du noyau de la chaleur. Soit $t>0$,
 Pour $ x \in X$,  $\exists N _{0}\quad \forall N\geq N_{0}$
on a (cf. paragraphe \ref{nych}):
$$p_t(x, x)= {\frac{1}{4\pi t}}\sum_{i=0}^{i=N}t^i\varphi_i(x, x;\Delta
_{h(u), 1})
+o(t^{N-1}).$$
 \end{prop}
On pose :
$$\Phi_{i}(x, x;\Delta_{h(u), 1})=\int_{X}\varphi_{i}(x, x;\Delta
_{h(u), 1})dx \, .$$
Fixons $N>N_{0}$,  on a donc :
  \begin{eqnarray}\label{chal}
\tr\alpha(\explapid{h(u)}{1}) &=&{\frac{1}{4\pi
t}}\sum_{i=0}^{i=N}t^i \tr(\alpha\Phi_i(x,
x;\Delta
_{h(u), 1}))
+ o(t^{N-1}) \nonumber\\ &=&\frac{1}{4\pi t}
\tr(\alpha\Phi_0(x, x;\Delta
_{h(u), 1}))+\frac{1}{4\pi }
\tr(\alpha\Phi_1(x, x;\Delta
_{h(u), 1}))\nonumber\\&+&{\frac{1}{4\pi
t}}\sum_{i=2}^{i=N}t^i \tr(\alpha\Phi_i(x,
x;\Delta
_{h(u), 1}))
+ o(t^{N-1}).\nonumber
  \end{eqnarray}

 Prenons Re $(s)>-1$.\\
Effectuons formellement la différence des deux
premiers termes du développement asymptotique et
de $I_{2, 0, 1, 1}$:
\begin{eqnarray*}I_{2, 0, 1, 1}&-&\frac{s}{\Gamma(s)}{\frac{1}{4\pi(s-1)}}\tr (\alpha \Phi_0(x, x;\Delta
_{h(u), 1}))-\frac{s}{\Gamma(s)}{\frac{1}{4\pi(s)}}\tr (\alpha \Phi_1(x, x;\Delta
_{h(u), 1}))\\
\end{eqnarray*}

 Cela
d\'efinit une fonction analytique de $s$ pour Re
$(s)>-1$ et comme $\Gamma$ possède un pôle simple
en zéro, on obtient finalement une fonction
holomorphe s'annulant en $s=0$. Le développement
analytique de  $I_{2, 0, 1, 1}$ nous donne
$$\ds(I_{2, 0, 1, 1})\vert_{s=0}={\frac{1}{4\pi }}\tr\alpha\Phi_1(x, x;\Delta
_{h(u), 1}).$$

%%%%%%%%%%%%%%%%%%%%%%%%%%%%%%%%%%%%%%%%%%%%%%%%%%%%%%%%
\item Pour $M$ au voisinage de l'infini,
L'application $s\mapsto
\int_1^M t^{s-1}
\tr\alpha\explapid{h(u)}{1}dt$ est une fonction
régulière pour $s=0$
. De plus :
$$\underset{s\lr 0}{Lim}\frac{1}{\Gamma(s)}=0.$$ Donc
 ${\frac{s}{\Gamma (s)}}\int_1^M t^{s-1}
\tr\alpha\explapid{h(u)}{1}dt$ possède un zéro double.
Il s'ensuit,
$$\ds \lbrack{\frac{1}{\Gamma (s)}}s\int_1^M t^{s-1}
\tr\alpha \explapid{h(u)}{1}dt\rbrack_{\vert s=0}=0
.$$
On a donc la proposition suivante :
\begin{prop}\label{a1}
$$\ds\lim_{\varepsilon\lr 0
\atop M\lr \infty}(I_{2, 0, \infty, 1})\vert_{s=0}={\frac{1}{4\pi }}\tr\alpha\Phi_1(x, x;\Delta
_{h(u), 1}).$$
\end{prop}
\end{itemize}
\subsection{Etude de $I_{2, \epsilon, M, 0}$}
Le travail est similaire au travail précédent. On
remarque que puisque $P_{n}^{\perp}=Id-P_{n}$ :
$$I_{2, \epsilon, 1, 0}=\frac{s}{\Gamma (s)}\int_\epsilon^1(t^{s-1} \tr(\alpha^{*}
\explapid{h(u)}{0})dt-\frac{s}{\Gamma (s)}
\int_\epsilon^1(t^{s-1} \tr(\alpha^{*} \explapid{h(u)}{0}P_{n})dt.$$
Or $$\tr(\alpha^{*}
\explapid{h(u)}{0}P_{n})=\tr(\alpha^{*}P_{n}),$$
 donc on
obtient :
\begin{eqnarray*}
I_{2,\epsilon, M, 0}&=&I_{2, \epsilon, 1,
0}+I_{2,1, M, 0}\\&=&
\frac{s}{\Gamma (s)}\int_\epsilon^1(t^{s-1}
\tr(\alpha^{*}
\explapid{h(u)}{0})dt\\&+&\frac{s}{\Gamma (s)}
\int_1^M(t^{s-1} \tr(\alpha^{*} \explap{h(u)}{0})dt
-
\frac{s}{\Gamma (s)}\int_\epsilon^1(t^{s-1} \tr(\alpha^{*}P_{n})dt.
\end{eqnarray*}Pour $I_{2,1, M, 0}$ on
utilise  la proposition \ref{majoi}.
   La fonction $$s\longmapsto
\int_1^\infty t^{s-1} \tr(\alpha^{*} \explap{h(u)}{0}dt$$  est régulière
pour $s=0$ et $I_{2,1, M, 0}$
 possède un zéro double . La dérivée pour
$s=0$  de $I_{2,1, M, 0}$  est donc nulle.\\
 Pour $\frac{s}{\Gamma (s)}\int_\epsilon^1(t^{s-1}
\tr(\alpha^{*}
\explapid{h(u)}{0})dt$, on utilise
l'analogue de la proposition \ref{chaleur} :
\begin{prop}(cf. \cite{berline})\label{chaleur1}\index{Th\'eor\`eme!Asymptote du noyau de la chaleur}Asymptote
du noyau de la chaleur. Soit $t>0$,
 Pour $ x \in X$,  $\exists N _{0}\quad \forall N\geq N_{0}$
on a:
$$p_t(x, x)= {\frac{1}{4\pi t}}\sum_{i=0}^{i=N}t^i\varphi_i(x, x;\Delta
_{h(u), 0})
+o(t^{N-1})$$
 \end{prop}
On pose :
$$\Phi_{i}(x, x;\Delta_{h(u), 0})=\int_{X}\varphi_{i}(x, x;\Delta
_{h(u), 0})dx \, .$$
Fixons $N>N_{0}$,  on a donc :
$$\tr\alpha^{*}(\explapid{h(u)}{0})
={\frac{1}{4\pi t}}\sum_{i=0}^{i=N}t^i \tr(\alpha^{*}\Phi_i(x, x;\Delta
_{h(u), 0}))
+ o(t^{N-1}) \ .$$ Un travail similaire à celui
fait dans le paragraphe \ref{i1} donne :
$$\lim_{\epsilon\lr 0}\frac{s}{\Gamma (s)}\int_\epsilon^1(t^{s-1}
\tr(\alpha^{*}
\explapid{h(u)}{0})dt={\frac{1}{4\pi }}\tr\alpha^{*}\Phi_1(x, x;\Delta
_{h(u), 0}).$$

Enfin, nous avons pour $Re (s)>-1$ :
\begin{eqnarray*}
\ds \lbrack \frac{s}{\Gamma (s)}\int_0
^1t^{s-1} \tr(\alpha^{*}P_{n})dt\rbrack_{\vert s=0}
&=&
\ds(\lbrack\frac{s.t^{s}}{\Gamma (s).s}\rbrack ^{1}_{0}
\tr(\alpha^{*}P_{n}))_{\vert s=0}.\\
\end{eqnarray*}

 D'où la
proposition

\begin{prop}\label{a2}.
$$\ds(I_{2, 0, \infty, 0})\vert_{s=0}={\frac{1}{4\pi }}\tr\alpha^{*}\Phi_1(x, x;\Delta
_{h(u), 0})- \tr(\alpha^{*}P_{n}).$$
\end{prop}
Les propositions \ref{a1} et \ref{a2} nous donnent
la proposition \ref{resul}.
%%%%%%%%%%%%%%%%%%%%%%%%%%%%%%%%%%%%%%%%%%%%%%%%%%%%%%%%
%%%%%%%%%%%%%%%%%%%%%%%%%%%%%%%%%%%%%%%%%%%%%%%%%%%%%%%%%%%%
\chapter{Calcul de $\Phi_{1}$}

\section{Introduction}

 Soient  deux m\'etriques
hermitiennes $h_{1}$ et $h_{0}$ sur $E_{n}$
. On consid\`ere le chemin $h (u)=uh_{1}+(1-u)h_{0}$
avec $u
\in
\lbrack0, 1\rbrack$. On reprend les définitions du paragraphe \ref{definit}.
 On a posé précédemment :
$$\alpha :=-{*^{-1}_{E_{n}}}{\du}*_{E_{n}}\hbox{ et
}\alpha^*:=-{\du}*_{E_{n}^*} *^{-1}_{E_{n}^*}.$$
On a montré que
:
$$A=\frac{1}{4\pi}(\tr\alpha^{*}\Phi_1(x, x;\Delta
_{h(u),0})+\tr\alpha\Phi_1(x,x;\Delta
_{h(u), 1}))-\tr(\alpha^{*}P_{n}),$$
avec $P_{n}$ la projection sur $H^{0}(X,E_{n})$.

La fin de la remarque p. \pageref{tttt} nous
fournit la reformulation de A :
\begin{prop}

\begin{equation}\label{calculA}
  A={\du}\tau_{h(u)}(X,
E_{n})=\frac{1}{4\pi}(\tr\alpha\Phi_1(x,x;\Delta
_{h(u), 1})-\tr\alpha\Phi_1(x, x;\Delta
_{h(u),0}))+\tr(\alpha P_{n}).
\end{equation}

\end{prop}
Nous allons maintenant  étudier les deux premiers
termes en utilisant la formule de Lichnerowicz et
montrer que $${\du}\tau_{h(u)}(X,
E_{n})=\frac{1}{2i\pi}
\int_{X}tr(\alpha R_{h(u)}(E_{n}))+
 \tr(\alpha
P_{n}).$$
%%%%%%%%%%%%%%%%%%%%%%%%%%%%%%%%%%%%%%%%%%%%%%%%%%%%%%%%%%%%%%%%%%%%%%%%
%%%%%%%%%%%%%%%%%%%%%%%%%%%%%%%%%%%%%%%%%%%%%%%%%%%%%%%%%%%%%%%%%%%%%%%
%%%%%%%%%%%%%%%%%%%%%%%%%%%%%%%%%%%%%%%%%%%%%%%%%%%%%%%%%%%%%%%%%%%%%%%%
%%%%%%%%%%%%%%%%%%%%%%%%%%%%%%%%%%%%%%%%%%%%%%%%%%%%%%%%%%%%%%%%%%%%%
%%%%%%%%%%%%%%%%%%%%%%%%%%%%%%%%%%%%%%%%%%%%%%%%%%%%%%%%%%%%%%%%%%%%%%
%%%%%%%%%%%%%%%%%%%%%%%%%%%%%%%%%%%%%%%%%%%%%%%%%%%%%%%%%%%%%%%%%%%%%%%%%%%%%%
%%%%%%%%%%%%%%%%%%%%%%%%%%%%%%%%%%%%%%%%%%%%%%%%%%%%%%%%%%%%%%%%%%%%%%%%%%%%

\section{Formule de Lichnerowicz}
On pourra se référer à \cite{berline} et
\cite{gilkey}. Puisque X est une surface de
Riemann, on peut la considérer comme une variété
kälhérienne. Considérons, comme d'habitude, le
fibré $E_{n}$ sur cette variété. Comme nous avons
une structure complexe, nous  avons le complexifié
du fibré tangent qui s'écrit
:
$$TX\otimes_{\mathbb{R}}\mathbb{C}=T^{1,0}X\otimes T^{0,1}X.$$
Comme la variété est kälhérienne, ces deux fibrés
sont préservés par la connexion de Levi-Civita
dont la dérivée covariante est notée $\cov$ (cf.
p. \pageref{derico}). De plus le fibré
$\Lambda(T^{0,1}X)^{*}\otimes E_{n}$ est un module
de Clifford. Soit
$$s=s^{0,1}+s^{1,0}\in
\Gamma(X,(TX)^{*})=\Gamma(X,(T^{0,1}X)^{*})\oplus\Gamma(X,(T^{1,0}X)^{*}).$$ L'action de
Clifford est définie par :
$$c(s)\beta=\sqrt{2}(\varepsilon(s^{0,1})-\iota(s^{1,0}))
\beta\quad \forall s\in\cin(X,(TX)^{*})\mbox{ et }
\forall \beta\in\Lambda(T^{0,1}X)^{*}\otimes  ,$$ avec $\varepsilon$ le produit extérieur et
$\iota(u)$ la contraction avec l'élement dual de
$u$ défini par la forme quadratique provenant de
la métrique riemannienne de X. Cette action est
auto-adjointe car $s^{0,1}=\bar s^{1,0}.$ On
définit la courbure de la connexion de Levi-Civita
$\cov$ notée $R_{LC}$  qui est donnée par la
relation suivante
:
\\ Soient
$Y,Z,S,W$ des champs de vecteurs sur $X$,
$$R_{LC}(S,Y)Z=\cov_{S}\cov_{Y}Z-\cov_{Y}\cov_{S}Z-\cov_{\lbrack S,Y\rbrack}Z.$$
 On
considère maintenant  le Laplacien de Bochner
défini par
:
\begin{dfn}\label{boce}Soit $\partial_{z}$ un repère local   de
$T^{0,1}X$ pour la métrique hermitienne sur $X$.
Le Laplacien de Bochner de $(E_{n},h)$,
$\Delta_{h}^{0,\bullet}$, est défini par :
$$\Delta^{0,\bullet}=\cov^{E_{n}}_{\cov_{\partial_{z}}\partial_{\bar z}}-
\cov^{E_{n}}_{\partial_{z}}\cov^{E_{n}}_{\partial_{\bar z}}.$$ En posant
$$ \left \lbrace
\begin{array}{l}
\partial_{z}=\partial _{x}-i \partial _{y}\\
\partial_{\bar z}=\partial _{x}+i\partial _{y}
\end{array}
\right. ,$$
où $\{\partial _{x},\partial _{y}\}$ est un repère
local orthonormé de $TX$.
\end{dfn} Nous avons alors le résultat suivant
:
\begin{prop}(Formule de Lichnerowicz. cf. \cite{berline} p140)\label{lichne} Soit $(E_{n},h)$ un fibré vectoriel hermitien sur X.
Dans un système de coordonnées locales, on a :
$$\detdid{h}+\ddetid{h}=\Delta^{0,\bullet}+\varepsilon(d\bar z)\iota(dz)
((R_{LC}(\partial_{z},\partial_{\bar{z}})\partial_{z},\partial_{\bar
z})+ R_{h}(E_{n})(\partial_{z},\partial_{\bar
z})),$$ Ici on utilise le complexe suivant
:
\begin{equation}\label{eqan1}
C^0\stackrel {\db}{\longrightarrow} C^1
\longrightarrow 0
\end{equation}
avec $C^i=\Omega_{X}^{0, i}(E_{n}).$
\end{prop}

\section{Calcul de $\Phi_{1}$}
Les éléments du développement asymptotique  du
noyau de la chaleur se définissent de façon
récurrente. En général ces coefficients deviennent
rapidement difficiles à calculer  mais dans notre
situation ce calcul est relativement simple. Soit
d'abord $\partial_{z}$ un repère local orthonormé
de $T^{0,1}X$ pour la métrique hermitienne sur
$X$. On pose  :\label{repere}
$$ \left \lbrace
\begin{array}{l}
\partial_{z}=\partial _{x}-i \partial _{y}\\
\partial_{\bar z}=\partial _{x}+i\partial _{y}
\end{array}
\right. ,$$
où $\{\partial _{x},\partial _{y}\}$ est un repère
local  de $TX$ :
\begin{dfn}
On appelle laplacien généralisé du fibré
$(E_{n},h)$ un opérateur sous la forme :
$$\Delta=\Delta^{E_{n}}+F,$$avec
$$\Delta^{E_{n}}=-\sum_{i=1,2}\cov^{E_{n}}_{e_{i}}\cov^{E_{n}}_{e_{i}}-
\cov_{\cov_{e_{i}}e_{i}}$$
avec $e_1:=\partial_{x}$ et $e_2:=\partial_{y}$
repère local pour $TX$ et $F$ une section du fibré
$End(E_{n})$.
\end{dfn}
Nous avons alors le résultat suivant :
\begin{prop}(\cite{berline} p. 84)\label{phi11}
$$\Phi_{1}(x,x,\Delta^{E_{n}}+F)=\frac{1}{3}R_{X}(x)-F(x),$$
où $R_{X}$ est la courbure scalaire de $X$.
\end{prop}
Donc si nous reprenons la proposition \ref{lichne}
nous allons expliciter la formule entre le
laplacien $\Delta^{E_{n}}$ et le laplacien
$\Delta^{0,\bullet}$. Cela donne pour le laplacien
$\Delta^{0,\bullet}$
:
\begin{eqnarray*}
  \Delta^{0,\bullet}&=&\cov^{E_{n}}_{\cov_{\partial_{z}}\partial_{\bar z}}-
\cov^{E_{n}}_{\partial_{z}}\cov^{E_{n}}_{\partial_{\bar z}}\\
&=&-\cov^{E_{n}}_{\partial
_{x}-i\partial _{y}}\cov^{E_{n}}_{\partial _{x}+i\partial _{y}}+
\cov^{E_{n}}_{\cov_{\partial _{x}-i\partial _{y}}\partial _{x}+i\partial _{y}}\\
&=&-(\cov^{E_{n}}_{\partial
_{x}}-\cov^{E_{n}}_{i\partial _{y}})(\cov^{E_{n}}_{\partial _{x}}+\cov^{E_{n}}_{i\partial _{y}})+
\cov^{E_{n}}_{(\cov_{\partial _{x}}-\cov_{i\partial _{y}})(\partial _{x}+i\partial _{y})}\\
&=&-(\cov^{E_{n}}_{\partial
_{x}}-\cov^{E_{n}}_{i\partial _{y}})(\cov^{E_{n}}_{\partial _{x}}+\cov^{E_{n}}_{i\partial _{y}})+
\cov^{E_{n}}_{(\cov_{\partial _{x}}(\partial _{x}+i\partial _{y})
-\cov_{i\partial _{y}}(\partial _{x}+i\partial _{y}))}\\
&=&\Delta^{E_{n}}-i(\lbrack
\cov^{E_{n}}_{\partial _{x}},\cov^{E_{n}}_{\partial _{y}}\rbrack+
\cov^{E_{n}}_{(\cov_{\partial _{y}}\partial _{x}-
\cov_{\partial _{x}}\partial _{y})})\\
&=&\Delta^{E_{n}}-iR_{h}(E_{n})(\partial
_{x},\partial _{y})+i\cov^{E_{n}}_{\lbrack
\partial _{x},\partial _{y}\rbrack-\cov_{\partial _{x}}\partial _{y}-\cov_{\partial _{y}}\partial _{x}}\\
&=&\Delta^{E_{n}}-iR_{h}(E_{n})(\partial
_{x},\partial _{y})+i\cov^{E_{n}}_{R_{h}(\partial _{x},\partial _{y})}\\
&=&\Delta^{E_{n}}-iR_{h,\R}(E_{n}).
\end{eqnarray*}

 On obtient donc, pour un fibré $(E_{n},h)$ en
utilisant la proposition \ref{lichne}
$$\detdid{h}+\ddetid{h}=\Delta^{E_{n}}-iR_{h,\R}(E_{n})
+\varepsilon(d\bar z)i(dz)
((R_{LC}(\partial_{z},\partial_{\bar{z}})\partial_{z},\partial_{\bar
z})+ R_{h}(E_{n})(\partial_{z},\partial_{\bar
z})).$$ On a par ailleurs la proposition suivante
:
\begin{prop}\label{delta}Soit $(E_{n},h)$ un fibré vectoriel holomorphe sur $X$. Avec les notations précédentes on a :
$$ \left \lbrace
\begin{array}{l}
\tilde \Delta_{h,0}=\Delta^{E_{n}}-iR_{h,\R}(E_{n})\\
\tilde \Delta_{h,1}=\Delta^{E_{n}}+iR_{h,\R}(E_{n})
+
(R_{LC}(\partial_{z},\partial_{\bar{z}})\partial_{z},\partial_{\bar
z}).\\
\end{array}
\right.$$
\end{prop}
Notons  que $\tilde\laplac{h}{q}$  est l'opérateur
de Laplace ou Laplacien agissant sur
$C^{q}=\Omega_{X}^{0, q}(E_{n})$ défini par :
$$\tilde  \laplac{h}{q}=\ddetid{h}+\detdid{h}.$$
\begin{remq}
Etant donné que nous considérons, pour le calcul
de $\Phi_{1}$, l'opérateur de Laplace
$\laplac{h}{q}$ agissant sur $C^{q}=\Omega_{X}^{0,
q}(E_{n})$, nous remarquons que :
$$ \left \lbrace
\begin{array}{l}
 \Delta_{h,0}=\Delta^{E_{n}}-iR_{h,\R}(E_{n})
\\
 \Delta_{h,1}=\Delta^{E_{n}\otimes}+iR_{h,\R}(E_{n})+
(R_{LC}(\partial_{z},\partial_{\bar{z}})\partial_{z},\partial_{\bar
z}),\\
\end{array}
\right.$$
avec les notations définies p.\pageref{zzzzz}.\\
Dans le lemme 6.1.7 p. 212 de \cite{don4} il y a
une petite erreur de signe dans la première
formule.
\end{remq}\\
\begin{demo} Nous utilisons le résultat \ref{lichne}.
En effet pour les 0-formes, la contraction par le
vecteur $\partial _{x}+i\partial _{y}$, ce qui
correspond à $\iota(dz)$ est nulle. On obtient
donc la première équation de la proposition. Pour
ce qui est des $(0,1)$-formes, il est immédiat que
$$\forall \eta\in \cin(X,\Lambda^{0,1}T^{*}X\otimes E_{n}),\quad\epsilon(d\bar z)\iota(dz)\eta=\eta.$$

Remarquons que \label{courburee}:
\begin{eqnarray*}
R_{h}(E_{n})(\partial_{z},\partial_{\bar
z})&=&R_{h,\C}(E_{n}).
\end{eqnarray*}
De plus
$$R_{h}(E_{n})(\partial_{z},\partial_{\bar z})=R_{h,\C}(E_{n})=-2iR_{h,\R}(E_{n}).$$ Nous
obtenons donc le résultat
.
\end{demo}

Finalement on a :
\begin{prop}\label{afin}
$${\du}\tau_{h(u)}(X,
E_{n})=\frac{1}{2i\pi}
\int_{X}tr(\alpha R_{h(u)}(E_{n}))+
 \tr(\alpha
P_{n}).$$
\end{prop}
\begin{demo}
Nous avons (cf. prop.\ref{calculA}) :
$${\du}\tau_{h(u)}(X,
E_{n})=\frac{1}{4\pi}(\tr\alpha\Phi_1(x,x;\Delta
_{h(u), 1})-\tr\alpha\Phi_1(x, x;\Delta
_{h(u),0}))+\tr(\alpha P_{n}).$$
Etudions le terme $\tr\alpha\Phi_1(x, x;\Delta
_{h(u),1}))-\tr\alpha\Phi_1(x, x;\Delta
_{h(u),0}))$.
 D'après la
proposition \ref{phi11} et la remarque p.
\pageref{delta} on a
:
$$ \left \lbrace
\begin{array}{l}
\Phi_1(x, x;\Delta
_{h(u),0}))=\frac{1}{3}R_{X}+iR_{h(u),\R}(E_{n})\\
\Phi_1(x, x;\Delta
_{h(u),1}))=\frac{1}{3}R_{X}-iR_{h(u),\R}(E_{n})+
((R_{\cov}(\partial_{z},\partial_{\bar{z}})\partial_{z},\partial_{\bar
z})\\
\end{array}
\right.$$
Comme $\alpha$ est à trace nulle $$\tr\alpha
((R_{\cov}(\partial_{z},\partial_{\bar{z}})\partial_{z},\partial_{\bar
z})=0.$$ D'où
$$\tr\alpha\Phi_1(x, x;\Delta
_{h(u),1}))-\tr\alpha\Phi_1(x, x;\Delta
_{h(u),0}))=-2i\int_{X}tr\alpha R_{h(u)}(E_{n}).  $$
 On obtient donc le résultat escompté.
\end{demo}
La proposition est le résultat exposé
p.\pageref{normap} .
\begin{prop}(cf. p.
\pageref{normap})\label{degrequel}
Soit $E_{n }$ un fibré sur $X$.%L'équation (cf. p.
%\pageref{normap})
On a pour des métriques $h$ telles que
$det_{E_{n}}(h)=1$ :\\ $\exists n_{0}, \forall
n>n_{0}$
\begin{equation}
Norm_{h, \omega}(\beta )=
(\det_{W_{n}}(L_{n}(h)))^{{\frac{1}{2}}}\times
T_{h}(X, E_{n}).
\end{equation}
\end{prop}
\begin{demo}Soient des métriques $h$  telles que $det_{E_{n}}(h)=1.$
Nous avons alors :

$$\ln (Norm_{h,\omega}(\beta))=\ln(Norm_{k,\omega}(\beta))+\frac{1}{8\pi}(
\int_{X}R_{2}(h,k)) \mbox { cf. équation \ref{normalpha}
p.\pageref{normalpha}},$$ et $$\du
\int_{X}{{R_{2}}}(h(u),k)=2i\int_{X}tr(v_{u}R_{h(u)}(E_{n})),$$
où $k\in Met(E_{n})$ fixé. De plus, d'après
l'équation \ref{alphavu}, $v_{u}=-\alpha$
 donc $$\du\ln
(Norm_{h,\omega}(\beta))=\frac{2}{8\pi}
\int_{X}tr(v_{u}R(h(u)))=\frac{1}{4i\pi}\int_{X}tr(\alpha R_{h(u)}(E_{n})).$$

On a également : $${\du}\tau_{h(u)}(X, E_{n})=2\du
T_h(u)(X, E_{n}) \mbox { cf. définition p.
\pageref{toranalyt}}.$$ On obtient alors :
\begin{eqnarray*}
\du \ln (Norm_{h(u),\omega}(\beta))&=&\frac{1}{4i\pi}\int_{X}tr(\alpha R_{h(u)}(E_{n}))\\
&&\mbox{ et en utilisant la proposition
\ref{afin},}\\ &=&\frac{1}{2}{\du}\tau_{h(u)}(X,
E_{n})-\frac{1}{2}\tr(\alpha P_{n})\\
\mbox{ (cf. pour la notation
p.\pageref{notapro})}.&&
\end{eqnarray*}
De plus
  \begin{eqnarray*}
\du\ln(\det_{W_{n}}(L_{n}(h(u))^{{\frac{1}{2}}}))&=&\frac{1}{2}\du\ln(\det_{W_{n}}(L_{n}(h(u)))\\
&=&\frac{1}{2}
\frac{\du\det_{W_{n}}(L_{n}(h(u))}{\det_{W_{n}}(L_{n}(h(u))}\\
&=&\frac{1}{2}tr\
[\det(L_{n}(h(u))).(L_{n}(h(u)))^{-1}\du
L_{n}(h_{u})\ ].det(L_{n}(h(u)))^{-1}\\
&=&\frac{1}{2}\int_{X}tr ((h(u))^{-1}\du(h_{u})
P_{n})\omega\\ &=&\frac{1}{2}\int_{X}tr(\alpha
P_{n})\omega
  \end{eqnarray*}
%%%%%%%%%%%%%%%%%%cf rideau p52
%\du \ln
%(Norm_{h(u),\omega}(\beta)) &=&{\du}\ln
%T_{h(u)}(X, E_{n})+\du
%\ln(\det_{W_{n}}(L_{n}(h(u))^{{\frac{1}{2}}}).
\label{remarquep}
(ce qui fut la remarque de P. Pansu cf. p.
\pageref{rempans}).
Nous obtenons alors, à une constante près :
\begin{equation}
Norm_{h, \omega}(\beta )=
(\det_{W_{n}}(L_{n}(h)))^{{\frac{1}{2}}}\times
T_{h}(X, E_{n}).
\end{equation}

\end{demo}
Il reste à étudier maintenant le terme $\tr(\alpha
P_{n}).$ C'est l'objet du chapitre suivant où nous
construisons des sections `` concentrées ".

\chapter{Sections ``concentrées"}

Dans toute cette partie $E_{n}$ est un fibré
holomorphe sur $X$ de
 rang deux.
Nous fixons une métrique lisse de référence notée
 $k$, l'espace $Met(E,C^{4})$ est l'espace des métriques $h$, telles que $h=kH$
  avec $H$ endomorphisme de $E$ autoadjoint défini positif de  $ C^{4}(End(E)).$
 Dans ce chapitre, nous  étudions le terme $\tr(\alpha
P_{n})$.
\section{Choix d'un repère holomorphe local sur $E_{n}$}\label{baser}
Soit $h\in Met(E,C^{4})$. On considère un point
$z_{0}\in X$. Nous prenons un système de
coordonnées locales  $(U,z)$ avec $z_{0}\in U$ de
telle sorte  que $z_{0}$ soit donné par
  $0$. On va donc choisir un repère $\epsilon$ de $\ox(1)$, ce qui donnera un repère
   $\epsilon^{\otimes n}$ sur $\ox(n)$. On notera parfois (avec abus) $h_{\ox(n)}(z)$ la quantité
   $\nor{\epsilon^{\otimes n}}_{h\ox(n)}$.  Nous construisons
maintenant un repère holomorphe dans un voisinage
de $0$ $\{e_{1},e_{2}\}$  de $E_{n}$, de telle
sorte que la métrique $h$ sur
$E_{n}=E_{0}\otimes\ox(n)$ s'écrive :
$$h(z)= h_{\ox(n)}({\bf {1}}-\begin{pmatrix}
d_{1}&0\\ 0&d_{2}\\
\end{pmatrix}z\bar z +{\bf {A}}z\bar z^{2}+{\bf {B}}z^{2}\bar z^{}+ \varphi(z))
 ,$$
 avec $${\bf {1}}=\begin{pmatrix}
1&0\\ 0&1\\
\end{pmatrix}\mbox{ et }{\bf {A}},{\bf {B}}\in
M_{2}(\C).$$ La métrique $h_{\ox(n)}$ est celle du
fibré $\ox(n)$. De plus  $\varphi
(z)=O(\mo{z}^{3})$ et les coefficients des
matrices ${\bf {A}}$ et ${\bf {B}}$ sont
indépendants de $n$.
  On pose également :
$$h^{ij}(z)=(e_{i}(z),e_{j}(z))_{h(z)},$$
et
\begin{equation}\label{mthy}
h^{ij}_{kl}(z)=h^{ij}_{kl}z^{k}\bar z^{l}.
\end{equation}
Commençons par construire un repère local
holomorphe   sur $\ox(1)$
 et un repère sur $E_{0}
 $.
 %%%%%%%%%%%%%%%%

 \subsection{Construction d'un repère sur $\ox(1)$}

 Nous construisons un repère local
holomorphe $\epsilon$ sur $\ox(1)$  et choisissons
une coordonnée  $z$ tels  que
\begin{equation}\label{hoi}
  h_{\ox(1)}(z)=1-z\bar z + m.z^{2}\bar z^{2} +O(\vert z \vert^{5}),
\end{equation}
avec $m\in \R$.
 Nous imposons donc à la coordonnée $z$ de vérifier
$R_{h_{\ox(1)}}(z_{0})=dz\wedge d\bar z$ et notons
également que $\omega =R_{h_{\ox(1)}}$(cf. p.
\pageref{volume}).
 On aura alors en utilisant le développement
limité en zéro de $\exp $ :\\
 $h_{\ox(1)}(z)=\exp(-\mo{z}^{2}+  O(\mo{z}^{4})).$
\\
 Ce travail s'effectue en plusieurs
étapes. On peut supposer que nous avons un repère
local $\tilde\epsilon$ sur $\ox(1)$ et une
coordonnée $\zeta$  tels que
:
\begin{equation}\label{metr}
  \tilde h_{\ox(1)}(\zeta)=1+O(\mo{\zeta})
\end{equation}
(quitte à multiplier $h$ par un scalaire).
 Nous
effectuons maintenant un changement de repère
$\epsilon$, holomorphe, de la forme
$$\epsilon=(1+a_{1}\zeta+a_{2}\zeta^{2}+a_{3}\zeta^{3}+a_{4}\zeta^{4}+a_{5}\zeta_{5})\tilde\epsilon\, ,\quad a_{i}\in \C,i=1..5.
$$
On obtient alors :
$$ h_{\ox(1)}=\epsilon^{*} \tilde h_{\ox(1)}\epsilon.$$
On peut vérifier que  $ h_{\ox(1)}$ satisfait
\ref{metr} et, en choisissant $a_{i}$ pour
$i=1..5$, on impose à $ h_{\ox(1)}$ de vérifier
$$ h_{\ox(1)}(\zeta)
=1-\zeta\bar \zeta+...+O(\vert \zeta \vert^{5})$$
Cela donne des conditions sur les coefficients
$a_{i}$ pour $i=1..5$. Par exemple, le
développement de Taylor-Young de $dh_{\ox(1)}$ en
zéro donne :
$$dh_{\ox(1)}(\zeta)=
dh_{\ox(1)}(0)+{\zeta\bar \zeta
}d^{2}h_{\ox(1)}(0)++O(\vert \zeta\vert
^{3}).$$ En comparant l'écriture précédente à
cette dernière on obtient
:
$$ dh_{\ox(1)}(0)=0 \mbox{ et }   \frac{\partial ^{2}}{\partial \zeta \partial \bar \zeta}h_{\ox(1)}(0)=-1.$$
Mais $$dh_{\ox(1)}(0)=d(\epsilon^{*} \tilde
h_{\ox(1)}\epsilon)(0).$$ Donc
$$0=d\epsilon^{*}(0)+d \tilde
h_{\ox(1)}(0)+ d\epsilon^{*}(0)=\bar a+d \tilde
h_{\ox(1)}(0)+a,$$ car
$$\epsilon^{*}(0)=\epsilon(0)=\tilde
h_{\ox(1)}(0)=1.$$  Cela donne

$$
a_{1}=-\partial\tilde h_{\ox(1)}(0).$$ De même on
a
$$d^{2}h_{\ox(1)}(0)=d^{2}(\epsilon^{*} \tilde
h_{\ox(1)}\epsilon)(0).$$ D'où des conditions sur
$a_{2}$.Nous faisons de même pour
$a_{3},a_{4},a_{5}$.
 Nous obtenons
ainsi le repère local $\epsilon$.
La métrique sur
$\ox(1)$ s'écrit alors
:
$$h_{\ox(1)}(\zeta)=1-\zeta\bar \zeta +....+
O(\vert \zeta\vert^{5}),$$ avec $p\in \C$. Posons
maintenant
$$z=\zeta+c\zeta^{2}+e\zeta^{3}\mbox{ avec }c,e\in \C$$
Alors,  :
\begin{eqnarray*}
h_{\ox(1)}(z)&=&1-z\bar z+mz^{2}\bar z^{2 }O(\vert
z\vert^{5})\\
&=&1-(\zeta+c\zeta^{2}+e\zeta^{3})\overline{(\zeta+c\zeta^{2}+e\zeta^{3})}+O(\vert
\zeta\vert^{4})\\
&=&1-\zeta\bar \zeta +p\zeta^{2}\bar
\zeta+\bar p\zeta\bar \zeta^{2}+....+O(\vert
\zeta\vert^{5}),
\end{eqnarray*}
avec $c=p$ et $m=2 c\bar c$.De même on trouve $e$.
Ce choix de coordonnée $z$ étant établi, nous
avons la métrique $h_{\ox(1)}$ écrite sous la
forme désirée.
\subsection{Remarque sur la forme volume de $X$}
Notons que la forme volume  $\omega$ va vérifier
d'une part :\\
$\omega(z_{0})=R_{h_{\ox(1)}}(z_{0})=dz\wedge
d\bar z$.\\ D'autre part, la forme volume vérifie
localement autour de $z_{0}$ (cf. remarque
p.\pageref{und})
:

$$\omega=\bar\partial( h_{\ox(1)}^{-1}\partial  h_{\ox(1)}).$$
Cela donne :

\begin{equation}\label{vol}
  \omega(z)=(1+ r z\bar z +...)dz\wedge d\bar z.
\end{equation}

  Nous
posons$$\omega(z)=w(z)dz\wedge
\bar z.$$
Remarquons également que en $z_{0}$ :
$$rdz\wedge d\bar z=R(T^{}X)(z_{0}).$$

\subsection{Construction d'un repère sur $E_{0}$}
Nous construisons maintenant un repère
$\{f_{1},f_{2}\}$ sur $E_{0}$ tel que
$$h_{E_{0}}(z)={\bf 1}-\begin{pmatrix} d_{1}&0\\
0&d_{2}\\
\end{pmatrix}z\bar z +Mz\bar z^{2}+Nz^{2}\bar z^{}+ O(\mo{z}^{3}))$$
où ${\bf 1}=\begin{pmatrix} 1&0\\ 0&1\\
\end{pmatrix} \mbox{ et }M,N\in M_{2}(\C).$\\
Nous procédons comme ci-dessus : quitte à faire un
changement de coordonnées de $M_{2}(\C)$, on peut
supposer que dans un repère $\{\tilde f_{1},\tilde
f_{2}\}$ :
$$\tilde h_{E_{0}}(z)= 1+O(\mo{z}),$$  et
la matrice $R_{\C}(h_{E_{0}})(z_{0})$ est
diagonale. Nous faisons également un changement de
coordonnées holomorphe de la forme :
$$\begin{pmatrix}
f_{1}&\\ f_{2}&\\
\end{pmatrix}=({\bf 1}+Az+Bz^{2}+Cz^{3})\begin{pmatrix}
\tilde f_{1}&\\\tilde  f_{2}&\\
\end{pmatrix}\, ,\quad  A,B,C\in M_{2}(\C).$$
 On vérifie également que
$$h_{E_{0}}(z)=g^{*} \tilde h_{E_{0}}g={\bf
1}+O(\mo{z}).$$En procédant comme précédemment
c'est à dire en choisissant $A,B,C$,  on a
maintenant :

$$\tilde h_{E_{0}}(z)=h_{E_{0}}={\bf 1}-\begin{pmatrix} d_{1}&0\\
0&d_{2}\\
\end{pmatrix}z\bar z +Mz\bar z^{2}+Nz^{2}\bar z^{}+ O(\mo{z}^{3})).$$
Nous  avons ainsi défini le repère
$\{f_{1},f_{2}\}$.
\subsection{Repère sur $E_{n}$}\label{seize}
On vient donc de construire un repère sur
$E_{0}\otimes \ox(1)$ de la forme $\{
f_{1}\otimes\epsilon ,f_{2}\otimes\epsilon\}$. La
métrique h s'écrit alors dans ce repère
$$h(z)= h_{\ox(n)}\otimes({\bf
1}-\begin{pmatrix} d_{1}&0\\ 0&d_{2}\\
\end{pmatrix}z\bar z +Mz\bar z^{2}+Nz^{2}\bar z^{}+ O(\mo{z}^{3})).$$
Et $R_{E_{1}}(h)(z_{0})=\bar
\partial\partial h(z_{0})=\begin{pmatrix} d_{1}+1&0\\ 0&d_{2}+1\\
\end{pmatrix}dz\wedge d\bar z $. \\
Sur $E_{n}$ il suffit de
considérer le repère
:
$$\{f_{1}\otimes\epsilon^{\otimes n},f_{2}\otimes\epsilon^{\otimes n}\}$$
que nous noterons $\{{e_{1},e_{2}}\}$. Nous avons
alors
:\label{reperp}
$$h(z)= h_{\ox(n)}({\bf 1}-\begin{pmatrix}
d_{1}&0\\ 0&d_{2}\\
\end{pmatrix}z\bar z +{\bf {A}}z\bar z^{2}+{\bf {B}}z^{2}\bar z^{}+ \varphi(z)),$$
avec $\varphi(z)=O(\mo{z}^{3})$ et les matrices
$\bf {A}$ et $\bf {B}$ sont  des matrices
indépendantes de $n$ d'après la formule  \ref{hoi}
.
Remarquons également que nous avons la majoration
entre matrices hermitiennes suivante:\\ $\exists\,
n_{0}$,$\exists\,C_{1},C_{2}>0$ indépendants de
$n$, telles que $n\geq n_{0}$
\begin{equation}\label{merthes}
  \exp(-n C_{1}\vert
z\vert ^{2}){\bf 1}\leq h_{\ox(n)}(z)\leq \exp(-n
C_{2}\vert z\vert ^{2}){\bf 1}.
\end{equation}
De la même façon : Il existe $n'_0$ et
$C_3,C_{4}>0$  indépendants de $n$, tels que si
$n> n'_0$
:

\begin{equation}\label{merthes1}
\exp (-C_{4} n | z|^2) {\bf 1}
\leq C_2 (h(z))w(z).
\end{equation}
 De plus
\begin{eqnarray}
R_{E_{n}}(h)(z_{0})&=&\bar
\partial\partial h(z_{0})\nonumber\\
&=&Id_{\ox(n)}\otimes R_{E_{0}}+I_{E_{0}}\otimes
R_{\ox(n)}\\ &=&\begin{pmatrix} d_{1}+n&0\\
0&d_{2}+n\\
\end{pmatrix}dz\wedge d\bar z\label{ssss}.
\end{eqnarray}

\section{Rappels sur les distributions}
Nous donnons tout d'abord les définitions sur un
ouvert de $R^{2}$. Nous généraliserons cela à $X$
en considérant les cartes locales de $X$.\\ Soit
ouvert $\Omega $ de $\R^{2}$. Nous allons définir
sur cet ouvert ce que l'on appelle des
distributions. On considère le sous espace
vectoriel ${\mathcal D}_{K}(\Omega)$ de
$\cin(\Omega)$ composé des fonctions dont le
support est contenu dans $K$ compact de $\Omega$.
On  munit ${\mathcal D}_{K}(\Omega)$ de sa
structure d'espace localement compact définie par
la famille de semi-normes $(\nor{.}_{k})_{k\in
\N}$ telles que :
$$\nor{f}_{k}=\sup\{D^{\alpha}f(z),\vert\alpha\vert\leq k,z\in \Omega\}.$$
Ainsi ${\mathcal D}_{K}(\Omega)$ est un espace de
Fréchet.
 Si $\mathcal K$ désigne l'ensemble de
toutes les parties compactes de $\Omega$, on note
$${\mathcal D}(\Omega)=\bigcup_{K\in \mathcal K}{\mathcal D}_{K}(\Omega).$$
\begin{dfn}
Une forme linéaire $u\,:\,{\mathcal D}(\Omega)\lr
\C$ est appelée une distribution sur $X$ si on a :
  \begin{eqnarray}
\forall K \in {\mathcal K}\, ,&&\quad\exists
C_{k}\geq 0\, ,\quad
\exists k\in \N ,\nonumber\\
&&\vert u(f)\vert \leq
C_{K}\nor{f}_{k}\,,\quad\forall f\in {\mathcal
D}_{K}(\Omega)\nonumber.
  \end{eqnarray}

On note ${\mathcal D}'(\Omega)$ cet espace
vectoriel.
\end{dfn}
On peut remarquer que l'espace $L^{2}(\Omega)$
peut être plongé dans cet espace. La distribution
qui nous intéresse est la distribution de Dirac
$\delta_{z},z\in \Omega$ définie par :
$$\delta_{z}\, :\, f\in \cin(\Omega)\lr f(z)\in \C$$
\begin{dfn}On définit
une distribution $u$  sur $X$  si pour toute carte
locale de $X$, $U\subset X$, il existe une
distribution $v$ sur $U$ considéré comme ouvert
$\Omega$ de $R^{2}$ et pour toute fonction  $f\in
C^{\infty}(U)$(prolongée par zéro sur $X$ tout
entier, ce qui en fait un élément de
$C^{\infty}(X)$) telle que : $$v(f)=u(f).$$ On
note $\mathcal {D}'(X)$ cet ensemble.
\end{dfn}

Nous définissons alors les distribution sur $X$
pour un fibré vectoriel $E_{n}$ :
\begin{dfn}
On définit une section distributionnelle  $u$  sur
$X$ pour un fibré vectoriel $E_{n}$ de rang $2$
par l'existence sur chaque  ouvert trivialisant de
la forme  $U\times \C ^{2}$ de $E_{n}$ de
distributions v sur $C^{\infty}(U)^{2}$ identifié
avec $C^{\infty}(U,E_{n})$ telles que pour tout
élément $f\in C^{\infty}(U,E)$(prolongée par zéro
sur $X$ tout entier, ce qui en fait un élément de
$C^{\infty}(X,E)$) on a
 $$v(f)=u(f).$$ On note $\Gamma
^{-\infty}(X,E_{n})$ cet ensemble.
\end{dfn}

\section{Choix d'une section correspondant à la première composante du noyau de $P_{n}$}\label{cho}
La proposition \ref{afin} donne :
$${\du}\tau_{h(u)}(X,
E_{n})=\frac{1}{2i\pi}
\int_{X}tr(\alpha R(E_{n}))+
 \tr(\alpha
P_{n}).$$ On étudie alors $
\tr(\alpha P_{n})$,  avec $P_{n}$ projection
orthogonale
  sur $H^{0}(X,E_{n})$. C'est la situation que
nous adopterons jusqu'à la fin de ce chapitre.\\
\label{distribu}
On se donne alors une section distributionnelle
$u_{1}\in \Gamma^{-\infty}(X,E_{n})$ telle que
:
$$u_{1}=\delta_{z_{0}} e_{1},$$
où $\delta $  est la  distribution de Dirac
définie ci-dessus. Nous allons maintenant
construire une approximation de la  section
holomorphe $\hat\sigma_{1}$ correspondant à la
projection orthogonale  pour $h$ de $u_{1}$. Cette
projection existe car  on projette sur le
sous-espace vectoriel fermé $H^{0}(X,E_{n})$. Par
définition du noyau,(cf. \cite{wells} chap. 4), on
peut remarquer également que
:
$$P_{n}u_{1}(z)=\int_{X}\pr_{n}(z,z_{1})\delta_{z_{0}} e_{1}(z_{1})dvol(z_{1})
=\pr_{n}(z,z_{0}) e_{1}(z_{0}),$$
avec $\pr_{n}$ le noyau de $P_{n}$. Nous prenons
une section qui vaut $(1,0)$ dans la base
$\{e_{1},e_{2}\}$ pour $z_{0}=0$, sous la forme
:$$f_{1}(z)= e_{1}(z).$$

 Cette section, $f_{1}$ (dont la notation n'a aucun lien avec celle définie
  dans les paragraphes précédents) est une fonction ``cloche" . Elle
  va nous servir à la construction
de la projection orthogonale de $u_{1}$ pour $h$
sur $H^{0}(X,E_{n})$.
 %Si $\hat \sigma^{1}$ est cette projection on a :
 %$$<\hat \sigma^{1}-u_{1},\eta>_{L^{2}(h)}=0\quad \forall \eta\in H^{0}(X,E_{n}).$$
\section{Extension à  $X$.}
  On définit ensuite une fonction plateau $\phi $ telle que les
dérivées soient portées par un anneau de centre
$z_{0}$, de rayon intérieur $r'=\frac{1}{2}$ et de
rayon extérieur $R=1$, noté $A(r',R)$, tels que
$$D(z_{0},R)\in U.$$
Pour cela soit :
$$
\tilde \phi (x)=\left\lbrace
\begin{array}{l}
1 \quad \forall x\leq \frac{1}{2}\\ 0\quad \forall
x\geq 1.\\
\end{array}
\right.
$$
On pose ensuite :
$$\phi (z)=\tilde \phi(R^{-1}\vert z\vert ),\, \forall z\in U.
$$Il s'ensuit alors,  après avoir
posé $\tilde f_{1}=\phi.f_{1}$ qui est une section
alors $\cin$:
$$
\bar\partial(\tilde f_{1})(z)=\bar\partial (\phi.f_{1})(z)=\bar\partial
\phi.f_{1}(z)= \bar\partial
\phi_{1}e_{1}(z).\\
$$

Nous effectuons maintenant le calcul de la norme
$L^{2}(h)$ de la dérivée de $\tilde f_{1}$. Nous
utiliserons ce calcul dans le paragraphe suivant.
  \begin{eqnarray*}
\normlh{\bar\partial (\phi.f_{1})}&=& \normlh{\bar\partial
\phi.e_{1}(z)}\\
&\leq& \norml{\bar\partial
\phi}\normca{ e_{1}(z)}\\
&\leq& C_{1}\exp(-n r\acute {\,}
^{\,2})\norml{\bar\partial
\phi}
 \end{eqnarray*}
 %%%%%%%%%%%%%%%%%%%%%%%
 \section{Propriété des sections holomorphes}\label{poiu}
Ce paragraphe établit deux résultats utilisés dans
le calcul  de la projection de $u_{1}$ sur
$W_{n}$. En effet nous voulons modifier
$\sigma_{1}$ pour définir cette projection.

\subsection{Première propriété}
Soit $\eta\in H^{0}(X,E_{n})$. Nous voulons
obtenir une majoration  de $\vert
\eta(z)\vert_{h}$ pour $z\in X$. La proposition
\ref{estimate1} donne   une estimation suffisante
pour la suite des  calculs.
\\ Nous avons besoin, pour la démonstration de cette proposition, du lemme suivant :
\begin{lem}Soit $C_{4}>0$.
\label{trois}
Il existe $C_5$  telle que pour tout $n\geq 1$, et
pour tout fonction holomorphe $f$ sur
$D(0,r')\subset \C$, on ait
$$
|f'(0)| \leq C_5 n \left(
\int _{D(0,r')}|f(z)|^2 \exp (-C_{4}n|z|^2)dz d\overline{z}
\right) ^{\frac{1}{2}}.
$$
\end{lem}
\begin{demo} On considère
l'espace des fonctions holomorphes sur $D(0,r')$
avec la métrique $L^2$ pond\'er\'ee par
$\exp(-C_{4}n|z|^2)$. Les $z^i$ forment une base
orthogonale . On peut alors
\'ecrire
$$
f(z) = \sum _i a_i z^i.
$$
On obtient :
$$
f'(0) = a_1,
$$
et
$$
\int _{D(0,r')}|f(z)|^2 \exp (-C_{4}n|z|^2)dz d\overline{z}
= \sum _i |a_i|^2 \int _{D(0,r')}|z|^{2i} \exp (-C_{4}n|z|^2)dz d\overline{z}.
$$
En particulier, on a
$$
|a_1|^2 \int _{D(0,r')}|z|^{2} exp
(-C_{4}n|z|^2)dz d\overline{z}
\leq \int _{D(0,r')}|f(z)|^{2} exp (-C_{4}n|z|^2)dz d\overline{z} .
$$
De plus $ \exists\, C_{6}>0$
$$\int _{D(0,r')}|z|^{2} exp (-C_{4}n|z|^2)dz d\overline{z}\geq C_{6}n^{-2}
,
$$
d'o\`u le résultat. Ce lemme est utilisé dans la
proposition suivante.
\end{demo}\\

\begin{prop}\label{estimate2}Soit $\eta\in H^{0}(X,E_{n})$ telle que $\eta(z_{0})=0$. On a alors :
$\exists \, C_{7}$ ,$\exists n_{0}$  tel que
\\$n\geq n_{0}$,$\forall z\in X$,
\begin{equation}\label{estimate1}
  \vert \eta(z)\vert _{h}\leq C_{7}nd(z,z_{0})\normlh{\eta}.
\end{equation}

\end{prop}
\begin{demo}
 On  suppose  $\eta(z_{0})=0$. Soit $z\in D(z_{0},r^{'})$
Notons par $\nabla _n$ la connexion  de Chern de
$E_{n}$.
 On a alors :
$$
| \eta (z) |_{h} \leq \int _{y= z_0}^{y=z} |
\nabla _n(\eta )(y) |_{h} d\ell
$$
o\`u $d\ell $ est l'\'el\'ement de longueur
g\'eod\'esique le long du chemin reliant $z_0$ \`a
$z$. Il suffit de montrer que pour $i=1,2$ et
$D(y,r'')\subset D(z_0, r^{'})$

\begin{equation}\label{sdf}
  \vert
\nabla_{n}\eta
_{i}(y)\vert_{h}
\leq C_{6}n (\int_{D(y,r'')}\vert
\eta
_{i}(z)\vert\exp (-C_{4} n | z|^2) dz d\overline{z})^{\frac{1}{2}}.
\end{equation}

Pour cela on applique le lemme \ref{trois}  à
$\eta_{i}$ en supposant $y=0$ et au disque
$D(y,r'')$.\\

%Si $\eta(z_{0})\not=0$, alors la démonstration
%donne :
%$$\vert \eta(z)-\eta(z_{0})\vert _{h} \leq C nd(z,z_{0})
%\normlh{\eta-\eta(z_{0})}\leq C_{3} nd(z,z_{0})(\normlh{\eta}+\normlh{\eta(z_{0})}) $$
%D'où
%$$\vert \eta(z)\vert_{h}\leq C_{3} nd(z,z_{0})(\normlh{\eta}+\normlh{\eta(z_{0})})
%+\vert \eta(z_{0})\vert .$$ Or $\vert
%\eta(z_{0})\vert_{h}\geq  n\normlh{\eta(z_{0})}$ en utilisant la formule \ref{merthes} et
%$\vert \eta(z_{0})\vert_{h}\leq n \normlh{\eta}$.
%La proposition est donc vérifiée.
\end{demo}
%Remarquons
%également que la formule \ref{sdf} et le lemme
%\ref{trois} nous donne également la majoration
%suivante :
%\begin{equation}\label{ghj}
%  \vert \eta '(z)\vert _{h}\leq C_{3}nd(z,z_{0})\normlh{\eta}.
%\end{equation}
%De plus on a?????????????????
%%%%%%%%%%%%%%%%
%p. 87 a la fin de la partie 8.6.1: donner a partir de la prop. 8.6.2,
%une majoration pour les fonctions $\eta _1(z)$ et $\eta _2(z)$.
%Pour cela il convient de remarquer que la partie $h_{E_0}$ est fixe,
%i.e. on a pour une section $(f_1,f_2)$ de $E_0$ que
%$$
%|f_i| \leq C |f|_{h_{E_0}}.
%$$
%Donc on a (pour $i=1,2$)
%$$
%ùù\|\eta _i(z)\| \cdot h_{\ox (n)} ^{1/2} \leq C \|
%\eta (z) \| _h .
%$$
%Ceci donne, donc:
%$$
%|\eta _i(z)| \leq Ch_{\Oo (n)} ^{-1/2}
%$$
%
%On peut remarquer par ailleurs que cette notation
%qui veut que $h_{\Oo (n)}$ est la fonction reel
%$|\epsilon |_{h_{\Oo (n)}}$, est un abus (car
%normalement $h_{\Oo (n)}$ est la metrique sur $\Oo
%(n)$). Ca va mais il vaut mieux remarquer quelque
%part (avant) qu'on fera cet abus.
On utilisera à plusieurs fois cette proposition
combinée au résultat suivant :
\begin{prop}\label{fghj}
Soit $\eta\in H^{0}(X,E_{n})$ et $z\in X$. On a
alors
: $\exists \,C>0$
$$\vert \eta_{i}(z)\vert\leq C h^{\frac{-1}{2}}_{\ox(n)}(z)\vert \eta(z)\vert_{h}\mbox{ pour }i=1,2.$$
\end{prop}
\begin{demo}
En effet on a :
$$\vert \eta_{i}(z)\vert\leq C.\vert \eta(z)_{}\vert h^{\frac{1}{2}}_{E_{0}}(z).$$
car  $h_{E_{0}}(z) $ est indépendant de de $n$. On
a donc :
$$\vert \eta_{i}(z)\vert h^{\frac{1}{2}}_{\ox(n)}(z)
\leq C.\vert \eta(z)_{}\vert h^{\frac{1}{2}}_{E_{0}}(z)h^{\frac{1}{2}}_{\ox(n)}(z)
\leq  C.\vert \eta(z)_{}\vert_{h}. $$

\end{demo}
\subsection{Deuxième propriété}\label{xcvb}
Soit $\eta\in H^{0}(X,E_{n})$ telle que
$\eta(z_{0})=0$. On l'écrit de la façon suivante :
$$\begin{pmatrix}
\eta_{1}(z)&\\ \eta_{2}(z)&\\
\end{pmatrix}$$ dans le repère
$\{e_{1},e_{2}\}$ local sur $E_{{n}}$ défini dans
un voisinage de zéro. Supposons également que
$\normlh{\eta}=1$. Nous voulons estimer
\begin{equation}\label{nbv}
  \int_{D(0,r')}
\eta_{1}(z)h_{\ox(n)}(z)dzd\bar z.
\end{equation}
Le corollaire \ref{nhg} nous donne l'estimation
souhaitée.\\ Dans ce paragraphe la lettre  $C$
désignera  de façon  générique une constante
indépendante de $n$. Nous utiliserons le lemme
suivant
:
\begin{lem}\label{dix}Soit $\mathcal L$ un fibré en droite sur $X$. Notons $k$ une métrique
sur $\mathcal L$. Soit $s$ une section holomorphe
de $\mathcal L$. Posons $g=k(s,s)$. On a alors :
\begin{equation}\label{onze}
  \nabla (\frac{s}{\normk{s}})=\frac
  {1}{2}\dc \ln g . \frac{s}{\normk{s}}.
\end{equation}

\end{lem}
\begin{demo} Comme
$\partial s=\frac{\partial g}{g}s,$ on a :
  \begin{eqnarray*}\nabla
\frac{s}{\normk{s}}&=&\frac{\partial s}{\normk{s}}-d
\ln \normk{s}. \frac{s}{\normk{s}}\\
&=& \frac{\partial
g}{g}\frac{s}{\normk{s}}-\frac{1}{2}d\ln g
\frac{s}{\normk{s}}\\&=&\frac{1}{2}\dc \ln g
 \frac{s}{\normk{s}}.  \end{eqnarray*}

\end{demo}
 Nous rendons tout d'abord "symétrique"
 $h_{\ox(n)}$. Pour cela soit une nouvelle
coordonnée(non holomorphe) $u(z)$ telle que
:\\ $u(z_{0})=0$ et  $h_{\ox(1)}=\exp(-\vert
u(z)\vert
^{2})$.
%Quitte à faire un changement de coordonnées
%holomorphe,(comme dans le paragraphe \ref{metr}),
%on peut supposer que
%:
%$$h_{\ox(1)}=\exp(-\vert
%z\vert^{2}+\beta  \vert z\vert^{4}+O(\vert
%z\vert^{5})).$$ Dans ce cas
%:
$$\vert u(z)\vert^{2}=\vert z\vert^{2}-\beta \vert z\vert^{4}+O(\vert
z\vert^{5})$$\label{hjkn}
 Donc  $u(z)$ est de
la forme
:
$$z=(1+\beta ' u(z)\bar u(z)+O(\vert u(z)\vert
^{3}))u(z),$$ avec$\beta '$ ne dépendant que de
$\beta$. Alors on a
:
$$dzd\bar z=(1+\tilde \beta u\bar u+\psi(u))dud\bar u,$$
où $\tilde \beta$ ne dépend que de $\beta$,\\ avec
$\psi(u)=O(\vert u\vert^{3})$.
 Nous
allons maintenant calculer l'intégrale définie par
la formule \ref{nbv}. Nous rappelons (cf.
proposition \ref{fghj}, formule  \ref{estimate1}
et
 $\normlh{\eta}=1$) que
:
\begin{equation}\label{trente}
  \vert \eta_{i}(z)\vert\leq C nh^{\frac{-1}{2}}_{\ox(n)}(z) \mbox{ pour }i=1,2.
\end{equation}

 Considérons
maintenant $v(z)=(1+\beta ' u(z)\bar u(z))u(z)$.
On a alors

$$         \eta^{}_{1}(z)=\eta^{}_{1}(v(z))+\int_{v(z)}^{z}\dt \eta_{1}(t)dt.$$

Notons que $arg(z)=arg(v(z)).$ Le chemin entre les
deux points est donc le rayon allant de $z$ à
$v(z)$. Nous avons :
\begin{equation}\label{douze}
\vert h^{\frac{1}{2}}(z)\eta_{1}(z)-h^{\frac{1}{2}}(v(z))\eta_{1}(v(z))
\vert \leq \int_{z}^{v(z)} \vert\dt h^{\frac{1}{2}}(t) \eta_{1}(t)\vert
dt.
\end{equation}
Ici $\dt$ est la dérivée directionnelle  dont
l'axe est le rayon allant de $z$ à $v(z)$.
Calculons $\dt h^{\frac{1}{2}}(t)
\eta_{1}(t)$. En utilisant le repère construit dans le paragraphe \ref{seize}
on a $\eta_{1}=<\eta, \check f_{1}\otimes
\check\epsilon^{\otimes n}>_{h}$ où $\check f_{1}$(resp. $\check \epsilon$) est l'élément dual
de $f_{1}$(resp. $ \epsilon$). De plus
$h_{\ox(n)}^{\frac{1}{2}}=\normhox{\epsilon}^{n}=\normhox{\check\epsilon}^{-n}$
d'où
$$h_{\ox(n)}^{\frac{1}{2}}\eta_{1}=<\eta, \check f_{1}\otimes \frac{\check\epsilon^{\otimes n}}{\normhox{\check\epsilon}^{n}}
>_{h}.$$
On obtient alors :
\begin{eqnarray*}\label{vingt}
\dt \hox^{\frac{1}{2}}(t)
\eta_{1}(t)&=&<\nabla_{\dt}\eta, \check f_{1}\otimes \frac{\check\epsilon^{\otimes n}}{\normhox{\check\epsilon}^{n}}
>_{h}\\&+&<\eta, \nabla_{\dt}\check f_{1}\otimes \frac{\check\epsilon^{\otimes n}}{\normhox{\check\epsilon}^{n}}
>_{h}\\&+&n<\eta, \check f_{1}\otimes \nabla_{\dt}\frac{\check\epsilon}{\normhox{\check\epsilon}}
\otimes \frac{\check\epsilon^{\otimes n-1}}{\normhox{\check\epsilon}^{n-1}}
>_{h}.
\end{eqnarray*}
Nous étudions maintenant chacun des trois termes
du membre de droite. On a :
\begin{itemize}
  \item[$\bullet$]
  $<\nabla_{\dt}\eta, \check f_{1}\otimes \frac{\check\epsilon^{\otimes n}}{\normhox{\check\epsilon}^{n}}
>_{h}\leq \normhsi{\nabla_{\dt}\eta}\normhsi{\check f_{1}\otimes \frac{\check\epsilon^{\otimes n}}{\normhox{\check\epsilon}^{n}}
}\leq C \normhsi{\nabla_{\dt}\eta}. $ La formule
 \ref{sdf} nous donne la
majoration
:
$$<\nabla_{\dt}\eta, \check f_{1}\otimes \frac{\check\epsilon^{\otimes n}}{\normhox{\check\epsilon}^{n}}
>_{h}\leq C .n.$$

  \item[$\bullet$]Nous avons également :
$$<\eta, \nabla_{\dt}\check f_{1}\otimes \frac{\check\epsilon^{\otimes n}}{\normhox{\check\epsilon}^{n}}
>_{h}\leq C \normhsi{\eta}\normhsi{\check f_{1}\otimes \frac{\check\epsilon^{\otimes n}}
{\normhox{\check\epsilon}^{n}} } \leq C
\normhsi{\eta}.$$  On applique la proposition
\ref{estimate2} pour obtenir :
$$<\eta, \nabla_{\dt
}\check f_{1}\otimes \frac{\check\epsilon^{\otimes n}}{\normhox{\check\epsilon}^{n}}
>_{h}\leq C .n.$$

  \item[$\bullet$]Enfin,  pour le dernier terme, on utilise le lemme \ref{dix}:
  $$\nabla\frac{\check\epsilon}{\normhox{\check\epsilon}}
=\frac{-1}{2}\dc\ln h.\frac{\check\epsilon}{\normhox{\check\epsilon}}.
$$
Donc
$$\nabla_{\dt}\frac{\check\epsilon}{\normhox{\check\epsilon}}=
\iota_{\dt}(\frac{-1}{2}\dc \ln h)\frac{\check\epsilon}{\normhox{\check\epsilon}},$$
où $\iota$ désigne l'opération de contraction avec
le champ de vecteur  $\dt$. Mais $\dc
\ln h$ est une forme angulaire à $O(\vert z\vert ^{4})$ près(cf. p.\pageref{hjkn}). Donc
$$\iota_{\dt}(\frac{-1}{2}\dc
\ln
h)=0.$$ Il s'ensuit :

$$n<\eta, \check f_{1}\otimes \nabla_{\dt}\frac{\check\epsilon}{\normhox{\check\epsilon}}
\otimes \frac{\check\epsilon^{\otimes n-1}}{\normhox{\check\epsilon}^{n-1}}
>_{h}=0.$$
\end{itemize}
On a donc :
$$\vert\int_{v(z)}^{z} \dt h_{\ox(n)}^{\frac{1}{2}}\eta_{1}(t)
dt\vert \leq Cn\vert z-v(z)\vert\leq  C  n
\vert u\vert^{4}.$$
Et $$\int_{D(0,r')}\vert u\vert^{4}\exp(-{n}
\vert u\vert ^{2})(1+\tilde
\beta u\bar u+O(\vert u\vert^{3}))dud\bar
u=O(\negli{3}).$$ Donc on a
$$\int_{D(0,r')}h_{\ox(n)}(\int_{v(z)}^{z}\dt h_{\ox(n)}^{\frac{1}{2}} \eta_{1}(t)dt)dzd\bar z=O(\negli{2}).$$
 Il reste donc à
étudier le terme en $\eta(v(z))$. Pour cela on a
la proposition suivante :
\begin{prop}Soit $f$ une fonction holomorphe telle que $f(z_{0})=0$, alors

\begin{equation}\label{rtyui}
  \int_{D(0,r')}f(v)
\exp(-n \vert u\vert ^{2})(1+\tilde \beta u\bar u))dud\bar u=0.
\end{equation}

\end{prop}
\begin{demo}
Comme $f$ est une fonction holomorphe, on a :
$$f(z)=\sum_{j\geq 1}a_{j}z^{j}$$
D'où
$$f(v)=\sum_{j\geq 1}a_{j}(1+\beta ' u(z)\bar
u(z))^{j}u(z)^{j}=\sum_{j\geq 1}a_{j}g(\vert
u\vert )u(z)^{j}.$$ L'expression \ref{rtyui} est
alors évidente .
\end{demo}
Cette proposition donne alors le corollaire
suivant :
\begin{cor}\label{nhg}
$$\int_{D(0,r')}
\eta_{1}(z)h_{\ox(n)}(z)dzd\bar z=O(\negli{2}).$$
\end{cor}
\begin{demo}
On applique la proposition  précédente à la
section $\eta_{1}$. Il reste donc à montrer que :
\begin{equation}\label{plj}
  \int_{D(0,r')}
\eta_{1}(z)(\exp-n\vert u\vert^{2})\psi(u)dzd\bar z=O(\negli{2}).
\end{equation}

Or $\psi(u)=O(\vert u\vert^{3})$ est indépendant
de $n$ et $$\vert \eta_{1}(z(u))\vert \leq C
n\vert u\vert \exp(n\frac{\vert u\vert
^{2}}{2}),$$ donc $
\int_{D(0,r')}\eta_{1}(z)(\exp-n\vert
u\vert^{2})\psi(u)dzd\bar z$ est majorée par
$$Cn\int_{D(0,r')}\vert u\vert^{4}
(\exp-n\frac{\vert u\vert^{2}}{2})dzd\bar z.$$ Ce
qui donne l'équation \ref{plj}.
\end{demo}
On obtient de même :
\begin{cor}\label{nhg1}
$$\int_{D(0,r')}
\eta_{1}(z)h_{\ox(n)}(z)z\bar zdzd\bar z=O(\negli{2}).$$
\end{cor}
 %%%%%%%%%%

\section{Holomorphisation de $\tilde f_{1}$}
 Nous avons donc maintenant une  section globale, $\tilde f_{1}$, mais
$\cin$, de $E_{n}$. On utilise alors la
proposition \ref{delta}, pour obtenir des objets
holomorphes. On a
:
$$\Delta_{\db}=\Delta^{0,\bullet}+\phi(R)-iR_{\R,h}({E_{n}}),$$
avec $\phi(R)$ un opérateur utilisant le tenseur
de Ricci de $X$. L'opérateur de Green $Gr$ a une
norme $L^{2}$ pour $h$
 qui est comparable à
$n^{-\frac{1}{2}}$(cf. \cite{don5}).
 La section $\cin$, $\tilde f_{1}$ de $E_{n}$, va
 être modifiée  pour avoir une section holomorphe. On
définit pour cela :
$$\zeta_{1}=-\db^{*}Gr\db \tilde f_{1}.$$
On obtient alors la première section globale
holomorphe $\sigma^{1}=\tilde f_{1}+\zeta_{1}$. De
plus :
\begin{eqnarray*}
\normlh{\zeta_{1}}^{2}&=&<\db^{*}Gr\db
\tilde f_{1},\db^{*}Gr\db \tilde f_{1}>_{h}
\\&=& <Gr\,\db
\tilde f_{1},\db \tilde f_{1}>_{h}
\\&\leq&
Cn^{-1}\normlh{\db \tilde f_{1}}^{2},
\end{eqnarray*}
d'où :
\begin{eqnarray}\label{eee} \normlh{\zeta_{1}}^{2}&\leq&n^{-1}C_{1}\exp(-2nr
^{2})\norml{\bar\partial
\phi}^{2}
.
\end{eqnarray}

Or  pour tout $a\in \N^{*}$ :
\begin{equation}\label{appp}
\quad \exp(-n r\acute {\,}
^{\,2})=o(\frac{1}{n^{a}}).
\end{equation}

Par la suite nous pourrons donc négliger de cette
manière les termes en $\exp(-nr
^{2}).$

 On a la proposition suivante
:
\begin{prop}\label{sigmaa}
Il existe une section  holomorphe que l'on note
$\sigma_{1}$ telle que :
\begin{equation*}
  \sigma^{1}(z)=
  \begin{cases}
  e_{1}(z)+\zeta_{1}(z)&  \forall z\in D(0,r')\\
  \phi(z)e_{1}(z)+\zeta_{1}(z)& \forall z\in A(R,r')\\
  \zeta_{1}(z)& \forall z\in X\backslash
D(0,R),
  \end{cases}
\end{equation*}
avec $$\normlh{\zeta_{1}}=o(\negli{a})\, ,\quad
\forall a\in
\N^{*}.$$ On a alors :
$$ \normlh{\sigma^{1}}^{2}=\frac{\pi}{n}-\frac{\pi d_{1}+ \pi r}{n^{2}}+O(\frac{1}{n^{3}}).$$
Ici $r$  ne dépend que de $X$(cf. formule
\ref{vol} p. \pageref{vol}).
\end{prop}
%%%%%%%%%%%%%%%%%%%%%%%
\begin{demo}
Nous considérons la forme volume définie par la
formule \ref{vol} p. \pageref{vol} c'est-à-dire la
forme volume est de la forme :
\begin{equation*}
  \omega(z)=(1+ r z\bar z +...)dz\wedge d\bar z=w(z)dz\wedge d\bar z.
\end{equation*}
 Nous avons
:
\begin{eqnarray*}
 \normlh{\sigma^{1}}^{2}&=&\normlh{\tilde f_{1}}^{2}+
 \normlh{\zeta_{1}}^{2}+<\tilde f_{1}(z),\zeta_{1}(z)>_{L^{2}(h)}
 +<\zeta_{1}(z),\tilde f_{1}(z)>_{L^{2}(h)}.\\
\end{eqnarray*}
\begin{enumerate}
  \item  Pour le premier terme du deuxième membre on a :
  $$<\tilde f_{1}(z),\tilde
f_{1}(z)>_{L^{2}(h)}=\int_{D(0,r')}h^{11}(z)w(z)dzd\bar
z+\int_{A(r',R)}<\tilde f_{1}(z),\tilde
f_{1}(z)>_{h}w(z) dzd\bar z.$$

  Comme
  $$\int_{A(r',R)}<\tilde f_{1}(z),\tilde f_{1}(z)>_{h}w(z) dzd\bar z=o(\negli{a}),$$
  et puisque $$h^{11}(z)= \exp(-n\mo{z}^{2})(1
-d_{1}
z\bar z +{\bf {A}}_{ 11}z\bar z^{2}+{\bf {B}}_{
11}z^{2}\bar z+
\varphi (z)),$$
On a, d'une part :
 $$ \int_{D(0,r')}exp(-n\mo{z}^{2})w(z)dzd\bar z=\int_{0}^{2\pi}
 \int_{0}^{r'}\exp(-n\rho^{2})w(\rho,\theta)\rho d\rho d\theta=\frac{\pi}{n}
 +\frac{r\pi}{n^{2}}+o(\negli{2})$$
 et d'autre part, par une intégration par partie, on a :
  \begin{eqnarray*} -\int_{D(0,r')}d_{1}\exp(-n\mo{z}^{2})w(z)
z\bar zdzd\bar z&=&-
\int_{0}^{2\pi} \int_{0}^{r'}d_{1}\exp(-n\rho^{2})w(\rho,\theta)\rho^{3} d\rho d\theta\\
&=&
-d_{1}\frac{\pi}{n^{2}}+o(\negli{2}).
\end{eqnarray*}
De même il est clair que
$$\int_{D(0,r')}exp(-n\mo{z}^{2}){\bf {A}}_{ 11}z\bar z^{2}w(z)dzd\bar z=0+O(\negli{3}),$$
et
$$\int_{D(0,r')}exp(-n\mo{z}^{2}){\bf {B}}_{11}z^{2}\bar zw(z)dzd\bar z=0+O(\negli{3})$$
Le dernier terme, $\int_{D(0,r')}
exp(-n\mo{z}^{2})\varphi (z)w(z)dzd\bar z$, vu les
conditions sur $h$ et vu $w(z)$, est en
$O(\frac{1}{n^{3}})$ car il est de la forme :
$\int_{D(0,r')}\psi_{1}(z) h_{\ox(n)}dzd\bar z$
avec $\psi_{1}(z)=O(\vert z\vert ^{4})$

  \item  Le deuxième terme vérifie $\normlh{\zeta_{1}}^{2}=o(\negli{a})$.

  \item Les  deux derniers termes  vérifient en utilisant
l'inégalité \ref{eee}
:
$$\vert<\zeta_{1}(z),\tilde f_{1}(z)>\vert=
\vert<\zeta_{1}(z),\tilde f_{1}(z)>
\vert\leq \normlh{\zeta_{1}}.\normlh{\tilde f_{1}}=o(\negli {a})$$

\end{enumerate}

On obtient le  résultat désiré.

\end{demo}

\section{Orthogonalisation de $\sigma^{1}$ par rapport à un sous espace vectoriel  $W^{1}_{n}$ de $W_{n}$}

\label{milkaa}
Nous rappelons que $W_{n}$ est de dimension notée
$p$. Montrons maintenant que la  section
$\sigma^{1}$ correspond ``assez bien" à la section
holomorphe
 globale cherchée. En effet, elle ne correspond pas pour l'instant
 à la projection orthogonale induite par la métrique  $h$ de la section $u_{1}$. Pour cela nous devons la modifier
 de façon à ce qu'elle définisse une projection orthogonale sur $H^{0}(X,E_{n})$. De ce fait on
étudie une perturbation de $\sigma^{1}$. Pour
commencer nous prenons une famille de $p-2$
sections orthonormées que nous notons $\eta^{i}$
pour $i=1..p-2$. Nous choisissons  $p-2$ sections
de telle sorte que : $
\eta^{i}(z_{0})=0\mbox{ pour }i=1..p-2.$
%Elles forment une famille notée $\mathcal {F}$.
%Quant à la section $\eta^{p-1}$, elle vérifie
%seulement que sa composante selon $e_{1}$ s'annule
%en $z_{0}$.
 La réunion de ces sections forment une
famille, appelée $\mathcal F$, de $p-2$ sections
orthonormées de
 $H^{0}(X,E_{n})$.
  Elles engendrent un  espace vectoriel de dimension $p-2$ que nous appellerons $W^{1}_{n}$.
  Dans le paragraphe suivant, nous calculons le "défaut" d'orthogonalité de la section
  $\sigma^{1}$ avec l'espace vectoriel engendré par la famille $\mathcal {F}$.
 % Le paragraphe \ref{vfe} correspond au calcul du "défaut" d'orthogonalité de la section
 % $\sigma^{1}$ avec la droite vectorielle engendrée par la section $\eta^{p-1}$.

\subsection{``Défaut"  d'orthogonalité de $\sigma^{1}$ avec
$Vect(\mathcal{F})$} \label {hgfd}Nous prenons
$p-2$ sections orthonormées  notées $\eta^{i}$
telles que pour $i=1..p-2$:
$$\begin{pmatrix}
\eta^{i}_{1}(z)&\\ \eta^{i}_{2}(z)&\\
\end{pmatrix}=
\begin{pmatrix}
\eta^{i}_{11}z+...+\eta^{i}_{1m}z^{m}+o(z^{m})&\\
\eta^{i}_{21}z+...+\eta^{i}_{2m}z^{m}+o(z^{m})&\end{pmatrix},$$
dans le repère $\{e_{1},e_{2}\}$ local sur
$E_{{n}}$ défini dans un voisinage de zéro.

 Elles vérifient $\eta^{i}(z_{0})=0$.

 Calculons maintenant le ``défaut"
d'orthogonalité de $\sigma^{1}$ avec
$vect(\mathcal{F}_{0})$.
\\
  Soit $i \in \{1..p-2\}$.
  \begin{eqnarray*}
\int_{X}< \sigma^{1},\eta^{i}>_{h}w(z)dzd\bar z&=&\int_{X}<\tilde f_{1}+\zeta_{1},
\eta^{i}>_{h}w(z)dzd\bar z\\
&=&\int_{X}<\tilde f_{1},\eta^{i}>_{h}w(z)dzd\bar
z+\int_{X}<\zeta_{1},\eta^{i}>_{h}w(z)dzd\bar z.
\end{eqnarray*}
Or
\begin{eqnarray*}
\vert\int_{X}<\zeta_{1},\eta^{i}>_{h}w(z)dzd\bar z\vert&\leq&\normlh{\zeta_{1}}.\normlh{\eta^{i}}\\
&\leq&o(\frac{1}{n^{a}})\mbox{ cf . p
\pageref{appp}}.
\end{eqnarray*}

De plus  pour $1\leq i\leq p-2$
\begin{eqnarray}\label{tilde1}
\int_{X}<\tilde f_{1},
\eta^{i}>_{h}w(z)dzd\bar z&=&\int_{D(0,r')}h^{11}(z)
\eta^{i}_{1}w(z)dzd\bar z+\int_{D(0,r')}h^{12}(z)
\eta^{i}_{2}w(z)dzd\bar z\nonumber\\
&+&\int_{A(r',R)}<\tilde f_{1},
\eta^{i}>_{h}w(z)dzd\bar z\nonumber.\\\nonumber
\end{eqnarray}
On trouve immédiatemment :
$$ \int_{A(r',R)}<\tilde f_{1},
\eta^{i}>_{h}w(z)dzd\bar z=o(\negli{a}).
$$
Etudions maintenant les deux intégrales suivantes.
Nous utilisons la même méthode dans les deux cas.
Nous les décomposons en somme d'intégrales où
chacune correspond au développement de
$h^{11}_{E_{0}}(z)w(z)$ ou $h^{12}_{E_{0}}(z)w(z)$
suivant le cas.
\begin{eqnarray}
\int_{D(0,r')}h^{11}
\eta^{i}_{1}(z)w(z)dzd\bar z&&=
 \int_{D(0,r')} \eta^{i}_{1}(z)h_{\ox(n)}(z)h^{11}_{E_{0}}(z)w(z)dzd\bar z\nonumber\\
=&&\int_{D(0,r')}
\eta^{i}_{1}(z)h_{\ox(n)}(z)dzd\bar
z\label{14}\\ +&&\int_{D(0,r')}
\eta^{i}_{1}(z)h_{\ox(n)}(z)(r+d_{1})z\bar z dzd\bar z\label{15}\\
+&&\int_{D(0,r')}
\eta^{i}_{1}(z)h_{\ox(n)}(z)O(\vert z\vert ^{3})dzd\bar z\label{16}
\end{eqnarray}
Ici   $O(\vert z\vert ^{3})$ est indépendant de
$n$. De même on a \begin{eqnarray}
\int_{D(0,r')}h^{12}
\eta^{i}_{2}(z)w(z)dzd\bar z&=
& \int_{D(0,r')} \eta^{i}_{2}h_{\ox(n)}(z)h^{12}_{E_{0}}(z)w(z)dzd\bar z\\
&=&\int_{D(0,r')}
\eta^{i}_{2}(z)h_{\ox(n)}(z)O(\vert z\vert ^{3})dzd\bar z\label{17}
\end{eqnarray}
Nous voulons maintenant obtenir des estimations
sur toutes ces intégrales. Pour ce faire, nous
utilisons le paragraphe \ref{poiu} et nous
rappelons que nous avons la majoration suivante
(cf. proposition \ref{fghj}, proposition
\ref{estimate2} et $\normlh{\eta}=1$):
$$\vert\eta^{i}_{j}(z)\vert \leq Cn h_{\ox(n)}^{\frac{-1}{2}}\mbox{ pour }i=1,2\,.$$ Nous
étudions de la même façon les intégrales \ref{16}
et \ref{17}. En effet,  en utilisant la
 la proposition \ref{estimate2} appliquée à $\eta^{i}_{1}(z)$ ou
$\eta^{i}_{2}(z)$ et formule \ref{merthes} , nous
obtenons :
\begin{eqnarray*}
\vert \int_{D(0,r')}
\eta^{i}_{1}(z)h_{\ox(n)}(z)O(\vert z\vert ^{3})dzd\bar z\vert &\leq&\int_{D(0,r')}
\vert \eta^{i}_{1}(z)h_{\ox(n)}(z)O(\vert z\vert ^{3})\vert dzd\bar z\\
&\leq &\int_{D(0,r')}
\vert \eta^{i}_{1}(z)\vert _{h} h^{\frac{1}{2}}_{\ox(n)}(z)O(\vert z\vert ^{3}) dzd\bar z\\
&\leq& C_{7} n\normlh{\eta^{i}}\int_{D(0,r')}
h^{\frac{1}{2}}_{\ox(n)}(z)O(\vert z\vert ^{3})
dzd\bar z\\ &=&O(\negli{2}).
\end{eqnarray*}
Ici $C_{7}$ est indépendante de $n$. On a donc :

$$\int_{D(0,r')}
\eta^{i}_{1}(z)h_{\ox(n)}(z)O(\vert z\vert ^{3})dzd\bar z=O(\negli{2})$$
et de la même façon :

$$\int_{D(0,r')}
\eta^{i}_{2}(z)h_{\ox(n)}(z)O(\vert z\vert ^{3})dzd\bar z=O(\negli{2}).$$
%%%%%%%%%%%%%%%%%%%%%%%%%%%%%%%%%%
D'après le paragraphe \ref{xcvb}, on a pour l'
intégrale \ref{14} :
$$\int_{D(0,r')}
\eta^{i}_{1}(z)h_{\ox(n)}(z)dzd\bar
z=O(\negli{2})\mbox{ (cf. corollaire \ref{nhg})
}$$ et de la même façon (cf. corollaire
\ref{nhg1}) l' intégrale \ref{14} vérifie
$$\int_{D(0,r')}
\eta^{i}_{1}(z)h_{\ox(n)}(z)(r+d_{1})z\bar z dzd\bar z=O(\negli{2}).$$
On obtient donc : $\forall \, 1\leq i\leq p-2$,
\begin{equation}\label{pml}
  <\sigma^{1},\eta^{i}>_{L^{2}(h)}=O(\negli{2}).
\end{equation}

%%%%%%%%%%%%%%%%%%%%%%%%%%%%%%%%%%%%%%%
\subsection{Construction de la projection $\check\sigma^{1}$ de
$u_{1}$ sur $W^{1}_{n}$} Une section holomorphe,
projection orthogonale pour $h$ de $u_{1}$ sur
$W_{n}$ vérifie
:

$$\forall \eta\in H^{0}(X,E_{n}) \quad <\check\sigma^{1},\eta>_{L^{2}(h)}=\eta_{1}(z_{0})$$
où $\eta_{1}$ est la première composante de $\eta$
dans le repère $\{e_{1},e_{2}\}$. Cela nous amène
à définir la projection $\check \sigma^{1}$ de
$u_{1}$ sur $W^{1}_{n}$, c'est à dire  $\check
\sigma^{1}$ vérifie
:
$$<\check\sigma^{1},\eta>_{L^{2}(h)}=0\quad \forall \eta\in W^{1}_{n}.$$
\begin{prop}\label{sigh1} Nous définissons une section holomorphe, projection  orthogonale
 de $\sigma^{1}$  sur $W^{1}_{n}$ par :

\begin{equation}\label{oij}
 \check\sigma^{1}=\sigma^{1}-\sum_{i=1..p-2}<\sigma^{1},\eta^{i}>_{L^{2}(h)}\eta^{i},
\end{equation}

avec
$$\normlh{\check\sigma^{1}}^{2}=\frac{\pi}{n}-\frac{\pi
d_{1}+\pi r}{n^2}+O(\frac{1}{n^{3}}).$$
$$\check\sigma^{1}(z_{0})=e_{1}(z_{0}).$$

Ici $r$  ne dépend que de $X$(cf. formule
\ref{vol} p. \pageref{vol}).
\end{prop}
\begin{demo}
On a en effet :
\begin{itemize}
  \item $\forall \eta^{i}\in \mathcal F$,
  \begin{eqnarray*}
  <\check\sigma^{1},\eta^{i}>_{L^{2}(h)}&=&
  <\sigma^{1},\eta^{i}>_{L^{2}(h)}-
\sum _{j=1..p-2}<\sigma^{1},\eta^{j}>_{L^{2}(h)}<\eta^{j},\eta^{i}>=0
  \end{eqnarray*}

  \item Nous utilisons la formule \ref{pml} et
  la proposition \ref{sigmaa} :
  \begin{eqnarray*}
  \normlh{\check \sigma^{1}}^{2}&=&< \sigma^{1}, \sigma^{1}>_{L^{2}(h)}
  -\sum_{i=1..p-2}\vert<\sigma^{1},\eta^{i}>\vert^{2}\\
&=&\frac{\pi}{n}-\frac{\pi d_{1}+\pi
r}{n^{2}}+O(\frac{1}{n^{3}})+
(p-2)(O(\negli{2}))^{2}
  \end{eqnarray*}
  et le fait que $p\sim n$.
  Pour la dernière propriété de $ \tilde \sigma^{1}$ on utilise l'équation \ref{oij}.
\begin{eqnarray*}
  \check\sigma^{1}(z_{0})=\sigma^{1}(z_{0})-\sum_{i=1..p-2}O(\negli{2})\eta^{i}(z_{0}),
\end{eqnarray*}
et le fait que $\eta^{i}(z_{0})=0$ pour
$i=1..p-2$.

\end{itemize}
\end{demo}
\subsection{Orthogonalisation de $\sigma^{2}$ par rapport à
 $W^{1}_{n}$}\label{che}

Nous effectuons le même travail que pour $u_{1}$
pour la section
$$u_{2}(z)=\delta_{z_{0}} e_{2}(z).$$
Nous résumons ici  les différentes étapes de la
construction. Nous avons :
$$P_{n}u(z)=\int_{X}\pr_{n}(z,z_{1})\delta_{z_{0}} e_{2}(z_{1})dvol(z_{1})
=\pr_{n}(z,z_{0}) e_{2}(z_{0}).$$
On considère maintenant la section
:$$f_{2}(z)=
e_{2}(z).$$  On obtient alors :
\begin{prop}
Il existe une section  holomorphe que l'on note
$\sigma^{2}$ telle que :
\begin{equation*}
  \sigma^{2}(z)=
  \begin{cases}
  e_{2}(z)+ \zeta_{2}(z) &\forall z\in D(0,r')\\
  \phi(z) e_{2}(z)+\zeta_{2}(z)& \forall z\in A(R,r')\\
  \zeta_{2}(z)& \forall z\in X\backslash
D(0,R),
  \end{cases}
\end{equation*}
avec
$$\zeta_{2}=-\db^{*}Gr\db \tilde f_{2}.$$
\end{prop}
Tout d'abord, nous construisons la section
holomorphe $\check\sigma^{2}$. Pour cela nous
considérons la famille $\mathcal F$ composée de
$p-2$ sections  orthonormées notées $\eta^{i}$
pour $i=\{1,...,p-2\}$ qui s'annulent en $z_{0}$.
(cf. p.
\pageref{milkaa}).
Les calculs sont identiques à ceux du paragraphe
\ref{milkaa}.
\begin{prop}\label{sigh2} Nous définissons une section holomorphe, projection  orthogonale pour $h$
de $\sigma^{2}$  sur $Vect(\mathcal F)$ par
:
$$\check\sigma^{2}=\sigma^{2}-\sum_{i=1..p-2}<\sigma^{2},\eta^{i}>\eta^{i}.$$
On a de plus
$$\normlh{\tilde\sigma^{2}}^{2}=\frac{\pi}{n}-\frac{\pi
d_{2}+\pi r}{n^2}+O(\frac{1}{n^{3}}),$$ et
$$\check\sigma^{2}_{2}(z_{0})=e_{2}(z_{0}),$$
 où $\check\sigma^{2}_{2}$ est la deuxième composante de
$\check\sigma^{2}$ dans le repère
$\{e_{1},e_{2}\}$.
\end{prop}

%%%%%%%%%%%%%%%%%%%%%%%%%%%%%%%%%%%%%%%%%%%%%%%%%%%%%%%
\section{Fin de l'orthogonalisation de $\check\sigma^{1}$ et $\check\sigma^{2}$ }
Nous complétons maintenant la famille $\mathcal F$
par deux sections de telle sorte que nous ayons
une base de $W_{n}$. Nous choisissons tout
naturellement les sections
 $\check\sigma^{1}$ et $\check \sigma^{2}$.
 En effet chacunes de ces  deux sections à une composante non nulle en $z_{0}$.
 Nous commençons par  normaliser $\check \sigma^{2}$. Nous effectuons le calcul pour la section
 $\check\sigma^{1}$. La démarche est identique pour $\check\sigma^{2}$.
 \subsection{Normalisation de  $\check{\sigma^{2}}$}\label{ghhhj}
 \begin{eqnarray*}
   <\check\sigma^{2},\check\sigma^{2}>_{L^{2}(h)}&=&
   <\sigma^{2}-\sum_{i=1..p-2}<\sigma^{2},\eta^{i}>\eta^{i},\sigma^{2}-\sum_{j=1..p-2}
   <\sigma^{2},\eta^{j}>\eta^{j}>_{L^{2}(h)}\\
&=&<\sigma^{2},\sigma^{1}>_{L^{2}(h)}-\sum_{i=1..p-2}<\sigma^{2},\eta^{i}><\sigma^{2},\eta^{i}>
\end{eqnarray*}
Il est évident que :
$$\sum_{i=1..p-2}<\sigma^{1},\eta^{i}><\sigma^{1},\eta^{i}>=O(\negli{3}),$$
car $p\sim n$ et
$<\sigma^{2},\eta^{i}>=O(\negli{2})$(cf.
paragraphe \ref{milkaa}). De plus
\begin{eqnarray*}
<\sigma^{2},\sigma^{1}>_{L^{2}(h)}&=&<e_{2}+\zeta_{2},e_{2}+\zeta_{2}>_{L^{2}(h)}\\
&=&<e_{2},e_{2}>_{L^{2}(h)}+2<\zeta_{2},e_{2}+\zeta_{2}>_{L^{2}(h)}
\end{eqnarray*}
Or (cf. proposition \ref{sigmaa})$$\vert
<\zeta_{2},e_{2}+\zeta_{2}>_{L^{2}(h)}
\vert\leq\normlh{\zeta_{2}}\normlh{e_{2}+\zeta_{2}}=O(\negli{a}).$$
On a également :
\begin{eqnarray}\label{mlkj}
<e_{2},e_{2}>_{L^{2}(h)}&=&\int_{X}h^{22}_{E_{0}}h_{\ox(n)}w(z)dzd\bar
z
=\frac{\pi}{n}+O(\negli{2})
\end{eqnarray}
On pose
$$\eta^{p-1}=\frac{1}{\normlh{\check{\sigma^{2}}}}\check{\sigma^{1}}.$$
 \subsection{''Défaut''  d'orthogonalité de
 $\check\sigma^{1}$ avec $Vect(\eta^{p-1})$}
 Nous calculons donc le ''défaut d'orthogonalité'' de
 $\check\sigma^{1}$ avec $Vect(\check \sigma^{2})$.
\begin{eqnarray*}
   <\check\sigma^{1},\check\sigma^{2}>_{L^{2}(h)}&=&
   <\sigma^{1}-\sum_{i=1..p-2}<\sigma^{1},\eta^{i}>\eta^{i},\sigma^{2}-\sum_{j=1..p-2}
   <\sigma^{2},\eta^{j}>\eta^{j}>_{L^{2}(h)}\\
&=&<\sigma^{1},\sigma^{2}>_{L^{2}(h)}-\sum_{i=1..p-2}<\sigma^{1},\eta^{i}><\sigma^{2},\eta^{i}>
\end{eqnarray*}
Il est évident que :
$$\sum_{i=1..p-2}<\sigma^{1},\eta^{i}><\sigma^{2},\eta^{i}>=O(\negli{3}),$$
car $p\sim n$ et
$<\sigma^{1},\eta^{i}>=O(\negli{2})$(cf.
paragraphe \ref{milkaa}). De plus
\begin{eqnarray*}
<\sigma^{1},\sigma^{2}>_{L^{2}(h)}&=&<e_{1}+\zeta_{1},e_{2}+\zeta_{2}>_{L^{2}(h)}\\
&=&<\zeta_{1},e_{2}+\zeta_{2}>_{L^{2}(h)}+
<e_{1}+\zeta_{1},\zeta_{2}>_{L^{2}(h)}
\end{eqnarray*}
Or (cf. proposition \ref{sigmaa})$$\vert
<\zeta_{i},e_{j}+\zeta_{j}>_{L^{2}(h)}
\vert\leq\normlh{\zeta_{i}}\normlh{e_{j}+\zeta_{j}}=O(\negli{a})\mbox{ pour }i,j=1,2\,.$$
On a donc
$$<\check\sigma^{1},\check\sigma^{2}>_{L^{2}(h)}=O(\negli{3}).$$
D'où
\begin{equation}\label{cdz}
  <\check\sigma^{1},\eta^{p-1}>_{L^{2}(h)}=O(\negli{2}).
\end{equation}
 \subsection{Projection de $u_{1}$ sur  $W^{1}_{n}+Vect(\eta^{p-1})$}\label{kkhj}
 On veut définir la projection orthogonale pour $h$ de $u_{1}$, notée $\tilde
\sigma^{1}$, sur le sous espace vectoriel engendré par $W^{1}_{n}$
et $Vect(\eta^{p-1})$.  On doit avoir pour cela
:

$$<\tilde \sigma^{1},\eta>_{L^{2}(h)}=0\quad \forall \eta\in W^{1}_{n}+ Vect(\eta^{p-1}).$$
\begin{prop}\label{sig1} Nous définissons une section holomorphe, projection  orthogonale
 de $\sigma^{1}$  sur $W_{n}$ par :

\begin{equation}\label{minj}
  \tilde\sigma^{1}=\sigma^{1}-\sum_{i=1..p-1}<\sigma^{1},\eta^{i}>_{L^{2}(h)}\eta^{i},
\end{equation}

avec
$$\normlh{\tilde\sigma^{1}}^{2}=\frac{\pi}{n}-\frac{\pi
d_{1}+\pi r}{n^2}+O(\frac{1}{n^{3}}).$$
$$\tilde\sigma^{1}(z_{0})=e_{1}(z_{0})+O(\negli{\frac{3}{2}})e_{2}(z_{0}).$$

Ici $r$  ne dépend que de $X$(cf. formule
\ref{vol} p. \pageref{vol}).
\end{prop}
\begin{demo}
On a en effet :
\begin{itemize}
  \item $\forall \eta^{i}$,$i=1..p-1$,
  \begin{eqnarray*}
  <\tilde\sigma^{1},\eta^{i}>_{L^{2}(h)}&=&
  <\sigma^{1},\eta^{i}>_{L^{2}(h)}-
\sum _{j=1..p-1}<\sigma^{1},\eta^{j}>_{L^{2}(h)}<\eta^{j},\eta^{i}>=0
  \end{eqnarray*}

  \item Nous utilisons la formule \ref{pml} et
  la proposition \ref{sigmaa} :
  \begin{eqnarray*}
  \normlh{\tilde \sigma^{1}}^{2}&=&< \sigma^{1}, \sigma^{1}>_{L^{2}(h)}
  -\sum_{i=1..p-1}\vert<\sigma^{1},\eta^{i}>\vert^{2}\\
&=&\frac{\pi}{n}-\frac{\pi d_{1}+\pi
r}{n^{2}}+O(\frac{1}{n^{3}})+
(p-1)(O(\negli{2}))^{2}
  \end{eqnarray*}
  et le fait que $p\sim n$.
  Pour la dernière propriété de $ \tilde \sigma^{1}$ on utilise l'équation \ref{minj}.
\begin{eqnarray*}
  \tilde \sigma^{1}(z_{0})=\sigma^{1}(z_{0})-\sum_{i=1..p-1}O(\negli{2})\eta^{i}(z_{0}),
\end{eqnarray*}
et le fait que seule
$\eta^{p-1}(z_{0})=O(n^{\frac{1}{2}})$.

\end{itemize}
\end{demo}
\subsection{Normalisation de $\tilde\sigma^{1}$}\label{zdfg}
Nous normalisons $\tilde\sigma^{1}$. Cela définit
la section  $\hat
\sigma^{1}$.  Alors $\hat
\sigma^{1}$ est la projection orthogonale de
$u_{1}$ sur $H^{0}(X,E_{n})$ pour la métrique $h$.
De plus elle vérifie :
$$\normlh{\hat \sigma^{1}}^{2}=1+O(\frac{1}{n})\mbox{ (cf. équation \ref{mlkj})}.$$
On a donc la proposition suivante :
\begin{prop}$\hat \sigma^{1}$ est la projection
orthogonale de $u_{1}$ sur $H^{0}(X,E_{n})$ pour
la métrique $h$. Elle vérifie  :

 $$\hat  \sigma_{1}^{1}(z_{0})=\frac{n+d_{1}+r}{\pi}e_{1}(z_{0})+c_{1}e_{1}(z_{0})
$$
où $\hat\sigma^{1}_{1}$ est la première composante
de $  \hat\sigma^{1}$ dans le repère
$\{e_{1},e_{2}\}$ et $c_{1}=O(\negli{}).$ De plus
$$\hat  \sigma_{2}^{1}(z_{0})=O(\negli{})e_{2}(z_{0})$$
 où $\hat\sigma^{1}_{2}$ est la deuxième composante de
$\hat\sigma^{1}$ dans le repère $\{e_{1},e_{2}\}$.

\end{prop}
On pose $$\hat\sigma^{1}=\eta^{p}.$$
\subsection{''Défaut''  d'orthogonalité de
 $\check\sigma^{2}$ avec $Vect(\eta^{p})$}
 On obtient de la même façon :

\begin{equation}\label{xsa}
<\check\sigma^{2},\eta^{p}>_{L^{2}(h)}=O(\negli{2}).
\end{equation}

%%%%%%%%%%%%%%%%%%%%%%%%%%%%%%%%%%%%%%
 \subsection{Projection de $u_{2}$ sur  $W^{1}_{n}+Vect(\eta^{p})$}
 On effectue le même travail que dans le paragraphe \ref{kkhj}

\begin{prop}\label{sig2} Nous définissons une section holomorphe, projection  orthogonale
 de $u_{2}$  sur $W^{1}_{n}+Vect(\eta^{p})$ par :

\begin{equation}\label{minj1}
  \tilde\sigma^{2}=\sigma^{2}-\sum_{i=1..p-2,p}<\sigma^{2},\eta^{i}>_{L^{2}(h)}\eta^{i},
\end{equation}

avec
$$\normlh{\tilde\sigma^{2}}^{2}=\frac{\pi}{n}-\frac{\pi
d_{1}+\pi r}{n^2}+O(\frac{1}{n^{3}}).$$
$$\tilde\sigma^{2}(z_{0})=e_{2}(z_{0})+O(\negli{\frac{3}{2}})e_{1}(z_{0}).$$

Ici $r$  ne dépend que de $X$(cf. formule
\ref{vol} p. \pageref{vol}).
\end{prop}
\begin{demo}
On a en effet :
\begin{itemize}
  \item $\forall \eta^{i}$,  $i=1..p-2,p$ ,
  \begin{eqnarray*}
  <\tilde\sigma^{2},\eta^{i}>_{L^{2}(h)}&=&
  <\sigma^{1},\eta^{i}>_{L^{2}(h)}-
\sum _{j=1..p-2,p}<\sigma^{1},\eta^{j}>_{L^{2}(h)}<\eta^{j},\eta^{i}>=0
  \end{eqnarray*}

  \item Nous utilisons la formule \ref{pml} et
  la proposition \ref{sigmaa} :
  \begin{eqnarray*}
  \normlh{\tilde \sigma^{2}}^{2}&=&< \sigma^{2}, \sigma^{2}>_{L^{2}(h)}
  -\sum_{i=1..p-2,p}\vert<\sigma^{2},\eta^{i}>\vert^{2}\\
&=&\frac{\pi}{n}-\frac{\pi d_{2}+\pi
r}{n^{2}}+O(\frac{1}{n^{3}})+
(p-1)(O(\negli{2}))^{2}
  \end{eqnarray*}
  et le fait que $p\sim n$.
  Pour la dernière propriété de $ \tilde \sigma^{1}$ on utilise l'équation \ref{minj1}.
\begin{eqnarray*}
  \tilde \sigma^{1}(z_{0})=\sigma^{1}(z_{0})-\sum_{i=1..p-2,p}O(\negli{2})\eta^{i}(z_{0}),
\end{eqnarray*}
et le fait que seule $\eta^{p}(z_{0})
=O(n^{\frac{1}{2}})$.

\end{itemize}
\end{demo}
\subsection{Normalisation de $\tilde\sigma^{2}$}
 On fait de même qu' au paragraphe \ref{zdfg}. Nous
normalisons $\tilde\sigma^{2}$. Cela définit la
section  $\hat
\sigma^{2}$.  Alors $\hat
\sigma^{2}$ est la projection orthogonale de
$u_{2}$ sur $H^{0}(X,E_{n})$ pour la métrique $h$.
De plus elle vérifie :
$$\normlh{\hat \sigma^{2}}^{2}=1+O(\frac{1}{n}).$$
On a donc la proposition suivante :
\begin{prop}$\hat \sigma^{2}$ est la projection
orthogonale de $u_{2}$ sur $H^{0}(X,E_{n})$ pour
la métrique $h$. Elle vérifie  :

 $$\hat  \sigma_{1}^{2}(z_{0})=\frac{n+d_{2}+r}{\pi}e_{2}(z_{0})+c_{1}e_{1}(z_{0})
$$
où $\hat\sigma^{1}_{1}$ est la première composante
de $  \hat\sigma^{1}$ dans le repère
$\{e_{1},e_{2}\}$ et $c_{1}=O(\negli{}).$ De plus
$$\hat  \sigma_{2}^{2}(z_{0})=O(\negli{})e_{1}(z_{0})$$
 où $\hat\sigma^{2}_{2}$ est la deuxième composante de
$\hat\sigma^{1}$ dans le repère $\{e_{1},e_{2}\}$.
\end{prop}
\section{Calcul de $tr(\alpha P_{n})$}
Dans le repère local  $\{e_{1},e_{2}\}$ , on
obtient
:
\begin{equation}\label{courbu1}
\pr_{n}(z_{0},z_{0})=\begin{pmatrix}
\hat \sigma^{1}(z_{0})& \hat
\sigma^{2}(z_{0})\\
  \end{pmatrix}=\begin{pmatrix}
\frac{n+d_{1}+r}{\pi}+O(\negli{})&O(\negli{})\\
O(\negli{})&
\frac{n+d_{2}+r}{\pi}+O(\negli{})\\
  \end{pmatrix}.
\end{equation}
Nous obtenons  alors :
\begin{eqnarray}\label{equar}
\tr(\alpha P_{n})&=&\int_{z\in
X}tr(\alpha \pr_{n}(z,z)d(vol(z))
\nonumber\\
\tr(\alpha P_{n})&=&\frac{1}{-2i\pi}\int_{X}tr(\alpha R_{h}(E_{n}))
+\int_{X}tr(\alpha
\mbox{diag}(\frac{r}{\pi}))\omega+O(\negli{})\nonumber\\
&=&\frac{1}{-2i\pi}\int_{X}tr(\alpha
R_{h}(E_{n}))+O(\negli{}),
 \end{eqnarray}
 en utilisant la formule
\ref{courbb} p. \pageref{courbb}.  On en déduit
alors :
\label{fini}
\begin{thm}{\label{torsiuon}}
 Soient  deux m\'etriques
hermitiennes $h_1$ et $h_{0}$ sur $E_{n}$.
 On consid\`ere le chemin $h (u)=uh_{1}+(1-u)h_{0}$
avec $u
\in
\lbrack0, 1\rbrack$.

$${\du}\tau_{h(u)}(X,
E_n)
=O(\negli{})
.$$
\end{thm}
\begin{demo}
Comme d'après la proposition \ref{afin}
p.\pageref{afin},
$${\du}\tau_{h(u)}(X, E_n)=\frac{1}{2i\pi}
\int_{X}tr(\alpha R_{h(u)}(E_{n}))+
\tr(\alpha P_{n}),$$  l'équation \ref{equar}
 donne :
 $${\du}\tau_{h(u)}(X, E_n)=O(\negli{}).$$

\end{demo}
\chapter{Calcul de $I_{n}L_{n}$}
Nous rappelons que :
\begin{enumerate}
  \item $I_{n}$ est un morphisme de $Met(W_n)$ dans $Met(E_{n})$ et à $m\in Met(W_n)$, $I_{n}(m)$
correspond  la métrique induite par le morphisme
surjectif $\varphi$ défini dans \ref{phix} c'est à
dire le  morphisme $I_{n}$ est défini de la façon
suivante
:
\begin{dfn}(cf.
définition \ref {i})
\begin{eqnarray*}
I_{n} : Met (W_{n})&\lr& Met(E_{n})\\
 I_{n}(m)(e_{z},e_{z})&=& min\{m(v,v)/v\in W_{n}\mbox{ et }v(z)=e_{z}\},
\end{eqnarray*}
avec $e_{z}\in E_{z}$, $m\in Met(W_{n})$.
\end{dfn}
  \item

 $L_{n}\, :\,Met(E_{n})\lr Met(W_n)$,
associe à une métrique $h$ sur $E_{n}$ la métrique
$L^{2}$ sur $W_n$ définie  par $h $ (cf.
définition \ref {l}).
\end{enumerate}

 Soit
$z_{0}\in X$, nous voulons calculer pour une
m\'etrique $h\in Met(E_{n})$, la quantit\'e :
$$I_{n}L_{n}(h)(e_{i},e_{j})(z_{0})\, ,\quad i,j=1..2\, ,$$
avec $e_{i}\, ,i=1..2$  défini précédemment p.
\ref{reperp}. Pour cela nous utilisons le travail
pr\'ec\'edent sur la cr\'eation de sections
``concentr\'ees". On
 a vu, dans le chapitre précédent,  que
  $\tilde \sigma^{1}$ est une section
holomorphe de norme $L_{n}(h)$ minimum et égale
 à $e_{1}(z_{0})$ à $O(\frac{1}{n^{3}})$ près. On considère $ v_{1}$ une section holomorphe du fibré
  $E_{n}$ telle que en $z_{0}$
 elle soit égale à $\tilde \sigma_{1}(z_{0})$ c'est à dire à
 $e_{1}(z_{0}) $ à $O(\frac{1}{n^{3}})$ près. On a alors
\begin{eqnarray*}
(I_{n}L_{n}(h)(v_{1},v_{1}))_{z_{0}} &=&<\tilde
\sigma^{1},\tilde\sigma^{1}>_{L_{n}(h)}\\
&=&<\tilde
\sigma^{1},\tilde\sigma^{1}>_{L^{2}(h)}\\
&=&\frac{\pi}{n}-\frac{\pi d_{1}+\pi
r}{n^{2}}+O(\frac{1}{n^{3}})\quad \mbox{ cf. prop.
\ref{sig1}}.
\end{eqnarray*}

Soit  $v_{2}$, une section holomorphe égale en
$z_{0}$ à $e_{2}(z_{0})$ à $O(\frac{1}{n^{a}})$ ,$\quad
\forall a\in
\N^{*}$ près.
 On considère alors la section $\tilde
\sigma^{2}$ comme défini dans la proposition
\ref{sig2} et on a
:
$$(I_{n}L_{n}(h)(v_{2},v_{2}))_{z_{0}}=<\tilde
\sigma^{2},\tilde\sigma^{2}>_{L_{n}(h)}=\frac{\pi}{n}-\frac{\pi
d_{2}+\pi r}{n^{2}}+O(\frac{1}{n^{3}})\quad \mbox{
cf. prop. \ref{sig2}}.$$ On a également
:
$$(I_{n}L_{n}(h)(v_{2},v_{1}))_{z_{0}}=<\tilde
\sigma^{2},\tilde\sigma^{1}>_{L_{n}(h)}=O(\negli{3})$$
et
$$I_{n}L_{n}(h)(v_{2},v_{1})_{z_{0}}=O(\negli{3}).$$
Donc finalement :
\begin{equation}\label{iln}
I_{n}L_{n}(h)_{z_{0}}=(\begin{pmatrix}
\frac{\pi}{n}-\frac{\pi
d_{1}+\pi r}{n^{2}}&0\\ 0&\frac{\pi}{n}-\frac{\pi
d_{2}+\pi r}{n^{2}}\\
\end{pmatrix}+\begin{pmatrix}
O(\negli{3})&O(\negli{3})\\
O(\negli{3})&O(\negli{3})\\
\end{pmatrix}).
\end{equation}
\label{end}
\begin{remq}
Comme $X$ est compact Les constantes
${\bf{O}(\negli{3})}$ sont indépendantes  du point
$z_{0}\in X$ choisi.
\end{remq}
Remarquons que  $\{v_{1},v_{2}\}$ est
$h_{z_{0}}$-orthonormée à $O(\negli{3})$. On
obtient donc, en utilisant l'équation \ref{ssss}
p.
\pageref{ssss}
:

$$h^{-1}I_{n}L_{n}(h)=\frac{\pi}{n^{2}}(R_{\C}(\ox(n)Id_{E_{0}}
-R_{\C,h}(E_{0})Id_{\ox(n)}+r Id_{E_{n}})+{\bf{O}(\negli{3})}.$$ On a
posé
$${\bf{O}(\negli{3})}:=\begin{pmatrix}
O(\negli{3})&O(\negli{3})\\
O(\negli{3})&O(\negli{3}))\\
\end{pmatrix}.$$Il en découle le
théorème suivant
\begin{thm}Soit $(E_{n},h)$ un fibré hermitien sur $X$, on a :
$$h^{-1}I_{n}L_{n}(h)=
\frac{\pi}{n^{2}}(R_{\C}(\ox(n)Id_{E_{0}}
-R_{\C,h}(E_{0})Id_{\ox(n)}+r Id_{E_{n}})+{\bf{O}(\negli{3})}.$$
Ici $r$ est une fonction ne dépendant que de
$X$(cf. formule \ref{vol} p. \pageref{vol}).
\end{thm}

Ce théorème  est utilisé dans la démonstration du
théorème \ref{princ} donnée p.
\pageref {dprinc}.

%%%%%%%%%%%%%%%%%%%%%%%%%%%%%%%%%%%%%%%%%%%%%%%%%%%%%%%%%%%%%%%%%%%%%%%%
%%%%%%%%%%%%%%%%%%%%%%%%%%%%%%%%%%%%%%%%%%%%%%%%%%%%%%%%%%%%%%%%%%%%%%%
%%%%%%%%%%%%%%%%%%%%%%%%%%%%%%%%%%%%%%%%%%%%%%%%%%%%%%%%%%%%%%%%%%%%%%%%
%%%%%%%%%%%%%%%%%%%%%%%%%%%%%%%%%%%%%%%%%%%%%%%%%%%%%%%%%%%%%%%%%%%%%
%%%%%%%%%%%%%%%%%%%%%%%%%%%%%%%%%%%%%%%%%%%%%%%%%%%%%%%%%%%%%%%%%%%%%%
%%%%%%%%%%%%%%%%%%%%%%%%%%%%%%%%%%%%%%%%%%%%%%%%%%%%%%%%%%%%%%%%%%%%%
%%%%%%%%%%%%%%%%%%%%%%%%%%%%%%%%%%%%%%%%%%%%%%%%%%%%%%%%%%%%%%%%%%%%%%
%%%%%%%%%%%%%%%%%%%%%%%%%%%%%%%%%%%%%%%%%%%%%%%%%%%%%%%%%thbib8.tex}

%%%%%%%%%%%%%%%%%%%%%%%%%%%%%%%bibliographie%%%%%%%%%%%%%%%%%%%%%%%%%%%

%%%%%%%%%%%%%%%%%%%%%%%%%%%%%%%%%%%%%%%%%%%%%%%%%%%%%%%%%%%%%%%%%%%%%
%%%%%%%%%%%%%%%%%%%%%%%%%%%%%%%%%%%%%%%%%%%%%%%%%%%%%%%%%%%%%%%%%%%%%%
\addcontentsline{toc}{chapter}{Index}
\printindex

\end{document}